%% file: popmdrv3.tex
\documentclass[12pt]{article}
\usepackage{amssymb}
\usepackage{amsmath}
\usepackage{graphics}

\usepackage{url}

\usepackage{hyperref}

\usepackage[numbers]{natbib}
\begin{document}

\newenvironment {proof}{{\noindent\bf Proof.}}{\hfill $\Box$ \medskip}

\newtheorem{theorem}{Theorem}[section]
\newtheorem{lemma}[theorem]{Lemma}
\newtheorem{condition}[theorem]{Condition}
\newtheorem{proposition}[theorem]{Proposition}
\newtheorem{remark}[theorem]{Remark}
\newtheorem{hypothesis}[theorem]{Hypothesis}
\newtheorem{corollary}[theorem]{Corollary}
\newtheorem{example}[theorem]{Example}
\newtheorem{definition}[theorem]{Definition}

\renewcommand {\theequation}{\arabic{section}.\arabic{equation}}
\def \non{{\nonumber}}
\def \hat{\widehat}
\def \tilde{\widetilde}
\def \bar{\overline}

\title{\large{\bf\ Genealogical constructions of population 
models}}
                                                       
\author{ \begin{small}
\begin{tabular}{ll}                              
Alison M. Etheridge \thanks{Research supported in part by EPSRC grants 
EP/I01361X and EP/K034316/1}
 & Thomas G. Kurtz \thanks{Research supported in part by NSF grant
DMS 11-06424}\\   
Department of Statistics & Departments of Mathematics and Statistics\\       
Oxford University & University of Wisconsin - Madison \\                   
24-29 St Giles & 480 Lincoln Drive\\                                                         
Oxford OX1 3LB & Madison, WI  53706-1388\\
UK & USA \\                        
etheridg@stats.ox.ac.uk & kurtz@math.wisc.edu     \\
\url{http://www.stats.ox.ac.uk/~etheridg/} & \url{http://www.math.wisc.edu/~kurtz/}  \\                      
\end{tabular}
\end{small}}

\date{March 2, 2018}

\maketitle

\begin{abstract}
Representations of population models in terms of 
countable systems of particles are constructed, in which 
each particle has a `type', typically recording both spatial 
position and genetic type, and a level.  For finite 
intensity models, the levels are distributed on $[0,\lambda ]$, 
whereas in the infinite intensity limit $\lambda\rightarrow\infty$, at each time $t$, 
the joint distribution of types and levels is 
conditionally Poisson, with mean measure $\Xi (t)\times \ell$ where $
\ell$ 
denotes Lebesgue measure and $\Xi (t)$ is a measure-valued 
population process.  The time-evolution of the levels captures 
the genealogies of the particles in the population.

Key forces of ecology and genetics 
can be captured within this common framework.  Models 
covered incorporate both individual and event based 
births and deaths, one-for-one replacement, immigration, 
independent `thinning' and independent or exchangeable 
spatial motion and mutation of individuals.  Since birth 
and death probabilities can depend on type, they also 
include natural selection.  The primary goal of the paper
 is to present particle-with-level or lookdown constructions 
 for each of these elements of a population model.  
 Then the elements can be combined to specify the desired model.  
 In particular, a non-trivial extension of the
spatial $\Lambda$-Fleming-Viot process is constructed.

\vspace{.1in}


\vspace{.1in}

\noindent {\bf MSC 20}10 {\bf Subject Classification:}  Primary:  
60J25, 92D10, 92D15, 92D25 92D40  Secondary:   60F05, 
60G09, 60G55, 60G57, 60H15, 60J68
 
\end{abstract}

\setcounter{equation}{0}
\input INTRO2

\setcounter{equation}{0}
\input COMPON2

\setcounter{equation}{0}
\input LDEXAMP1

\input LDEXAMP2

\input POPAPP2

\bibliography{eth}

\end{document}

%% file: INTRO2.tex


\section{Introduction}
\label{introduction}

\subsection{Background}
\label{background}

There is now a vast mathematical literature devoted to 
modeling the dynamics of biological populations.  The 
models employed generally fall into one of two classes:  
ecological models, that aim to elucidate the interactions 
within and between populations, and between those 
populations and the environment; and models of 
population genetics, that aim to explain the patterns of 
genetic variation observed in samples from a population.  
Ecological models typically take into account (some of) 
spatial structure, competition for resources, 
predator-prey interactions and changing environmental 
conditions.  Often they assume infinite populations, 
allowing one to concentrate on fluctuations in growth 
rates and ignore demographic stochasticity.  Models from 
population genetics, by contrast, often concentrate on 
the demographic stochasticity (known in that context as 
random genetic drift) which arises from the randomness 
due to reproduction in a finite population and assume 
that the population from which one is sampling is 
{\em panmictic\/} (that is, there are no group structures or 
mating restrictions) and of constant size.  The `size' 
however, is not taken to be the census population size, 
but rather an {\em effective\/} population size, which is 
intended to capture the effects of things like varying 
population size and 
spatial structure.  In particular, the underlying ecology  
is supposed to be encapsulated in this single parameter.  
This strategy has been surprisingly effective, but in 
most situations, notably when the population is 
geographically dispersed, the influence of different 
evolutionary and ecological forces on the value of the 
effective population size remains unresolved.  To address 
these effects one must combine ecological and genetical models.  

Whereas in ecological models one usually asks about the 
existence of equilibria or the probability that a species 
can invade new territory, in population genetics, data on 
the differences between genes sampled from a finite 
number of individuals in the population is used to infer 
the `genealogical trees' that relate those genes, and so 
from a practical point of view, it is the distribution of
 these trees that one would like to describe.  
As a result, 
we require a framework for modeling populations which 
allows one to combine ecology and genetics in such a 
way that the genealogical trees relating individuals in a 
sample from the population are retained.  Our goal in 
this paper is to provide just such a framework.  

Mathematical population genetics is concerned with 
models that capture, for large populations, the key 
forces of evolution that are acting on the population, but 
which are robust to changes in the fine detail of local 
reproduction mechanisms.  Diffusion limits lie at the 
heart of the theory.  The prototypical example is the 
Wright-Fisher diffusion which arises as an approximation 
to the dynamics of allele frequencies in large panmictic 
populations of neutral genes whose dynamics can be 
governed by a plethora of different models.  In this 
situation, the genealogical trees relating individuals in a 
sample are approximated by Kingman's coalescent, in 
which each pair of ancestral lineages coalesces into a 
common ancestor at a rate inversely proportional to the 
effective population size.  Na\"ively one obtains the 
Kingman coalescent as a `moment dual' to the diffusion.  
However, this is not sufficient to guarantee that it 
really approximates the genealogy of a sample from one 
of the individual based models.  Indeed, there are 
examples of systems of individual based models for 
which the allele frequencies are approximated by a 
common diffusion, but for which the genealogical trees 
relating individuals in a sample from the limiting 
populations have {\em different\/} distributions 
\citep{taylor:2009}.  Whereas the structure of the 
genealogical trees is usually implicit in the description 
of individual based models, in the diffusion limit the 
individuals have disappeared and with them their 
genealogies.  Our approach allows us to retain 
information about the genealogies as we pass to the 
limit.  

The framework that we shall present here is very 
general.  It will allow us to construct population models 
that capture the key ecological forces shaping the 
population as well as demographic stochasticity.  Many 
`classical' examples will emerge as special cases.  We 
shall use it to pass from individual based models to 
continuous approximations, but while retaining 
information about the way in which individuals in a 
random sample from the population are related to one 
another.  In particular, we shall fulfill one of our 
primary aims when we began this project, by 
constructing the spatial $\Lambda$-Fleming-Viot process (that 
was introduced in \cite{BEV10,Eth08}) as a high-density 
limit of a class of individual based models that 
generalize those considered by \cite{BEH09} (Section 
\ref{sectslfv}).  We also 
present a different construction, equivalent in the 
high-density limit to that of \cite{VW15}, but requiring 
somewhat weaker conditions. Moreover, we present a generalisation
of the spatial $\Lambda$-Fleming-Viot process which incorporates
fluctuations of the local population density (Section 
\ref{sectslfv2}).

\subsection{Approach}
\label{approach}

Our approach belongs to the family of `lookdown 
constructions'.  Building on the ideas of \cite{DK96} and 
\cite{DK99a}, a number of authors have developed 
constructions of population models that incorporate 
information about genealogical relationships. 
These constructions typically involve assigning each 
individual in the population to a (non-negative) integer or 
real-valued `level', with connections between the levels 
determining the genealogical trees.  They are generically 
referred to as `lookdown' constructions since, in most 
cases, during reproduction events, offspring inserted at a 
given level `look down' to individuals at lower levels to 
determine their parent.  

Lookdown constructions are simplest if the spatial 
locations or types of individuals in the population do not 
affect the reproductive dynamics.  In that setting, the 
`levels' can be taken to be nonnegative integer-valued. 
The processes are constructed in such a way that at 
each time $t$, the types, elements of an appropriate space 
$E$, 
of the individuals indexed by 
their levels $\{X_i(t)\}$ are exchangeable, that is, the joint 
distribution does not change if we permute the indices, 
and in an infinite 
population limit, the measure that gives the state of the 
limiting measure-valued process is simply the de Finetti 
measure of the infinite exchangeable family $\{X_i(t)\}$.

We illustrate the key 
idea for the simple example originally considered in 
\cite{DK96}.  Consider a population of constant size $N$.  
Individuals are assigned levels $1,\ldots ,N$ by choosing 
uniformly at random among all possible assignments.  
The dynamics are as follows:  we attach an independent 
Poisson process $\pi_{(i,j)}$, of rate $\lambda$, to each pair $(
i,j)$ of 
levels.  At a point of $\pi_{(i,j)}$, the individual with the 
higher of the two levels $i$ and $j$ dies and is replaced by 
a copy of the individual with the lower level.  In 
between these replacement events, individuals 
(independently) accumulate mutations.  Since the level of 
an individual has such a strong impact on its evolution, 
it is not at all obvious that this description gives rise 
to a sensible population model.  To see that it does, one 
must show that if $\{X_i(0)\}$ is exchangeable, then for each 
$t>0$, $\{X_i(t)\}$ is exchangeable, and that the 
probability-measure-valued process $Z_N$ given by the
empirical measure $Z_N(t)=\sum_{i=1}^NX_i(t)/N$ has the same 
distribution as the probability-measure-valued process 
$\hat {Z}_N$ obtained from a sensible population model.  

Ignoring the possibility of mutations, the 
generator of the process described above is 
\[A^Nf(x)=\sum_{1\leq i<j\leq N}\lambda (f(\Phi_{ij}(x))-f(x)),\]
where $\Phi_{ij}(x)$ is obtained from $x$ by replacing $x_j$ by $
x_i$.
A sensible population model, specifically, a simple Moran 
model, has generator 
\[\hat {A}^Nf(x)=\frac 12\sum_{1\leq i\neq j\leq N}\lambda (f(\Phi_{
ij}(x))-f(x)).\]
In \cite{DK96},  it was shown that if $X^N$ is a solution of 
the martingale problem for $A^N$ and $\hat {X}^N$ is a solution of the 
martingale problem for $\hat {A}^N$ such that $X^N(0)$ and $\hat {
X}^N(0)$ 
have the same exchangeable initial distribution, then $Z^N$ 
and $\hat {Z}^N$ have the same  distribution as ${\cal P}(E)$-valued 
processes.

The proof in \cite{DK96} is based on an explicit 
construction and a 
filtering argument.  This filtering argument, along with a 
similar argument used in \cite{KKM90} in a proof of 
Burke's theorem in queueing theory, motivated the 
development of the Markov mapping theorem in \cite{Kur98}, 
Theorem \ref{mf}\ in the Appendix of this paper, which is 
a fundamental tool in the present work.

To apply the Markov mapping theorem in the setting of 
\cite{DK96}, for $x\in E^N$, let $z_N\in {\cal P}(E)$  be given by 
$z_N=\frac 1N\sum_{i=1}^N\delta_{x_i}$.  For $f\in B(E^N)$, the bounded, measurable 
functions on $E^N$, 
define
\[\alpha f(z_N)=\frac 1{N!}\sum_{\sigma}f(x_{\sigma (1)},\ldots ,
x_{\sigma (N)}),\]
where the sum is over all permutations of $\{1,\ldots ,N\}$.  In 
other words, we {\em average out\/} over the (uniform) 
distribution of the assignment of individuals to levels.  
We then observe that for $f\in B(E^N)$
\[\alpha A^Nf(z_N)=\alpha\hat {A}^Nf(z_N)=\frac 12\sum_{1\leq i\neq 
j\leq N}\lambda (\alpha f(z_N+N^{-1}(\delta_{x_i}-\delta_{x_j}))-
\alpha f(z_N))\equiv C_N\alpha f(z_N)\]
for any choice of $x$ satisfying $z_N=\frac 1N\sum_{i=1}^N\delta_{
x_i}$.  Theorem 
\ref{mf}\ then implies that for any solution $\tilde {Z}$ of the 
martingale problem for $C_N$ there exist solutions $X^N$ and 
$\hat {X}^N$ of the martingale problems for $A$ and $\hat {A}$ respectively, 
such that 
$Z^N$ and $\hat {Z}^N$ have  the same distribution as $\tilde {Z}^
N$. 
In other words, our model is really just the classical 
Moran model, but augmented with a very particular 
labeling of the individuals in the population.  A nice 
property of this labeling, is that the model for a 
population of size $N$ is embedded in that for a population 
of size $M$ for any $M>N$, and so it is straightforward to 
identify what will happen in the limit as $N\rightarrow\infty$.

Finally, observe that for $f\in\cup_NB(E^N)$, we can define
\[Af(x)=\sum_{1\leq i<j}\lambda (f(\Phi_{ij}(x))-f(x)),\]
that is, if $f\in B(E^N)$, $Af(x)=A^Nf(x)$.

Let $\{X_i(0)\}$ be an infinite exchangeable sequence in $E$, and 
construct a process $X(t)=\{X_i(t)\}$ using independent 
Poisson processes $\pi (i,j)$ as above.  Then $\{X_1,\ldots ,X_N\}$ is a 
solution of the martingale problem for $A^N$ and hence $X$ 
is a solution of the martingale problem for $A$.  The limit 
$Z(t)$
of $Z^N(t)$ is the de Finetti measure for $\{X_i(t)\}$ and 
averaging implies that for each $f\in B(E^N)$,
\[\langle f,Z^{(N)}(t)\rangle -\langle f,Z^{(N)}(0)\rangle -\int_
0^t\langle Af,Z^{(N)}(s)\rangle ds,\]
where $Z^{(N)}(t)$ is the $N$-fold product measure of $Z(t)$, is a 
$\{{\cal F}^Z_t\}$-martingale.  That, in turn, implies $Z$ is a 
Fleming-Viot process. These observations give an explicit 
construction of a process given implicitly in \cite{DH82}.

From the construction of $X$, it is also a 
simple matter to see that the genealogical trees relating 
individuals in the population are governed by the 
Kingman coalescent, just as for the Moran model.  
In addition, the genealogy of a 
sample of size $n$, that is, the particles at the $n$ lowest 
levels, does not change as we increase the 
population size since, by construction, the processes at 
the $n$ lowest level are not affected by the processes at 
the higher levels.

In order to extend the lookdown construction to the 
setting in which the locations or types of individuals in 
the population affect their reproductive dynamics, 
\cite{K00} introduced the idea of taking random levels in 
$[0,\infty )$.  More precisely, writing $E$ for the space in which 
the population evolves, conditional on the empirical 
measure of the population configuration being $K(t)$ at 
time $t$, `individuals' are assigned types and levels 
according to a Poisson distribution on $E\times [0,\infty )$ with 
mean measure $K(t)\times\ell$, where $\ell$ is Lebesgue measure.  If 
we `average out' over the distribution of the levels we 
recover $K(t)$.  Under appropriate conditions, the most 
important of which is that the generator governing the 
dynamics of the labeled population respects the 
conditionally Poisson structure (the analogue of the 
exchangeability in the case of fixed levels), the Markov 
mapping theorem, Theorem \ref{mf}, allows us to 
conclude that by `removing the levels' we recover the 
Markov process whose generator is obtained through this 
process of averaging.  In particular, existence of a 
solution to the martingale problem for the unlabeled 
population process is enough to guarantee existence of a 
solution to the martingale problem for the labeled 
population, from which a solution to that for the 
unlabeled population can be read off by averaging.  
Moreover, uniqueness of the solution of the labeled 
martingale problem guarantees that of the solution to 
the unlabeled martingale problem.  In \cite{K00}, this 
approach was used to construct measure-valued 
population models with spatially dependent birth and 
death rates:  for a given spatial location, offspring can 
be inserted at rates that depend on the local 
configuration without destroying the conditionally 
Poisson structure.  Poisson levels have been used 
extensively since (e.g.~\cite{Buh02, DEFKZ00, GLW05, 
VW15}).  In \cite{KR11}, levels are again conditionally 
Poisson, but now they are allowed to evolve continuously 
with time, a device which we shall also exploit in this 
work.  The main novelty in the examples presented here 
is that we are able to (flexibly) incorporate `event-based' 
updating mechanisms in the lookdown construction.  

Our approach in this article will be to define population 
models in which individuals are assigned levels, to 
average out over those levels in order to identify the 
unlabeled population model, and to pass to an infinite 
population limit.  Justification of this approach to 
constructing the unlabeled population model is based 
upon filtering arguments, that is, the `averaging out' 
corresponds to conditioning on all information about the 
past of the process {\em except\/} the levels of the particles.
Ensuring the validity of this conditioning argument 
requires that 
the assignment of individuals to levels be done in such a 
way that past observations of the distribution of spatial 
positions and genetic types does not give any 
information about the current levels of individuals in the 
population.  It is important to realize that such 
assignments are far from unique.  For example, in 
\S\ref{pure death process} we provide three possible ways 
for levels to evolve in a simple pure death process 
and in \S\ref{sectslfv}, we give two different 
particle constructions 
of the spatial $\Lambda$-Fleming-Viot process.

In the models we consider, new individuals have a single 
parent.  This assumption is common in both genetic and 
ecological models.  To specify a model, one must specify 
the rules by which a parent is selected, the rules 
by which the number and types of offspring are 
determined, the rules that determine the time of death 
of an individual, and the rules by which types change 
through movement, mutation, or other process.  One can 
also include such processes as immigration.  The 
primary goal of the paper is to outline how one can 
obtain a lookdown construction for any such model, and 
hence determine the genealogy of a sample of 
individuals from the population.  In
\S\ref{components}, we consider each of the pieces 
separately.  One then constructs a model by selecting 
``one of these'' and ``one of those'' and ``one of something 
else.''  Since we are considering Markov models, each 
piece corresponds to a generator, and the final model is 
essentially obtained by adding the generators.  Since each 
piece has a lookdown representation built in, the 
lookdown representation of the final model is obtained.  
This description of the construction is formal and 
additional work must be done to ensure that the 
generator obtained uniquely determines a process. One 
useful approach to proving uniqueness is to show that 
the martingale problem is equivalent to a system of 
stochastic equations (c.f., Theorem \ref{sumeq}) and then 
prove uniqueness for the system of equations.  For 
example, see Lemma \ref{lamexun}, which gives a new 
proof of uniqueness for the spatial $\Lambda$-Fleming-Viot process under 
conditions given in \cite{VW15}.

\subsection{Structure of paper}
\label{structure}

The rest of the paper is laid out as follows.  In 
\S\ref{notation} we lay out the notation that we need for 
our discrete and continuous population models and for 
the `averaging' operations that we apply when we use 
the Markov mapping theorem.  In order to construct our 
general models, we exploit the fact that sums of 
generators are typically generators (see, for example, 
Problem 32 in Section 4.11 of \cite{EK86}), and so we can break 
our models apart into component pieces.  In 
\S\ref{components}, we examine each of these components 
in turn.  In \S\ref{examples}, we draw these together into 
a collection of familiar, and not so familiar, examples.  
For convenience, some useful identities for Poisson 
random measures are gathered together in 
Appendix~\ref{poisson identities}, and the Markov 
mapping theorem is stated in Appendix~\ref{mpsect}.  We 
refer to Appendix~A.2 of \cite{KR11} for necessary 
results on conditionally Poisson systems.  

We need to emphasize that although \S\ref{components} is 
the main focus of the paper, it contains 
calculations, not proofs.  These calculations give the 
first step in the application of the Markov mapping 
theorem, Theorem \ref{mf}, which ensures that the 
lookdown constructions actually represent the desired 
processes, but additional details must be checked for 
particular applications. 
We spell this out in the simplest example of a pure death 
process in \S\ref{pure death process}
and in the novel setting of the spatial $\Lambda$-Fleming-Viot process
and its extensions in \S\ref{sectslfv} and \S\ref{sectslfv2}.
In addition, the discrete particle 
models, indexed by $\lambda >0$, should converge to 
measure-valued models as $\lambda\rightarrow\infty$.  For many of the 
models, convergence of the lookdown constructions is 
obvious while in other cases, convergence follows easily 
by standard generator/martingale problem arguments.  It 
is then useful to know that convergence of the lookdown 
constructions implies convergence of the corresponding 
measure-valued processes.  Appendix~A.3 of 
\cite{KR11} provides the results needed to verify this 
convergence.

The results given in \S\ref{examples}\ are intended to be 
rigorous unless otherwise indicated.

\subsection{A note of appreciation} We are impressed by and thankful for the time and effort the Associate Editor and referees have put into this paper.  We hope that, in then end, they feel it was worthwhile.

\section{Notation}
\label{notation}

We will consider continuous-time, time-homogeneous, 
Markov models specified by their generators.  Each 
individual will have a {\em type\/} chosen from a complete 
separable metric space $(E,d)$.  We emphasize that here 
we are using `type' as shorthand for both spatial location 
and genetic type.  The distribution of types over $E$ may 
be discrete, that is, given by a counting measure that 
``counts'' the number of individuals in each subset of $ $$E$, 
or continuous, that is, the distribution of types is given 
by a measure on $E$ as in the classical examples of 
Dawson-Watanabe and Fleming-Viot. 
In addition, each individual will be 
assigned a `level' which in the discrete case will be 
sampled from an interval $[0,\lambda ]$ and in the continuous 
case from $[0,\infty )$.  No two individuals will have the same 
level, and in the continuous case, the types along with 
their levels give a countable collection of particles that 
determines the measure.

A state of one of our discrete population models will be 
of the form $\eta =\sum\delta_{(x,u)}$, where $(x,u)\in E\times [
0,\lambda ]$.  We shall 
abuse notation and treat $\eta$ both as a set and a counting 
measure, with the understanding that multiple points are 
treated as distinct individuals.  In other words 
\[\sum_{(x,u)\in\eta}g(x,u)=\int g(x,u)\eta (dx,du)\mbox{\rm \  and  }
\prod_{(x,u)\in\eta}g(x,u)=\exp\{\int\log g(x,u)\eta (dx,du)\}.\]
The projection of $\eta$ on $E$ will be denoted $\bar{\eta }=\sum_{
(x,u)\in\eta}\delta_x$ 
and $\eta$ will have the property that conditional on $\bar{\eta}$, the 
levels of the individuals in the population are 
independent uniform random variables on $[0,\lambda ]$.  It will 
be crucial that this conditioning property be preserved 
by the transformations of $\eta$ induced by the components 
in our generator.  Notice that allocating levels as 
independent uniform random variables is the natural 
continuous analogue of the way in which we allocated 
discrete levels through a uniform random sample from 
all possible permutations.  We shall write $\alpha (\bar{\eta },\cdot 
)$ for the 
joint distribution of independent uniform $[0,\lambda ]$ random 
variables $U_x$ indexed by the points $x\in\bar{\eta}$.  If $f$ is a 
function of the $U_x$, then $\alpha f$ will denote the 
corresponding expectation.  

When there is a need to be precise about the 
state space ${\cal N}_0$ for the counting measures $\eta$, 
we will  assume that ${\cal N}_0$ 
satisfies the following condition.
\begin{condition}\label{statecnd}
There exist $c_k\in C(E\times [0,\infty ))$ (or $c_k\in C(E\times 
[0,\lambda ])$ if 
$\lambda <\infty$), $k=1,2,\ldots$, $c_k\geq 0$, 
$\sum_{k=1}^{\infty}c_k(x,u)>0$, $(x,u)\in E\times [0,\infty )$ such that  $
	\eta\in {\cal N}_0$ if and 
only if
$\int_{E\times [0,\infty )}c_k(x,u)\eta (dx,du)<\infty$ for each $
	k$, and for $\eta_n,\eta\in {\cal N}_0$, 
$\eta_n\rightarrow\eta$ if and only if $\int_{E\times [0,\infty )}
fd\eta_n\rightarrow\int_{E\times [0,\infty )}fd\eta$ for each 
$f\in C(E\times [0,\infty ))$ such that $|f|\leq a_fc_k$ for some $
k$ and some 
$a_f\in (0,\infty )$.
\end{condition}
We note that ${\cal N}_0$ defined in this way will be a Polish 
space, and if all the $c_k$ have compact support, convergence 
is just vague convergence.

Under appropriate conditions (which we make explicit for 
the examples in \S\ref{examples}) we can pass from the 
discrete population models to an infinite density limit.  
The resulting continuous population models arise as 
limits of states $\eta_{\lambda}$ under assumptions that imply 
$\lambda^{-1}\eta_{\lambda}(\cdot ,[0,\lambda ])$ converges (at least in distribution) to a 
(possibly random) measure $\Xi$ on $E$.  This is the analogue 
of convergence of the empirical distribution in the 
simple case of a fixed number of discrete levels 
described in \S\ref{approach}.  Since we require that the 
levels in $\eta_{\lambda}$ be conditionally independent uniform random 
variables given $\bar{\eta}_{\lambda}$, it follows that $\eta_{\infty}$, the limit of the 
$\eta_{\lambda}$, will be a counting measure on $E\times [0,\infty 
)$ that is 
conditionally Poisson with Cox measure $\Xi\times\ell$, $\ell$ being 
Lebesgue measure.  That is, for example, 
\[E[e^{-\int_{E\times [0,\infty )}f(x,u)\eta (dx,du)}|\Xi ]=e^{-\int_
E\int_0^{\infty}(1-e^{-f(x,u)})du\Xi (dx)}.\]
To mirror our notation in the discrete setting, in the continuous case, 
$\alpha (\Xi ,\cdot )$ will denote the distribution of a conditionally
Poisson random measure $\eta$  on $E\times [0,\infty )$ with mean measure $\Xi (
dx)\times\ell$.  See Appendix \ref{poisson identities} and Appendices A.1, A.2, and A.3 of \cite{KR11}.

To describe the generators of our population models, we take the 
domain to consist of functions of the form 
\begin{equation}f(\eta )=\prod_{(x,u)\in\eta}g(x,u)=\exp\{\int\log 
g(x,u)\eta (dx,du)\},\label{domfrm}\end{equation}
where $g$ is 
continuous in $(x,u)$, differentiable in $u$, and $0\leq g\leq 1$.  
In order for the generator to be defined in specific examples, $g$ may, for
example, be required to satisfy additional regularity conditions, but
the key point is that the collection 
of $g$ employed will be large enough to ensure that the 
domain is separating. In what follows, we will 
frequently write expressions in which $f(\eta )$ is 
multiplied by one or more factors of the form 
$1/g(x,u)$. It should be understood that if $g(x,u)=0$, 
it simply cancels the corresponding factor in $f(\eta )$. Since 
linear combinations of martingales are martingales, we 
could, of course, extend the domain to include finite 
linear combinations of functions of the form 
(\ref{domfrm}).

In the discrete case, if a transformation moves the level of an individual above
$\lambda$, then the individual dies.  We therefore impose the condition $g(x,u)=1$ if 
$u\geq\lambda$. 
In this case $\alpha f(\bar{\eta })=\prod_{x\in\bar{\eta}}\bar{g}(x)$, where 
$\bar {g}(x)=\lambda^{-1}\int_0^{\lambda}g(x,u)du$.

In the continuous case, we assume that there exists some $u_ g$ 
such that $g(x,u)=1$ for $u\geq u_g$.  
Consequently,  $h(x)=\int_0^{\infty}(1-g(x,u))du$ is finite, and we have
\[\alpha f(\Xi )=e^{-\int_E\int_0^{\infty}(1-g(x,u))du\Xi (dx)}=e^{
-\int_Eh(x)\Xi (dx)}.\]

%% file: COMPON2.tex


\section{Components of our generators}
\label{components}

Having established our notation, we now turn to the building blocks of our population 
models.  By combining these components, we will be able 
to consider models which incorporate a wide
range of reproduction mechanisms.

\subsection{Pure death process}
\label{pure death process}

In this subsection we introduce a component which, 
when we average over levels, corresponds to each 
individual in the population, independently, dying at an 
instantaneous rate $d_0(x)\geq 0$ which may depend on its 
type, $x$.  We reiterate that $x$ encodes both spatial 
position and genetic type.  In particular, we do not 
require the population to be selectively neutral.  

We assume that the level of an individual of type
$x$ evolves according to the differential equation
$\dot {u}=d_0(x)u$.  The individual will be killed when its level first reaches $\lambda$.
Note that since the initial level $u(0)$ of an 
individual must be uniformly distributed on $[0,\lambda ]$, if 
nothing else affects the level, the 
lifetime of the individual (that is the time $\tau$ until the 
level hits $\lambda$) is exponentially distributed,
\[P\{\tau >t\}=P\{u(0)e^{d_0(x)t}<\lambda \}=P\{u(0)<\lambda e^{-
d_0(x)t}\}=e^{-d_0(x)t},\]
and conditional on $\{\tau >t\}=\{u(0)e^{d_0(x)t}<\lambda \}$, 
$u(0)e^{d_0(x)t}$ is uniformly distributed on $[0,\lambda ]$.  

The generator of this process is
\[A_{pd}f(\eta )=\int_{E\times [0,\lambda ]}f(\eta )d_0(x)u\frac {
\partial_ug(x,u)}{g(x,u)}\eta (dx,du).\]
Note that $g(x,u)$  in the denominator cancels the 
corresponding factor in $f(\eta )$.  Consequently,
\[\alpha A_{pd}f(\bar{\eta })=\alpha f(\bar{\eta })\int_E\frac 1{
\bar {g}(x)}\lambda^{-1}\int_0^{\lambda}d_0(x)u\partial_ug
(x,u)du\,\bar{\eta }(dx).\]
Observing that 
\[\lambda^{-1}\int_0^{\lambda}u\partial_ug(x,u)du=\left.\lambda^{
-1}u(g(x,u)-1)\right|_0^{\lambda}-\lambda^{-1}\int_0^{\lambda}(g(
x,u)-1)du=1-\bar {g}(x),\]
we see that
\begin{equation}\alpha A_{pd}f(\bar{\eta })=\alpha f(\bar{\eta })
\int_Ed_0(x)(\frac 1{\bar {g}(x)}-1)\bar{\eta }(dx),\label{pdgen}\end{equation}
so that in this case, the projected population model is indeed
just a pure death process in which the death rates may
depend on the types of the individuals.

The calculation above was purely formal. It is instructive to illustrate
the work required to apply Theorem~\ref{mf} in the context of this 
simple example. The key is that we must be able to check~(\ref{opest});
that is, we restrict the domain of $A_{pd}$ and exhibit a function
$\psi\geq 1$ for which, for each $f$ in this smaller domain, we can find 
a constant $c_f$ such that
$|A_{pd}f(\eta)|\leq c_f\psi(\eta)$.  
To this end, suppose 
$K_1\subset K_2\subset\cdots$ are subsets of $E$ such that $E=\bigcup_
kK_k$ (for 
example, if $E={\Bbb R}^d$, we might take $K_k=B_k(0)$).  Then let 
${\cal D}(A)$ be the collection of $f$ of the form 
$f(\eta )=\prod_{(x,u)}g(x,u)$ for $g(x,u)\in C_b(E\times [0,\lambda 
])$ satisfying 
$\partial_ug(x,u)\in C_b(E\times [0,\lambda ])$ and $g(x,u)=1$ for $
(x,u)\notin K_k\times [0,u_g]$ for 
some $k\in {\mathbb N}$ and $0\leq u_g\leq\lambda$.  We can take 
$\psi$ in Theorem~\ref{mf}\ to be of the form 
\[\psi (\eta )=\int_{E\times [0,\lambda ]}\sum_kd_0(x)\delta_k{\bf 1}_{
K_k}(x)\eta (dx,du) +1 \]
for some $\{\delta_k\}$ satisfying $\sum_k\delta_k\sup_{x\in K_k}
d_0(x)<\infty$. (The `$+1$' is just to guarantee that $\psi\geq 1$.) Then 
for $g(x,u)=1$ outside $K_k\times [0,u_g]$, we can take
$c_f=u_g\Vert\partial_ug\Vert\delta_k^{-1}$ in (\ref{opest}), where $
\Vert\cdot\Vert$ denotes the 
sup norm.
The function $\tilde{\psi}$ of Theorem~\ref{mf}, which is just
$\alpha\psi(\bar{\eta})$, takes the form
\[\tilde{\psi }(\bar{\eta })=\int_E\sum_kd_0(x)\delta_k{\bf 1}_{K_
k}(x)\bar{\eta }(dx) +1.\]
Then, by Theorem \ref{mf}, any solution of the 
martingale problem for $\alpha A_{pd}$ satisfying 
$E[\int_0^t\tilde{\psi }(\bar{\eta}_s)ds]<\infty$ (which will hold provided 
$E[\tilde{\psi }(\bar{\eta}_0)]<\infty$) can be obtained from a solution of the 
martingale problem for $A_{pd}$.

Other choices of the dynamics of the process with levels 
would have projected onto the same population model on 
averaging out the levels.  For example, we could equally 
have obtained (\ref{pdgen}) by starting with 
\[\tilde {A}_{pd}f(\eta )=\int_{E\times [0,\lambda ]}f(\eta )d_0(
x)(\frac 1{g(x,u)}-1)\eta (dx,du)\]
(the levels don't move; the particles just disappear) or
\[\hat{A}_{pd}f(\eta )
=\int_{E\times [0,\lambda ]}f(\eta )d_0(x)\left(\frac{\lambda^2}{2}-\frac {
u^2}2\right)\frac {\partial^2_ug(x,u)}{g(x,u)}\eta (dx,du)\]
for $g$ such that $\left.\partial_ug(x,u)\right|_{u=0}=0$
(the levels diffuse and absorption at $\lambda$ corresponds to 
death of the particle). Checking that $\alpha\hat {A}_{pd}f=\alpha 
A_{pd}$ 
given in (\ref{pdgen}) is an exercise in integration by 
parts.

For the continuous population limit, conditionally Poisson 
as described in Section \ref{approach}, 
it is immediate that
\[A_{pd}f(\eta )=\int_{E\times [0,\infty )}f(\eta )d_0(x)u\frac {
\partial_ug(x,u)}{g(x,u)}\eta (dx,du).\]
Recall that $g(x,u)=1$ for $u$ above some $u_g$.
Defining $h(x)=\int_0^{\infty}(1-g(x,u))du$, and using the identities
of Lemma~\ref{app1},
\begin{eqnarray}
\nonumber\alpha A_{pd}f(\Xi )&=&\alpha f(\Xi )\int_E\int_0^{\infty}
d_0(x)u\partial_ug(x,u)du\Xi (dx)\\
\label{pdgeninfinity}&=&\alpha f(\Xi )\int_Ed_0(x)h(x)\Xi (dx),\end{eqnarray}
where $\alpha f(\Xi )=e^{-\int_Eh(x)\Xi (dx)}$.  Define 
\begin{equation}\Xi_t(dx)=e^{-d_0(x)t}\Xi_0(dx),\label{detev}\end{equation}
and note that 
\[\frac d{dt}\alpha f(\Xi_t)=\alpha f(\Xi_t)\int_Ed_0(x)h(x)\Xi_t
(dx),\]
so $\alpha A_{pd}$ 
is the generator corresponding to the evolution of $\Xi$ 
given by (\ref{detev}).

\subsection{Multiple deaths}
\label{multdeath}

Whereas in the pure death process of the previous 
subsection, individuals are removed from the population 
one at a time, we now turn to a model that allows for 
multiple simultaneous deaths.  Moreover, in place of 
individual based death rates, deaths in the population 
will be driven by a series of `events' at which a 
specified number of deaths occur.  Since in the discrete 
setting, death occurs when the level of an individual 
crosses level $\lambda$, in order to have multiple simultaneous 
deaths, 
the levels must evolve through a series of jumps (c.f.~the thinning 
transformation in \S\ref{sectthin}).

We parametrize the multiple death events by points 
from some abstract space ${\Bbb U}_d$.  Corresponding to each 
$z\in {\Bbb U}_d$ is a pair $(k(z),d_1(\cdot ,z))$, where $k(z)$ is an integer 
and $d_1(\cdot ,z)$ is a nonnegative function on $E$, which allows 
us to weight each individual's relative probability of 
death during an event according to its type and spatial 
position.  We shall focus on the case in which events 
happen with intensity determined by a measure $\mu_d$ on 
${\Bbb U}_d$, but exactly the same approach applies if we demand 
that the events occur at discrete times.  

For a given pair $(k,d_1(\cdot))$,
let 
\[\tau (k,d_1,\eta )=\inf\{v:\eta \{(x,u):e^{vd_1(x)}u\geq\lambda 
\}\geq k\},\]
where the infimum of an empty set is infinite.  After 
the death event, the configuration becomes
\[\theta_{k,d_1}\eta\equiv \{(x,e^{\tau (k,d_1,\eta )d_1(x)}u):(x
,u)\in\eta\mbox{\rm \ and }e^{\tau (k,d_1,\eta )d_1(x)}u<\lambda 
\}.\]
Note that $k$ individuals will die if
$\eta (\{(x,u):d_1(x)>0\})\geq k$.  Otherwise, all individuals in 
$\{(x,u):d_1(x)>0\}$ are killed.  

Now assuming that $k$ and $d_1$ depend on $z\in {\Bbb U}_d$,
the generator for the model in which discrete 
death events occur with intensity $\mu_d(dz)$ then takes the form
\[A_{md}f(\eta )=\int_{{\Bbb U}_d}(\prod_{(x,u)\in\eta}g(x,ue^{\tau 
(k(z),d_1(\cdot ,z),\eta )d_1(x,z)})-f(\eta ))\mu_d(dz).\]
Since, conditional on $\bar{\eta}$, the levels of individuals in the 
population are independent uniformly distributed random variables on 
$[0,\lambda]$, $\tau_{x,z}$ given by $U_xe^{d_1(x,z)\tau_{x,z}}=\lambda$  
is exponential with parameter $d_1(x,z)$. The lack of memory property of
the exponential distribution guarantees that the levels of individuals 
in the population after the event are still uniformly distributed 
on $[0,\lambda ]$.  Moreover, since the $\tau_{x,z}$ are independent,
\begin{eqnarray*}
\alpha A_{md}f(\bar{\eta })=\int_{{\Bbb U}_d}{\bf 1}_{\{k(z)\leq 
|\bar{\eta }|\}}\sum_{S\subset\bar{\eta },|S|=k(z)}d(S,z)\alpha f
(\bar{\eta })(\frac 1{\prod_{x\in S}\bar {g}(x)}-1)\mu_d(dz)\\
+\int_{{\Bbb U}_d}{\bf 1}_{\{k(z)>|\bar{\eta }|\}}(1-\alpha f(\bar{
\eta }))\mu_d(dz),\end{eqnarray*}
where 
\[d(S,z)=P\{\max_{x\in S}\tau_{x,z}<\min_{x\in\bar{\eta}\backslash 
S}\tau_{x,z}\}.\]

Note that while all the points in $\eta$ will be distinct (no 
two points can have the same level), we have not ruled 
out the possibility that multiple points may have the 
same type.  Consequently, the same value of $x$ may 
appear multiple times in $S$, that is, we allow $\bar{\eta}$ and $
S$ to 
be multisets.  

As particular examples, if $k(z)=1$ and $S=\{x\}$, then
\[d(S,z)=\frac {d_1(x,z)}{\int_Ed_1(y,z)\bar{\eta }(dy)},\]
and if $d_1(x,z)=\zeta_z{\bf 1}_{C_z}(x)$ and $S\subset C_z$, 
$|S|=k(z)$ (so that $k(z)$ individuals will be chosen at random from the
region $C_z$ to die), then 
\[d(S,z)={{\bar{\eta }(C_z)}\choose {k(z)}}^{-1}.\]

Many interesting high density limits require a balance 
between birth and death events.  However, we close this 
subsection with a high density limit for the discrete 
death process above when there are no balancing births.  
Suppose that $\lambda^{-1}\bar{\eta}_{\lambda}\rightarrow\Xi$ as $
\lambda\rightarrow\infty$.  At an event of type 
$z\in {\Bbb U}_d$, 
\[P\{\tau (k(z),d_1(\cdot ,z),\eta_{\lambda})>\lambda^{-1}c\}=P\{
\eta_{\lambda}\{(x,u):e^{cd_1(x,z)/\lambda}u\geq\lambda \}<k(z)\}
.\]
Now 
since, conditional on $\bar{\eta}_{\lambda}$, the levels $u$ are independent 
uniform random variables on $[0,\lambda ]$, for a single $(x,u)\in
\eta_{\lambda}$, 
the probability that $u\geq\lambda e^{-cd_1(x,z)/\lambda}$ is $1-
e^{-cd_1(x,z)/\lambda}$ and 
the events $\{u\geq\lambda e^{-cd_1(x,z)/\lambda}\}$ are independent.  
Consequently, a Poisson 
approximation argument implies that $P\{\tau (k(z),d_1(\cdot ,z),
\eta_{\lambda})>\lambda^{-1}c\}$ 
converges to $P\{Z_c<k(z)\}$ where, conditional on $\Xi$, $Z_c$ is 
Poisson distributed with parameter $\int cd_1(x,z)\Xi (dx)$.

Consider the motion of a single level.  The jumps (of 
size $(e^{\tau (k(z),d_1(\cdot ,z),\eta_{\lambda})d_1(x,z)}-1)u$) that it experiences 
whenever a death event falls are independent (by lack of 
memory of the exponential distribution) and so if we 
speed up time by $\lambda$ and apply the law of large numbers, observing
that $\lambda E[\tau(k(z), d_1(\cdot, z), \eta_\lambda)]=
\int_0^\infty P[Z_c<k(z)]dc=k(z)/\int d_1(x,z)\Xi (dx)$,
we see that, in the limit as $\lambda\rightarrow\infty$, the motion of a single 
level converges to 
\[\dot {u}=\int_{{\Bbb U}_d}\frac {k(z)}{\theta (z,\Xi (t))}d_1(x
,z)\mu_d(dz)u,\]
where $\theta (z,\Xi )=\int d_1(x,z)\Xi (dx)$.  The limit of $\lambda 
A$ is
\[A^{\infty}_{md}f(\eta )=\int_{{\Bbb U}_d}f(\eta )\int\frac {k(z
)}{\theta (z,\Xi )}d_1(x,z)u\frac {\partial_u g(x,u)}{g(x,u)}\eta (
dx,du)\mu_d(dz).\]
Integrating the limiting form of the generator by parts, 
exactly as we did to obtain~(\ref{pdgeninfinity}), yields 
\[\alpha A_{md}^{\infty}f(\Xi )=\alpha f(\Xi )\int_{{\Bbb U}_d}\int_
E\frac {k(z)d_1(x,z)h(x)}{\theta (z,\Xi )}\Xi (dx)\mu_d(z).\]
Note that there is a time change relative to the 
generator~(\ref{pdgeninfinity}) even in the case when 
$k(z)\equiv 1$ and $d_1(x,z)\equiv d_1(x)$, since deaths are driven by 
`events' and not linked to individuals.  

\subsection{Discrete birth events} \label{sectdb}

We shall consider two different approaches to birth 
events.  Just as in the case of deaths, a fundamental 
distinction will be that in the approach outlined in this 
subsection, births will be based on events and particle 
levels will evolve in a series of jumps, whereas in the 
next subsection, births will be individual based and levels 
will evolve continuously, according to the solution of a 
differential equation.  To emphasize this point, we shall 
refer to discrete and continuous birth events.  

A discrete birth event involves the selection of a 
parent, the determination of the number of offspring, 
and the placement of the offspring.  Selection of the parent 
is controlled by a function $r$ with $r(x)\geq 0$ (the larger $r(x)$,
the more likely an individual of type $x$ is to be the 
parent); the number of offspring is specified by an 
integer $k$; and the placement of the offspring is 
determined by a transition function $q(x,dy)$ from $E$ to 
$E^k$.  In this discrete model, we can either assume that the parent 
is eliminated from the population or that it is identified with 
the offspring at level $v^{*}$ defined below (in which case it 
jumps according to $q(x,dy)$ as a result of the event).  

For a birth event to occur for a given triple $(r,k,q)$, we
must have $\int r(x)\bar{\eta }(dx)>0$, otherwise no individual is 
available to be the parent.  If there is a parent 
available, then
 $k$ points, $v_1,\ldots ,v_k$, are chosen 
independently and uniformly on $[0,\lambda ]$.  These will be the 
levels of the offspring of the event.  Let $v^{*}$ denote the 
minimum of the $k$ new levels.  For old points $(x,u)\in\eta$ 
with $u>v^{*}$ and $r(x)>0$, let $\tau_x$ be defined by 
$e^{-r(x)\tau_x}=\frac {\lambda -u}{\lambda -v^{*}}$ and for $(x,
u)\in\eta$ satisfying $u<v^{*}$ and 
$r(x)>0$, let $\tau_x$ be determined by $e^{-r(x)\tau_x}=\frac u{
v^{*}}$.  Note 
that conditioned on $u>v^{*}$, $\frac {\lambda -u}{\lambda -v^{*}}$ is uniformly distributed 
on $[0,1]$ and similarly, conditioned on $u<v^{*}$, $\frac u{v^{*}}$ is 
uniformly distributed on $[0,1]$, so in both cases, $\tau_x$ is 
exponentially distributed with parameter $r(x)$.  Take 
$(x^{*},u^{*})$ to be the point in $\eta$ with $\tau_{x^{*}}=\min_{
(x,u)\in\eta}\tau_x$.  
This point will be the parent.  We have 
\[P\{x^{*}=x'\}=\frac {r(x')}{\int r(x)\bar{\eta }(dx)},\quad x'\in
\bar{\eta }.\]
After the event, the configuration $\gamma_{k,r,q}\eta$ of levels and types in
the population is obtained by 
assigning types $(y_1,\ldots ,y_k)$ with joint distribution $q(x^{
*},dy)$ 
uniformly at random to the $k$ new levels and 
transforming the old levels so that
\begin{eqnarray*}
\gamma_{k,r,q}\eta =\{(x,\lambda -(\lambda -u)e^{r(x)\tau_{x^{*}}}
):(x,u)\in\eta ,\tau_x>\tau_{x^{*}},u>v^{*}\}\\
\cup \{(x,ue^{r(x)\tau_{x^{*}}}):(x,u)\in\eta ,\tau_x>\tau_x^{*},
u<v^{*}\}\\
\cup \{(y_i,v_i),i=1,\ldots ,k\}.\end{eqnarray*}
Notice that the parent has been removed from the 
population and that if $r(x)=0$, the point $(x,u)$ is 
unchanged.  

Since $x^{*}$ and $\tau_{x^{*}}$ are deterministic functions of $\eta$ and 
$v^{*}\equiv\wedge_{j=1}^kv_k$, 
for $(x,u)\in\eta$, $(x,u)\neq (x^{*},u^{*})$, that is an `old' individual which is
not the parent, we can write the new level 
as ${\cal J}^{\lambda}_r(x,u,\eta ,v^{*})$. 
Then
\[f(\gamma_{k,r,q}\eta )=\prod_{(x,u)\in\eta ,u\neq u^{*}}g(x,{\cal J}^{
\lambda}_r(x,u,\eta ,v^{*}))\prod g(y_i,v_i).\]

The crucial feature of this construction is captured by 
the following lemma.  
\begin{lemma}\label{uniflev} 
Conditional on $\{(y_i,v_i)\}$ and $\bar{\eta}$, $\{{\cal J}^{\lambda}_
r(x,u,\eta ,v^{*}):(x,u)\in\eta ,u\neq u^{*}\}$ 
are independent and uniformly distributed on $[0,\lambda ]$.  
\end{lemma}

\begin{proof}
Conditioned on $\bar{\eta}$ and the vector $v=(v_1,\ldots ,v_k)$, 
the levels $u$ are independent and uniformly distributed on $[0,\lambda]$. 
Conditioned further on $u<v^{*}$, $u$ is uniform on $[0,v^{*}]$, whereas
conditioned on $u>v^{*}$, $u$ is uniform on $[v^{*},\lambda ]$.

Now if $u<v^{*}$ and $u\neq u^{*}$, then, by definition,
$ue^{r(x)\tau_{x^*}}<v^*$, that is $u<v^*e^{-r(x)\tau_{x^*}}$ and, 
conditional on $\tau_{x^*}$ and this event, $u$ is uniform on 
$[0, v^*e^{-r(x)\tau_{x^*}}]$.

Similarly, conditioning on $u>v^*$ and $u\neq u^*$, knowing $\tau_{x^*}$,
$u$ is uniform on $[\lambda -(\lambda -v^*)e^{-r(x)\tau_{x^*}},\lambda]$.

Consequently, for $(x,u)\in\eta$, we compute 
\begin{eqnarray*}
&&E[g({\cal J}_r^{\lambda}(x,u,\eta ,v^{*}))|\bar{\eta },u\neq u^{
*},v^{*}]\\
&&\qquad =E[g(\lambda -(\lambda -u)e^{r(x)\tau_{x^{*}}}),u>v^{*}|
\bar{\eta },u\neq u^{*},v^{*}]\\
&&\qquad\qquad +E[g(ue^{r(x)\tau_{x^{*}}}),u<v^{*}|\bar{\eta },u\neq 
u^{*},v^{*}]\\
&&\qquad =E[g(\lambda -(\lambda -u)e^{r(x)\tau_{x^{*}}})|\bar{\eta }
,\lambda -(\lambda -u)e^{r(x)\tau_{x^{*}}}>v^{*},u>v^{*},v^{*}]\frac {
\lambda -v^{*}}{\lambda}\\
&&\qquad\qquad +E[g(ue^{r(x)\tau_{x^{*}}})|\bar{\eta },ue^{r(x)\tau_{
x^{*}}}<v^{*},u<v^{*},v^{*}]\frac {v^{*}}{\lambda}\\
&&\qquad =\frac 1{\lambda -v^{*}}\int_{v^{*}}^{\lambda}g(z)dz\frac {
\lambda -v^{*}}{\lambda}+\frac 1{v^{*}}\int_0^{v^{*}}g(z)dz\frac {
v^{*}}{\lambda}\\
&&\qquad =\frac 1{\lambda}\int_0^{\lambda}g(z)dz,\end{eqnarray*}
where the conditioning $\lambda-(\lambda -u)e^{r(x)\tau_{x^*}}>v^*$ in the 
first line of the second equality captures
$u\neq u^*$ and to pass from the third line to the fourth we partition over
$\tau_{x^*}$ (and perform a simple change of variable in the integral).
\end{proof}

By Lemma  \ref{uniflev},
\[E[f(\gamma_{k,r,q}\eta )|\bar{\eta }]=\sum_{x'\in\bar{\eta}}\frac {
r(x')}{\int r(x)\bar{\eta }(dx)}\alpha f(\bar{\eta })\frac 1{\bar {
g}(x')}\int\prod_{i=1}^k\bar {g}(y_i)q(x',dy).\]
Of course if $\int r(x)\bar{\eta }(dx)=0$, there is no parent and 
$\gamma_{k,r,q}\eta =\eta$.

This tells us how a configuration will be transformed by 
a single discrete birth event.  Now, just as in the 
previous subsection, we suppose that the events are 
parametrized by some abstract space, this time denoted 
by ${\Bbb U}_b$, equipped with a measure $\mu_{db}$ that determines the 
intensity of events.  The discrete birth generator will 
then be of the form 
\[A_{db}f(\eta )=\int_{{\Bbb U}_b}
{\bf 1}_{\{\int r(x,z)\bar{\eta }(dx)>0\}} (H_{k(z),r(\cdot ,z),q(\cdot ,z
,\cdot )}(g,\eta )-f(\eta ))\mu_{db}(dz),\]
where
\begin{eqnarray*}
H_{k(z),r(\cdot ,z),q(\cdot ,z,\cdot )}(g,\eta )=\lambda^{-k(z)}\int_{
[0,\lambda ]^{k(z)}}\prod_{(x,u)\in\eta ,u\neq u^{*}(\eta ,v^{*})}
g(x,{\cal J}^{\lambda}_{r(\cdot ,z)}(x,u,\eta ,v^{*}))\\
\int_{E^{k(z)}}\prod_{i=1}^{k(z)}g(y_i,v_i)q(x^{*}(\eta ,v^{*}),z
,dy)dv_1\ldots dv_{k(z)}.\end{eqnarray*}
Integrating out the levels gives
\begin{eqnarray*}
\alpha A_{db}f(\bar{\eta })&=&\int_{{\Bbb U}_{db}}{\bf 1}_{\{\int 
r(x,z)\bar{\eta }(dx)>0\}}\\
&&\left(\sum_{x'\in\bar{\eta}}\frac {r(x',z)}{\int r(x,z)\bar{
\eta }(dx)}\frac{\alpha f(\bar\eta )}{\bar {g}(x')}\big(\int_{E^{
k(z)}}\prod_{i=1}^{k(z)}\bar g(y_i)q(x',z,dy)\big)-\alpha f(\bar\eta 
)\right)\mu_{db}(dz).
\end{eqnarray*}

If we wish to pass to a high density limit, we must 
control the size and frequency of the jumps in the level 
of an individual, so that the level process converges as 
we increase $\lambda$.  To investigate the restriction that this 
will impose on the discrete birth events, we examine 
${\cal J}^{\lambda}_r(x,u,\eta ,v^{*})$ more closely.  Recall that for $
u>v^{*}$, $\tau_x$ is 
defined by $e^{-r(x)\tau_x}=(\lambda -u)/(\lambda -v^{*})$.  Evidently, we are 
only interested in the case when $\lambda^{-1}\bar{\eta}_{\lambda}$ converges to a 
non-trivial limit, and in changes in those levels that we 
actually `see' in our limiting model, that is, to levels 
that are of order one.  For such changes, $v^{*}$ will also be 
order one, and then it is easy to see that for $\lambda$ 
sufficiently large, we will have $u^{*}>v^{*}$ and, since 
$(\lambda -u^{*})/(\lambda -v^{*})\rightarrow 1$, $\tau_{x^{*}}\rightarrow 
0$ as $\lambda\rightarrow\infty$.  Now 
\begin{eqnarray}
{\cal J}_r^{\lambda}(x,u,\eta ,v^{*})&=&{\bf 1}_{\{u>v^{*}\}}(\lambda 
-(\lambda -u)e^{r(x)\tau_{x^{*}}})+{\bf 1}_{\{u<v^{*}\}}ue^{r(x)\tau_{
x^{*}}}\label{levtran}\\
&=&{\bf 1}_{\{u>v^{*}\}}(ue^{r(x)\tau_{x^{*}}}-\lambda (e^{r(x)\tau_{
x^{*}}}-1))+{\bf 1}_{\{u<v^{*}\}}ue^{r(x)\tau_{x^{*}}},\non\end{eqnarray}
and so it follows that for 
$u<v^{*}$, ${\cal J}_r^{\lambda}(x,u,\eta ,v^{*})=ue^{r(x)\tau_{x^{
*}}}\rightarrow u$. However, for $(x,u)\in\eta$ with $u>v^{*}$,
$u\neq u^{*}$, and $r(x)>0$, 
\begin{eqnarray}
{\cal J}_r^{\lambda}(x,u,\eta ,v^{*})&=&ue^{r(x)\tau_{x^{*}}}-\lambda 
(e^{r(x)\tau_{x^{*}}}-1)\non\\
&=&ue^{r(x)\tau_{x^{*}}}-\lambda\frac {(e^{r(x)\tau_{x^{*}}}-1)}{
(e^{r(x^{*})\tau_{x^{*}}}-1)}(e^{r(x^{*})\tau_{x^{*}}}-1)\non\\
&=&ue^{r(x)\tau_{x^{*}}}-\lambda\frac {(e^{r(x)\tau_{x^{*}}}-1)}{
(e^{r(x^{*})\tau_{x^{*}}}-1)}(\frac {\lambda -v^{*}}{\lambda -u^{
*}}-1)\non\\
&=&ue^{r(x)\tau_{x^{*}}}-\frac {\lambda}{\lambda -u^{*}}(u^{*}-v^{
*})\frac {(e^{r(x)\tau_{x^{*}}}-1)}{(e^{r(x^{*})\tau_{x^{*}}}-1)}\non\\
&\rightarrow&u-(u^{*}-v^{*})\frac {r(x)}{r(x^{*})}.\label{hlim}\end{eqnarray}
Thus if a level jumps, then that jump will be order one.
It 
is clear that, regardless of balancing death events, to 
have stable behavior of the levels as $\lambda\rightarrow\infty$, we must 
have $v^{*}<u$ and $r(x)>0$ only finitely often per unit 
time.
Since for a given $k$, the probability that $v^{*}$ will be less 
than $u$ is $1-(\frac {\lambda -u}{\lambda})^k$, we need
\begin{equation}\lim_{\lambda\rightarrow\infty}\int_{{\Bbb U}_b}(
1-(\frac {\lambda -u}{\lambda})^{k(z)}){\bf 1}_{\{r(x,z)>0\}}\mu_b^{\lambda}
(dz)=u\lim_{\lambda\rightarrow\infty}\frac 1{\lambda}\int_{{\Bbb U}_
b}k(z){\bf 1}_{\{r(x,z)>0\}}\mu_b^{\lambda}(dz)<\infty\label{brthmnt}\end{equation}
for each $x$.  If the limit were infinite for some $x$, then 
each individual of that type would instantaneously 
become a parent and be removed from the population.

\subsection{Continuous birth}
\label{continuous births}

As an alternative to the birth process described above, 
in which levels move by discrete jumps, in this 
subsection we consider a birth process in which births are 
based on individuals and levels 
move continuously.  Our aim is to obtain a construction 
of a pure birth process in which an individual of type $x$ 
gives birth to $k$ offspring at a rate $r(x)$.  For simplicity, 
we assume that offspring adopt the type of their parent.  
In the model with levels, an individual $(x,u)\in\eta$ gives 
birth to $k$ offspring at rate $r(x,u)=(k+1)(\lambda -u)^k\lambda^{
-k}r(x)$.  
The parent remains in the population and the offspring 
are assigned levels independently and uniformly 
distributed above the level of the parent.  Evidently this 
will result in an increase in the proportion of 
individuals with higher levels and so to preserve the 
conditionally uniform distribution of levels, we make 
them move downwards.  We shall do this by making them 
evolve according to a differential equation $\dot {u}=r(x)G_k^{\lambda}
(u)$, for 
an appropriate choice of the function $G_k^{\lambda}$.  

At first sight, there is something arbitrary about the 
choice of the dependence of branching rate on level.  It 
is, of course, essential that $\lambda^{-1}\int_0^{\lambda}r(x,u)
du=r(x)$, so that 
when we average out over its level, the expected 
branching rate of an individual of type $x$ is indeed $r(x)$.  
However, in principle, other choices of $r(x,u)$ with this 
property would work, provided we change the 
differential equation driving the levels.  This particular 
choice has the advantage that it makes calculation of the 
averaged generator, and hence identification of $G_k^{\lambda}$, very 
straightforward.  
 
The generator of the process with levels is of the form
\begin{eqnarray}
A_{cb,k}f(\eta )=f(\eta )\sum_{(x,u)\in\eta}r(x)\bigg[\frac {(k+1
)}{\lambda^k}\int_u^{\lambda}\cdots\int_u^{\lambda}(\prod_{i=1}^k
g(x,v_i)-1)dv_1\cdots dv_k\non\\
+G_k^{\lambda}(u)\frac {\partial_ug(x,u)}{g(x,u)}\bigg].\label{cbdef}
\end{eqnarray}

For brevity, for the rest of this subsection, we drop the 
subscript $k$ in the generator.  In order to calculate 
$\alpha A_{cb}$, for each $x\in\bar{\eta}$, write $\bar{\eta}_x$ for $
\bar{\eta}\backslash x$.  Then 
\begin{eqnarray*}
\alpha A_{cb}f(\eta )=\sum_{x\in\bar{\eta}}r(x)f(\bar{\eta}_x)\Bigg
[\frac 1{\lambda}\int_0^{\lambda}g(x,u)\frac {(k+1)}{\lambda^k}\int_
u^{\lambda}\cdots\int_u^{\lambda}\big(\prod_{i=1}^kg(x,v_i)-1\big
)dv_1\cdots dv_kdu\\
+\frac 1{\lambda}\int_0^{\lambda}G_k^{\lambda}(u)\partial_ug(x,u)
du\Bigg].\end{eqnarray*}
Now observe that
\[\frac {k+1}{\lambda^{k+1}}\int_0^{\lambda}g(x,u)\int_u^{\lambda}
\cdots\int_u^{\lambda}\prod_{i=1}^kg(x,v_i)dv_1\cdots dv_kdu=\left
(\frac 1{\lambda}\int_0^{\lambda}g(x,u)du\right)^{k+1}.\]
To see this, notice that on the right side we have 
the result of averaging over $k+1$ independent uniform 
levels, while on the left we have $(k+1)$ times the result 
of averaging over those levels if we specify that the 
first level is the smallest, and by symmetry any of the 
$k+1$ uniform variables is equally likely to be the 
smallest.  This deals with the first term of the 
averaged generator.  All that remains of the expression 
in square brackets is
\begin{equation}\frac 1{\lambda}\int_0^{\lambda}\left(G_k^{\lambda}
(u)\partial_ug(x,u)-\frac {(k+1)(\lambda -u)^k}{\lambda^k}g(x,u)\right
)du.\label{remainder}\end{equation}
Now we make a judicious choice of $G_k^{\lambda}$.  Suppose that
\begin{equation}
\label{gklu}
G_k^{\lambda}(u)=\lambda^{-k}(\lambda -u)^{k+1}-(\lambda -u).
\end{equation}
Then, noting that $G_k^{\lambda}(0)=G_k^{\lambda}(\lambda )=0$, and integrating by parts,
we see that~(\ref{remainder}) reduces to
$$-\frac{1}{\lambda}\int_0^\lambda g(x,u) du$$
and so we obtain
\[\alpha A_{cb}f(\bar{\eta })=\alpha f(\bar{\eta })\sum_{x\in\bar{
\eta}}r(x)\left[\bar g(x)^k-1\right],\]
which is the generator of a branching process, as required.

Now consider what happens as $\lambda\rightarrow\infty$.
Since $g(x,u)\equiv 1$ for $u>u_g$, 
$$\frac{1}{\lambda}\int_u^\lambda g(x,v)dv\rightarrow 1\quad\mbox{as }
\lambda\rightarrow\infty,$$
and so the first term on the right hand side of~(\ref{cbdef}) vanishes,
and observing that 
\[G_k^{\lambda}(u)=u\frac {(1-\frac u{\lambda})^{k+1}-(1-\frac u{
\lambda})}{\frac u{\lambda}}\rightarrow -ku,\]
we obtain
\[A_{cb}^{\infty}f(\eta )=-f(\eta )\sum_{(x,u)\in\eta}r(x)ku\frac {
\partial_ug(x,u)}{g(x,u)}.\]
Assuming $\lambda^{-1}\eta_{\lambda}(\cdot\times [0,\lambda ]
)\rightarrow\Xi$, we have 
$\alpha f(\Xi )=e^{-\int_Eh(x)\Xi (dx)}$ and, using~(\ref{pois1fac}),
\begin{eqnarray*}
\alpha A^{\infty}_{cb}f(\Xi )=-e^{-\int_Eh(x)\Xi (dx)}\int_Er(x
)k\int_0^{\infty}u\partial_ug(x,u)du\Xi (dx)\\
=e^{-\int_Eh(x)\Xi (dx)}\int_Er(x)k\int_0^{\infty}(g(x,u)-1)du\Xi 
(dx)\\
=-e^{-\int_Eh(x)\Xi (dx)}\int_Er(x)kh(x)\Xi (dx),\end{eqnarray*}
(where to perform the integration by parts we have used that
$\partial_ug(x,u)=\partial_u(g(x,u)-1)$)
which corresponds to the evolution of $\Xi$ given by
\[\Xi_t(dx)=e^{r(x)kt}\Xi_0(dx).\]

\subsection{One for one replacement} \label{sectoneforone}

So far, we have considered separately the births and 
deaths of individuals.  In some models, 
it is natural to think of offspring as replacing 
individuals in the population and, thereby, maintaining 
constant population size.  In this section we consider 
three different models of one-for-one replacement.  For 
$\lambda <\infty$ we shall suppose that the population size is finite.  
In the first model, we specify a number $k<|\eta |$ of 
individuals to be replaced.  Those individuals are then 
sampled uniformly at random from the population.  In 
the second model, the probability, $r(x)$, that an individual 
of type $x$ is replaced is specified.  
Either of these models can be modified in such a way that events affect
only a subset $C\subset E$ and this allows us to replace the requirement that
$\eta$ be finite by a local condition (for example $\eta(C)<\infty$). Both are
special cases of a third model in which there is a probability distribution
$p(S)$ over subsets $S\subset \eta$ 
that determines the subset to be replaced.  We require 
$p(S)$ to depend only on the types of the members of $S$ 
and not on their levels.  
By focusing on the case in which events happen with intensity determined
by a measure $\mu_{dr,3}$ on ${\Bbb U}_{dr,3}$ we can replace $p(S)$, a 
probability, by $r(S)$, a rate, giving the intensity for a replacement
event involving the individuals in $S$.

In all three cases, we can take the levels to be fixed.
The parent $(x^{*},u^{*})$ is taken 
to be the individual chosen to be replaced that has the 
lowest level.  We assume that the types of the new 
individuals are chosen independently with distribution 
given by a transition function $q(x^{*},dy)$, but we could 
allow dependence provided the new individuals are 
assigned to the chosen levels uniformly at random.  

For the first model, it is natural to take a generator of the form
\begin{eqnarray*}
A_{dr,1}f(\eta )=\int_{{\Bbb U}_{dr,1}}\binom{|\eta|}{k(z)}^{-1}\sum_{S\subset\eta ,|S|=k(z)}f(\eta 
)(\prod_{(x,u)\in S}\frac {\int g(y,u)q(x^{*}(S),z,dy)}{g(x,u)}-1
)\mu_{dr,1}(dz)\end{eqnarray*}
where $x^{*}(S)=x'$ if $(x',u')\in S$ and $u'=\min\{u:(x,u)\in S\}$.  
As usual, ${\Bbb U}_{dr,1}$ parametrizes the events and they occur with 
intensity $\mu_{dr,1}$.  The levels are fixed and the individuals 
chosen to be replaced `look down', just as in the simple 
example of \S\ref{approach}, to identify their parental 
type.  Averaging over levels yields 
\[\alpha A_{dr,1}f(\bar{\eta })=\int_{{\Bbb U}_{dr,1}}
\binom{|\bar{\eta}|}{k(z)}^{-1}\sum_{S\subset\bar{
\eta },|S|=k(z)}\alpha f(\bar{\eta })\frac 1{k(z)}\sum_{x'\in S}(
\prod_{x\in S}\frac {\int\bar {g}(y)q(x',z,dy)}{\bar {g}(x)}-1)\mu_{
dr,1}(dz).\]

~

In the second case, let $\xi_{x,u}$ be independent random 
variables with $P\{\xi_{x,u}=1\}=1-P\{\xi_{x,u}=0\}=r(x)$.  Then
 $(x^{*},u^{*})\in\eta$ is the parent if $u^{*}=\min\{u:\xi_{x,u}
=1\}$.  Let 
\begin{equation}\hat {g}(x,z,u)=\int_Eg(y,u)q(x,z,dy).\label{ghatdef}\end{equation}
Once again we can fix the levels, in which case the generator will take the form
\begin{equation}A_{dr,2}f(\eta )=\int_{{\Bbb U}_{dr,2}}\Bigg(E\big[\prod_{(x
,u)\in\eta}(\xi_{x,u}\hat {g}(x^{*},z,u)+(1-\xi_{x,u})g(x,u))\big]-f(
\eta )\Bigg)\mu_{dr,2}(dz),\label{dr2lam}\end{equation}
where the expectation is with respect to the $\xi_{x,u}$ and $x^{
*}$ 
is a function of $\eta$ and the $\{\xi_{x,u}\}$.  (More precisely, $x^{*}$ is 
a function of $\eta$ and the subset $S$ of the individuals for 
which $\xi_{x,u}=1$.)  Recall our assumption that there is a $u_g$ such that
$g(x,u)=1$ for all $u>u_g$. This property is inherited by $\hat{g}$.
Thus the factor in the product 
in~(\ref{dr2lam}) is $1$ if $u\geq u_g$, and so 
the expectation in the integral can be 
written as
\begin{eqnarray}
H(g,\hat {g},\eta ,z)&\equiv&E[\prod_{(x,u)\in\eta ,u\leq u_g}(\xi_{
x,u}\hat {g}(x^{*},z,u)+(1-\xi_{x,u})g(x,u))]\label{expn}\\
&=&\sum_{S\subset\eta|_{E\times [0,u_g]}}E[\prod_{(x,u)\in 
S}(\xi_{x,u}\hat {g}(x^{*}(S),z,u)\prod_{(x,u)\notin S, u\leq u_g}(1-\xi_{x,u}
)g(x,u)]\nonumber\\
&=&\sum_{S\subset\eta|_{E\times [0,u_g]}}\prod_{(x,u)\in 
S}(r(x,z)\hat {g}(x^{*}(S),z,u))\prod_{(x,u)\notin S,u\leq u_g}(1
-r(x,z))g(x,u)\non.
\end{eqnarray}
Partitioning on the lowest level particle, we see that this expression can 
also be
written as 
\begin{eqnarray}\label{partition}
&&\sum_{(x^{*},u^{*})\in\eta}r(x^{*},z)\hat {g}(x^{*},z,u^{*})\prod_{
(x,u)\in\eta ,u<u^{*}}(1-r(x,z))g(x,u)\\
&&\qquad\qquad \times\prod_{(x,u)\in\eta ,u>u^{*}}\Big(r(x,z)\hat {g}(x^{*},z,u)
+(1-r(x,z))g(x,u) \Big)\non.
\end{eqnarray}

It will be useful to  write $A_{dr,2}$  as a sum of two terms,
\begin{eqnarray}
\label{dr2gen}A_{dr,2}f(\eta )&=&f(\eta )\int_{{\Bbb U}_{dr,2}}\bigg
(\sum_{(x^{*},u^{*})\in\eta_{|E\times [0,u_g]}}\frac {r(x^{*},z)(
\hat {g}(x^{*},z,u^{*})-g(x^{*},u^{*}))}{g(x^{*},u^{*})}\\
&&\qquad\qquad\qquad\qquad\qquad\qquad\times\prod_{(x,u)\in\eta_{
|E\times [0,u_g]},u\neq u^{*}}(1-r(x,z))\non\\
&&\quad\qquad +\sum_{S\subset\eta_{|E\times [0,u_g]},|S|\geq 2}\Big
(\frac {\prod_{(x,u)\in S}\hat {g}(x^{*}(S),z,u)}{\prod_{(x,u)\in 
S}g(x,u)}-1\Big)\prod_{(x,u)\in S}r(x,z)\non\\
&&\qquad\qquad\qquad\qquad\qquad\qquad\times\prod_{(x,u)\in\eta_{
|E\times [0,u_g]}-S}(1-r(x,z))\bigg)\mu_{dr,2}(dz).\non\end{eqnarray}
We separate the first term, in which only one individual is replaced in 
the event, because it looks like the 
generator for simple, {\em almost} independent evolution of the 
particle types.  We exploit this observation in Section 
\ref{sectslfv}.

Since
\[\lambda^{-k}\int_0^{\lambda}g(u')\left(\int_{u'}^{\lambda}g(u)d
u\right)^{k-1}du'=\frac 1k\left(\lambda^{-1}\int_0^{\lambda}g(u)d
u\right)^k,\]
(c.f. the calculations in \S\ref{continuous births}) it follows 
from~(\ref{expn}) that 
\begin{eqnarray*}
\alpha A_{dr,2}f(\bar{\eta })=\int_{{\Bbb U}_{dr,2}}\bigg(\sum_{\bar {S}
\subset\bar{\eta}}\left(\prod_{x\in\bar {S}}r(x,z)\right)\left(\frac 
1{|\bar {S}|}\sum_{y\in\bar {S}}\bar g(y,z)^{|S|}\right)\left(\prod_{
x\in\bar{\eta }-\bar {S}}(1-r(x,z))\bar g(x)\right)\\-\alpha f(\bar{
\eta })\bigg)\mu_{dr,2}(dz),
\end{eqnarray*}
where $\bar {g}(y,z)=\lambda^{-1}\int_0^{\lambda}\hat {g}(y,z,u)d
u$. 

~

In the third case, in which we specify the rate at which subsets 
of individuals are replaced, we can write
\begin{equation}
\label{adr3}
A_{dr,3}f(\eta )= \int_{{\Bbb U}_{dr,3}} \sum_{S\subset\eta}
r(S,z)f(\eta 
)(\prod_{(x,u)\in S}\frac {\int g(y,u)q(x^{*}(S),z,dy)}{g(x,u)}-1
)\mu_{dr,3}(dz)
\end{equation}
and 
\[\alpha A_{dr,3}f(\bar{\eta })=\int_{{\Bbb U}_{dr,3}}\sum_{\bar {S}\subset
\bar{\eta}}r(\bar {S},z)\alpha f(\bar{\eta })\frac 1{|\bar {S}|}\sum_{
x'\in\bar {S}}(\prod_{x\in\bar {S}}\frac {\int\bar {g}(y)q(x',z,d
y)}{\bar {g}(x)}-1)\mu_{dr,3}(dz).\]

So far we have dealt with finite population models with 
one-for-one replacement.  We now turn our attention to 
infinite population limits.  For the first model, $A_{dr,1}$, 
there are two natural ways to pass to an infinite 
population limit.  In one, the rate at which birth events 
occur remains the same, but the size of the event (by which we mean
the number of individuals replaced) grows 
with $\lambda$, that is, 
\[\lambda^{-1}k_{\lambda}(z)\rightarrow\kappa (z)<|\Xi |=\lim_{\lambda
\rightarrow\infty}\lambda^{-1}|\eta_{\lambda}|.\]
Asymptotically, this model behaves in the same way as $A_{dr,2}$ in 
the special case in which $r(x,z)\equiv\kappa (z)/|\Xi |$ and so we don't
consider it here. 

The other possibility is for $k(z)$ to remain fixed, but for 
$\mu_{dr,1}$ to increase with $\lambda$, that is, to have replacement 
events occur at an increasingly rapid rate.  (For 
example, this is the approach when we pass from a 
Moran model to a Fleming-Viot process.)  

First we identify the appropriate scaling.  Assume that 
$\lambda^{-1}\eta_{\lambda}(t,\cdot )\Rightarrow\Xi (t,dx)$, where $
\Xi (t,E)<\infty$.  
(Of course, unless other factors are acting, $\Xi(t,E)$ is constant in 
time, but recall that we are thinking of our components as 
`building blocks' of population models.)
If a discrete birth 
event $z$ occurs at time $t$, then conditional on 
$\eta_{\lambda}(t,\cdot\times [0,\lambda ])$ and $z$, the number of individuals selected 
with levels below $a$, where $0<a<\lambda$, is binomial with 
parameters $k(z)$ and $\frac {\eta_{\lambda}(t,E\times [0,a])}{\eta_{
\lambda}(t,E\times [0,\lambda ])}=O(\lambda^{-1})$.  Since the 
probability of selecting two levels below $a$ is $O(\lambda^{-2})$, if 
we are to see any interaction between levels in the 
limiting model, we need to scale $\mu_{dr,1}$ by $\lambda^2$.  On the 
other hand, if we scale $\mu_{dr,1}$ by $\lambda^2$, the rate at which the 
individual at a fixed level is selected is of order $\lambda$.  
When this happens, unless it is one of the (finite rate) 
events in which more than one level below $a$ is selected, 
the individual at the selected level
will necessarily be the parent of the event and so 
will jump to a new position determined by the transition 
density $q$.  If the limiting model is to make sense, we 
must therefore rescale $q$ in such a way that in the limit, the 
motion of a fixed level will be well defined.  

To make this more precise, suppose that an event of 
type $z$ occurs at time $t$.  If an individual has level $u$, the 
probability that they are the parent of the event is 
\[\frac {{{\eta_{\lambda}(t,E\times (u,\lambda ])}\choose {k(z)-1}}}{{{
\eta_{\lambda}(t,E\times [0,\lambda ])}\choose {k(z)}}}\approx\frac {
k(z)}{\eta_{\lambda}(t,E\times [0,\lambda ])}.\]
Assume that $q$ depends on $\lambda$. Then the motion of a 
particle at level $u$ due to its being chosen as a parent is 
essentially (since we ignore the asymptotically negligible 
number of times when a particle with level  below $u$ is also 
chosen) 
Markov with generator
\[\tilde {B}_{\lambda}g(x)=\frac {\lambda^2}{\eta_{\lambda}(t,E\times 
[0,\lambda ])}\int_{{\Bbb U}_{dr,1}}k(z)(g(y)-g(x))q_{\lambda}(x,
z,dy)\mu_{dr,1}(dz).\]
We assume that the 
Markov process with generator
\[B_{\lambda}g(x)=\lambda\int_{{\Bbb U}_{dr,1}} k(z)(g(y)-g(x))q_{\lambda}(x,z,dy)\mu_{
dr,1}(dz)\]
converges in distribution to a Markov process with 
generator $B$.  Then (up to a time change which, in the limit, will be
$1/\Xi(E)$) this Markov process will describe 
the motion of a particle at a
 fixed level that results from it being 
selected as parent of a replacement event.  Note that 
this convergence implies that for each $\epsilon >0$, 
\begin{equation}\int k(z){\bf 1}_{\{d(y,x)>\epsilon \}}q_{\lambda}
(x,z,dy)\mu_{dr,1}(dz)=O(\lambda^{-1}).\label{jmpneg}\end{equation}

Similarly, we identify the interaction between distinct 
levels in the limiting process.  If there are individuals 
at levels $u_1<u_2$ and an event of type $z$ occurs at time 
$t$, then the probability that $u_1$, $u_2$ are the lowest two 
levels selected is 
\[\frac {{{\eta_{\lambda}(t,E\times (u_2,\lambda ])}\choose {k(z)
-2}}}{{{\eta_{\lambda}(t,E\times [0,\lambda ])}\choose {k(z)}}}\approx\frac {
k(z)(k(z)-1)}{\eta_{\lambda}(t,E\times [0,\lambda ])^2}=O(\lambda^{
-2}).\]
We chose our rescaling in such a way that events 
involving two levels below a fixed level $a$ will occur at 
a rate $O(1)$, and by (\ref{jmpneg}), after the event, asymptotically, both 
the parent and the offspring will have the 
type of the parent immediately before the event.  In this limit, we 
will never see events involving three or more levels below a fixed level $a$.

If the replacement process is the only 
process affecting the population, then 
\[|\Xi |=\Xi (t,E)=\lim_{\lambda\rightarrow\infty}\frac {\eta_{\lambda}
(t,E\times [0,\lambda ])}{\lambda}\]
is constant 
in time and (recalling that $g(x,u)=1$ for $u>u_g$)
the limiting model will have generator
\begin{eqnarray*}
A_{dr,1}^{\infty}f(\eta )&=&\int_{E\times [0,\infty )}\frac 1{|\Xi 
|}f(\eta )\frac {Bg(x,u)}{g(x,u)}\eta (dx,du)\\
&&\qquad +\int_{{\Bbb U}_{dr,1}}\frac {k(z)(k(z)-1)}{|\Xi |^2}\sum_{
(x_1,u_1),(x_2,u_2)\in\eta ,u_1<u_2}f(\eta )(\frac {g(x_1,u_2)}{g
(x_2,u_2)}-1)\mu_{dr,1}(dz)\\
&=&\int_{E\times [0,\infty )}\frac 1{|\Xi |}f(\eta )\frac {Bg(x,u
)}{g(x,u)}\eta (dx,du)\\
&&+\int_{{\Bbb U}_{dr,1}}\frac {k(z)(k(z)-1)}{|\Xi |^2}\\
&&\times\sum_{(x_1,u_1),(x_2,u_2)\in\eta}\frac {f(\eta )}{g(x_1,u_
1)g(x_2,u_2)}\left[{\bf 1}_{\{u_1<u_2\}}(g(x_1,u_2)g(x_1,u_1)-g(x_
2,u_2)g(x_1,u_1))\right]\mu_{dr,1}(dz)\end{eqnarray*}
Applying (\ref{pois1fac}) and (\ref{prws}), the averaged 
generator becomes 
\begin{eqnarray*}
\alpha A_{dr,1}^{\infty}f(\Xi )&=&e^{-\int_Eh(x)\Xi (dx)}\Big[-\int_
E\frac 1{|\Xi |}\int_EBh(x)\Xi (dx)\\
&&\qquad +\int_{{\Bbb U}_{dr,1}}\frac {k(z)(k(z)-1)}{2|\Xi |^2}\int_{
E\times E}(h(x_1)^2-h(x_1)h(x_2))\Xi (dx_1)\Xi (dx_2)\mu_{dr,1}(d
z)\Big].\end{eqnarray*}
The dependence of the first term on $\mu_{dr,1}$ is absorbed into 
our definition of $B$.  If $|\Xi |\equiv 1$ and $k(z)=2$ for all $
z$, then 
we recognize the generator of a Fleming-Viot diffusion. 
(See, for example, Section 1.11 of \cite{Eth00}.) 

It is elementary to identify the limit of our second 
model as $\lambda$ tends to infinity.  Since
$g(x,u)=1$ for $u>u_g$, 
the only changes that we `see' are those that affect 
$\eta^{u_g}=\sum_{(x,u)\in\eta ,u\leq u_g}\delta_{(x,u)}$ and these are determined by 
the generator when $\lambda =u_g$, so that 
\begin{equation}A^{\infty}_{dr,2}f(\eta )=\int_{{\Bbb U}_{dr,2}}\Bigg(H(
g,\hat {g},\eta ,z)-f(\eta )\Bigg)\mu_{dr,2}(dz),\label{dr2inf}\end{equation}
with $H$ given by~(\ref{expn}).
If $\eta$ is conditionally Poisson with Cox measure 
$\Xi (dx)du$, $\{\xi_{x,u,z}\}$ are independent with 
$P\{\xi_{x,u,z}=1\}=1-P\{\xi_{x,u,z}=0\}=r(x,z)$, 
and 
\[\eta_1=\sum_{(x,u)\in\eta}\xi_{x,u,z}\delta_{(x,u)},\quad\eta_2=
\sum_{(x,u)\in\eta} 
(1-\xi_{x,u,z})\delta_{(x,u)},\]
then $\eta_1$ and $\eta_2$ are conditionally independent given $\Xi$, $
\eta_1$ 
and $\eta_2$ are conditionally Poisson with Cox measures 
$r(x,z)\Xi (dx)du$ and $(1-r(x,z))\Xi (dx)du$ respectively and the 
cumulative distribution function of the level of the 
lowest particle to be replaced is $1-e^{-u\int r(x,z)\Xi(dx)}$.  
We now recall that the $x$ coordinates of the points in $\eta_1$, 
ordered according to the $u$ coordinates, are exchangeable 
with de Finetti measure 
\[\frac {r(x,z)\Xi (dx)}{\int_Er(y,z)\Xi (dy)},\]
and partition on the lowest level particle as in~(\ref{partition}).
Using~(\ref{pois1fac}), this yields
\begin{eqnarray*}
E[H(g,\hat {g},\eta ,z)|\Xi ]&=&e^{-\int_Eh(x)(1-r(x,z))\Xi (dx)}
\int_E\int_0^{\infty}r(x^{*},z)\hat {g}(x^{*},z,u)e^{-u\int_Er(x,
z)\Xi (dx)}\\
&&\qquad\qquad\qquad\qquad\qquad\qquad\qquad\times e^{-\int_u^{\infty}
(1-\hat {g}(x^{*},z,v))dv\int_Er(x,z)\Xi (dx)}\Xi (dx^{*})du\\
&=&e^{-\int_Eh(x)(1-r(x,z))\Xi (dx)}\int_E\int_0^{\infty}r(x^{*},
z)\hat {g}(x^{*},z,u)e^{-\int_0^u\hat {g}(x^{*},z,v)dv\int_Er(x,z
)\Xi (dx)}\\
&&\qquad\qquad\qquad\qquad\qquad\qquad\qquad\times e^{-\hat {h}(x^{
*},z)\int_Er(x,z)\Xi (dx)}\Xi (dx^{*})du\\
\\
&=&e^{-\int_Eh(x)(1-r(x,z))\Xi (dx)}\int_E\frac {r(x^{*},z)}{\int_
Er(x,z)\Xi (dx)}e^{-\hat {h}(x^{*},z)\int_Er(x,z)\Xi (dx)}\Xi (dx^{
*})\\
&\equiv&{\Bbb H}(h,\hat {h},\Xi ,z),\end{eqnarray*}
where the factor $\int r(x,z)\Xi(dx)$ in the density function of the lowest
level has canceled with the denominator in the de Finetti measure of $\eta_1$
on the right hand
side in the first line and to get from the second line to the third we 
integrated with respect
to $u$ and used that $\hat{g}=1$ for $u>u_g$.
Thus 
\[\alpha A^{\infty}_{dr,2}=\int_{{\Bbb U}_{dr,2}}({\Bbb H}(h,\hat {h},\Xi 
,z)-\alpha f(\Xi ))\mu_{dr,2}(dz).\]

Evidently, 
since $A_{dr,1}$ and $A_{dr,2}$ are special cases of $A_{dr,3}$ and their 
continuous density limits are quite different, we can't expect a 
general result for the continuous density limit of $A_{dr,3}$, but a 
large class of limits should retain the discrete model form
\[A^{\infty}_{dr,3}f(\eta )=\sum_{S\subset\eta}\int_{{\Bbb U}_{dr,3}}r(S
,z)f(\eta )(\prod_{(x,u)\in S}\frac {\int g(y,u)q(x^{*}(S),z,dy)}{
g(x,u)}-1)\mu_{dr,3}(dz),\]
provided there is a sufficiently large class of functions 
$g$ satisfying
\begin{equation}\sum_{S\subset\eta}\int_{{\Bbb U}_{dr,3}}r(S,z)\sum_{(x,
u)\in S}|\int (g(y,u)-g(x,u))q(x^{*}(S),z,dy)|\mu_{dr,3}(dz)<\infty 
\label{sumcnd}\end{equation}
with $0\leq g\leq 1$ and 
$g(x,u)\equiv 1$ for $u>u_g$.  In Section \ref{sectslfv}, we 
consider an example in which we can center 
$g(y,u)-g(x,u)$ in order to weaken the condition in 
(\ref{sumcnd}).  The form of the averaged generator is 
problem dependent, but convex combinations of $\alpha A^{\infty}_{
dr,1}$ 
and $\alpha A^{\infty}_{dr,2}$ can arise.

\subsection{Independent thinning} \label{sectthin}

Independent thinning will work in essentially the same 
way as the pure death process.  However, whereas in 
the pure death process the levels grew continuously, 
here we scale them up by a (type-dependent) factor at 
discrete times.  Levels which are above level $\lambda$ after 
this multiplication are removed.  The generator with 
finite $\lambda$ is then of the form 
\[A_{th}f(\eta )=\int_{{\Bbb U}_{th}}(\prod_{(x,u)\in\eta}g(x,u\rho 
(x,z))-f(\eta ))\mu_{th}(dz),\]
for some $\rho (x,z)\geq 1$.
Setting $\rho (x,z)=\frac 1{1-p(x,z)}$, we see that the probability that $
\rho (x,z)U_x>\lambda$,
for $U_x$ uniformly distributed on $[0,\lambda ]$, is 
$P\{U_x>\lambda /\rho (x,z)\}=p(x,z)$.
Recalling that $g(x,u)=1$ for $u\geq\lambda$ and 
integrating out the levels gives
\[\alpha A_{th}f(\bar{\eta })=\int_{{\Bbb U}_{th}}(\prod_{x\in\bar{
\eta}}((1-p(x,z))\bar {g}(x)+p(x,z))-\alpha f(\bar{\eta }))\mu_{t
h}(dz),\]
which says that when a thinning event of type $z$ occurs, 
individuals are independently eliminated with (type-dependent) probability 
$p(x,z)$.

In the continuous population limit, the form of $A_{th}$ 
remains unchanged, and the projected operator becomes 
\[\alpha A_{th}f(\Xi )=\int_{{\Bbb U}_{th}}(e^{-\int_E\frac 1{\rho 
(x,z)}h(x)\Xi (dx)}-\alpha f(\Xi ))\mu_{th}(dz),\]
where as usual $h(x)=\int_0^{\infty}(1-g(x,u))du$ and $\alpha f(\Xi 
)=e^{-\int_Eh(x)\Xi (dx)}$.

\subsection{Event based models}\label{sectevent}
Motivated by the model considered in \citet*{BEH09}, we 
combine independent thinning and discrete birth so that 
both transformations take place at the same time.  
Event times and types $(t,z)$ are determined by a Poisson 
random measure with mean measure $dt\mu_{th,db}(dz)$.  The 
value of $z$ determines the number of offspring $k(z)$, the 
relative chance $r(x,z)$ that an individual of type $x$  will 
be the parent, and the parameter $\rho (x,z)$ that determines 
the probability 
\[p(x,z)=\frac {\rho (x,z)-1}{\rho (x,z)}\]
that an individual of type $x$ is killed. Let
\[\bar{\eta }(r,z)=\int r(x,z)\bar{\eta }(dx),\]
and note that for there to be a parent, we must have 
$\bar{\eta }(r,z)>0$.  We will assume that the parent is killed, 
although alternatively, we could interpret the model as 
saying the parent jumps to the location of the particle 
at level $v^{*}$.

The 
form of the generator will be 
\[A^{\lambda}_{th,db}f(\eta )=\int_{{\Bbb U}}{\bf 1}_{\{\bar{\eta }
(r,z)>0\}}(H_z^{\lambda}(g,\eta )-f(\eta ))\mu_{th,db}(dz),\]
where, for ${\cal J}^{\lambda}_r$ given by (\ref{levtran}), if $\bar{
\eta }(r,z)>0$, 
\begin{eqnarray*}
H_z^{\lambda}(g,\eta )&=&\lambda^{-k(z)}\int_{[0,\lambda ]^{k(z)}}
\prod_{(x,u)\in\eta ,u\neq u^{*}(\eta ,v^{*})}g(x,\rho (x,z){\cal J}^{
\lambda}_{r(\cdot ,z)}(x,u,\eta ,v^{*}))\\
&&\qquad\qquad\times\prod_{i=1}^{k(z)}\int_Eg(y_i,v_i)q(x^{*}(\eta 
,v^{*}),z,dy_i)dv_1\ldots dv_{k(z)}.\end{eqnarray*}
The first product in the integral accounts for the 
thinning of the existing population (after the removal of 
the parent), and the second product accounts for the 
births.  Note that $(x^{*},u^{*})$ is a function of $\eta$ and $v^{
*}$, and 
if an event $z$ occurs at time $t$ and $\bar{\eta}_{t-}(r,z)>0$, then 
\begin{eqnarray*}
\eta_t=\sum_{(x,u)\in\eta_{t-},u\neq u^{*}}{\bf 1}_{\{\rho (x,z){\cal J}^{
\lambda}_{r(\cdot ,z)}(x,u,\eta_{t-},v^{*})<\lambda \}}\delta_{(x
,\rho (x,z){\cal J}^{\lambda}_{r(\cdot ,z)}(x,u,\eta_{t-},v^{*}))}
+\sum_{i=1}^{k(z)}\delta_{(y_i,v_i)}.\end{eqnarray*}

Averaging gives
\[\alpha A^{\lambda}_{th,db}f(\bar{\eta })=\int_{{\Bbb U}}{\bf 1}_{
\{\bar{\eta }(r,z)>0\}}\sum_{x^{*}\in\bar{\eta}}\frac {r(x^{*},z)}{
\int r(x,z)\bar{\eta }(dx)}(\bar {H}_z^{\lambda}(g,\bar{\eta },x^{
*})-\alpha f(\bar{\eta }))\mu_{th,db}(dz),\]
where, recalling that $p(x,z)=\frac {\rho (x,z)-1}{\rho (x,z)}$ and $\bar{
\eta}_{x^{*}}=\bar{\eta }-\delta_{x^{*}}$, 
\begin{eqnarray*}
\bar {H}_z^{\lambda}(g,\bar{\eta },x^{*})&=&\prod_{x\in\bar{\eta}_{
x^{*}}}((1-p(x,z))\bar {g}(x)+p(x,z))\\
&&\qquad\times\prod_{i=1}^{k(z)}\int_E\bar {g}(y_i)q(x^{*},z,dy_i
).\end{eqnarray*}

Note that if $\frac {k(z)}{\lambda}\rightarrow\zeta$ as $\lambda\rightarrow
\infty$, then calculating as in 
Section~\ref{sectdb}, 
\begin{eqnarray*}
&&H_z^{\lambda}(g,\eta )=\lambda^{-k(z)}\int_{[0,\lambda ]^{k(z)}}\bigg
[\prod_{(x,u)\in\eta ,u\neq u^{*}(\eta ,v^{*})}g(x,\rho (x,z){\cal J}^{
\lambda}_{r(\cdot ,z)}(x,u,\eta ,v^{*}))\\
&&\qquad\qquad\qquad\qquad\qquad\qquad\qquad\times\prod_{i=1}^{k(
z)}\int_Eg(y_i,v_i)q(x^{*},z,dy_i)\bigg]dv_1\ldots dv_{k(z)}\\
&&\rightarrow\int_0^{\infty}\bigg[\zeta e^{-\zeta v^{*}}\prod_{(x
,u)\in\eta ,u\neq u^{*}(\eta ,v^{*})}g(x,\rho (x,z)(u-{\bf 1}_{\{
u>u^{*}\}}(u^{*}-v^{*})\frac {r(x,z)}{r(x^{*},z)}))\\
&&\qquad\qquad\times\int_Eg(y,v^{*})q(x^{*},z,dy)\\
&&\qquad\qquad\times\exp\{-\zeta\int_E\int_{v^{*}}^{\infty}(1-g(y
,v))q(x^{*},z,dy)dv\}\bigg]dv^{*}\\
&&\equiv H^{\infty}_{\zeta}(g,\eta ).\end{eqnarray*}
Consequently, at least in the simple setting when 
$\mu^{\lambda}_{th,db}({\Bbb U})<\infty$ and the various parameters are continuous,
if we assume that as
$\lambda\rightarrow\infty$, for each $\varphi\in C_b({\Bbb R}\times 
{\Bbb U})$,
\[\int_{{\Bbb U}}\varphi (\frac {k(z)}{\lambda},z)\mu^{\lambda}_{
th,db}(dz)\rightarrow\int_{{\Bbb U}}\int_0^{\infty}\varphi (\zeta 
,z)\mu_{\zeta}(d\zeta ,z)\mu^{\infty}_{th,db}(dz),\]
where $\mu_{\zeta}(d\zeta ,z)$ is a probability distribution on $
[0,\infty )$,  
then $A_{th,db}f(\eta )$ converges to
\[A_{th,db}^{\infty}f(\eta )=\int_{{\Bbb U}}\int_0^{\infty}{\bf 1}_{
\{\bar{\eta }(r,z)>0\}}(H^{\infty}_{\zeta}(g,\eta )-f(\eta ))\mu_{
\zeta}(d\zeta ,z)\mu^{\infty}_{th,db}(dz).\]
\[\]
If $\int r(x,z)\Xi (dx)>0$, define
\[\beta (x^{*},\Xi )=\frac {r(x^{*},z)}{\int_Er(x,z)\Xi (dx)},\quad\]
and
\begin{eqnarray*}
&&{\cal H}_z(g,\Xi )=\int_0^{\infty}\int_0^{\infty}\bigg[\exp\{-\int_
E\frac 1{\rho (x,z)}h(x)\Xi (dx)\}\\
&&\qquad\times\int_E\beta (x^{*},\Xi )\exp\{-\zeta\int_Eh(y)q(x^{
*},z,dy)dv)\}\Xi (dx^{*})\bigg]\mu_{\zeta}(d\zeta ,z).\end{eqnarray*}
The projected generator 
then becomes
\[\alpha A^{\infty}_{th,db}f(\Xi )=\int_{{\Bbb U}}{\bf 1}_{\{\Xi 
(r,z)>0\}}({\cal H}_z(g,\Xi )-\alpha f(\Xi ))\mu^{\infty}_{th,db}
(dz).\]

\subsection{Immigration}
Immigration can be modeled by simply assigning each 
new immigrant a randomly chosen level.  This approach 
gives a generator of the form
\[A_{im}f(\eta )=\int_{{\Bbb U}_{im}}f(\eta )(\lambda^{-1}\int_0^{
\lambda}g(x(z),v)dv-1)\mu_{im}(dz)=\int_{{\Bbb U}_{im}}f(\eta )(\bar {
g}(x(z))-1)\mu_{im}(dz),\]
which gives
\[\alpha A_{im}f(\bar{\eta })=\int_{{\Bbb U}_{im}}\alpha f(\bar{\eta }
)(\bar {g}(x(z))-1)\mu_{im}(dz).\]
Again setting $h(x)=\int_0^{\infty}(1-g(x,u))du$, 
replacing $\mu_{im}$ by $\lambda\mu_{im}$, and passing to the limit as 
$\lambda\rightarrow\infty$ 
gives
\[A_{im}f(\eta )=-\int_{{\Bbb U}_{im}}f(\eta )h(x(z))\mu_{im}(dz)
,\]
and integrating out the levels
\[\alpha A_{im}f(\bar{\eta })=-\int_{{\Bbb U}_{im}}\alpha f(\bar{
\eta })h(x(z))\mu_{im}(dz)\]
which implies
\[\frac d{dt}\int_Eh(x)\Xi_t(dx)=\int_{{\Bbb U}_{im}}h(x(z))\mu_{
im}(dz),\]
as we would expect.

\subsection{Independent and exchangeable motion} 
Typically, population models assume independent motion 
or mutation causing individual types to change between 
birth/death events.  Some models allow common 
stochastic effects to influence type changes so that 
particle types evolve in an exchangeable fashion.  In 
either case, we assume the existence of a collection of 
process generators $\{B_n\}$, where $B_n$ determines a process 
with state space $E^n$, $B_n$ is exchangeable in the sense 
that if $(X_1,\ldots ,X_n)$ is a solution of the martingale problem 
for $B_n$, then any permutation of the indices $(X_{\sigma_1},\ldots 
,X_{\sigma_n})$ 
also gives a solution of the martingale problem for $B_n$, 
and the $B_n$ are consistent in the sense that if 
$(X_1,\ldots ,X_{n+1})$ is a solution of the martingale problem for 
$B_{n+1}$, then $(X_1,\ldots ,X_n)$ is a solution of the martingale 
problems for $B_n$.  Of course, if $B_n$ is the generator for 
$n$ independent particles, each with generator $B_1$, then the 
collection $\{B_n\}$ has the desired properties.  

To combine motion with the other possible elements of a 
model described above, we need a sufficiently rich class 
of function $g(x,u)$ such that for each $n$, and fixed 
$u_1,\ldots ,u_n$, $\prod_{i=1}^ng(x_i,u_i)$ gives a function in the domain of 
$B_n$.  In the independent case, this requirement simply 
means that $g(x,u)$ is in the domain of $B\equiv B_1$, and 
\[B_{|\eta |}f(\eta )=f(\eta )\sum_{(x,u)\in\eta}\frac {Bg(x,u)}{
g(x,u)}.\]

For finite $\lambda$, if $\bar{\eta }(E)<\infty$, then the motion generator is 
just given by
\[\hat {B}f(\eta )=B_{|\eta |}\prod_{(x,u)\in\eta}g(x,u).\]
For $\lambda =\infty$, since we assume that $g(x,u)\equiv 1$ for $
u\geq u_g$, 
the same formula works provided $\eta (E\times [0,u_g])<\infty$.

For models with infinitely many particles with levels 
below a fixed level, we can require the existence of a 
sequence $K_k\subset E$ such that $\cup_kK_k=E$ and  
$\eta (K_k\times [0,u_0])<\infty$ for each $k$ and $u_0$.  Requiring $
g(x,u)=1$ 
and
$B_1g(x,u)=0$ for $(x,u)\notin K_{k_g}\times [0,u_g]$ for some $k_g$ 
would give
\begin{equation}\hat {B}f(\eta )=B_{\eta (K_k\times [0,u_g])}\prod_{
(x,u)\in\eta_{|K_k\times [0,u_g]}}g(x,u).\label{mot2}\end{equation}
Note that this condition simply places restrictions on 
the size or direction of jumps by the motion process.

For finite $\lambda$ and $\bar{\eta }(E)<\infty$,
\[\alpha\hat {B}f(\bar{\eta })=B_{|\eta |}\prod_{x\in\bar{\eta}}\bar {
g}(x),\]
and similarly for (\ref{mot2}).  For $\lambda =\infty$, a general 
derivation for exchangeable but not independent motion is not 
clear, but for independent motion, observing that
$Bg=B(g-1)$ 
we have
\[\alpha\hat {B}f(\Xi )=-e^{-\int_Eh(x)\Xi (dx)}\int Bh(x)\Xi (dx
).\]

\subsection{Selecting a random sample}\label{sectransam}
The various recipes described above allow one to  
construct population models in a way that 
parent-offspring relationships can be identified knowing 
the evolution of the state in the model.  In particular, 
one can select a random ``sample'' from an appropriately
 finite region of the type space (even in 
the $\lambda =\infty$ case) and trace its genealogy.  For example, 
let $C\subset E$ satisfy $\bar{\eta }(t,C)<\infty$ in the $\lambda 
<\infty$ case and 
$\Xi (t,C)<\infty$ in the $\lambda =\infty$ case.  Then the set of particles with 
types in $C$ at the $n$ lowest levels is a uniform random 
sample of size $n$ drawn from the subpopulation of 
particles with types in $C$ and the genealogies of these 
$n$ particles can be traced by following the evolution of 
the levels back in time.

If the levels are constant in time, then as noted in 
Remark \ref{geneaeq} and Section 5 of \cite{DK99a}, one 
can define a family of counting processes and a system 
of stochastic equations driven by these counting 
processes whose solution gives the desired genealogy.
Tracing the genealogy for a model with moving levels is 
much less elegant; however, complete genealogical 
information is present in the levels and the stochastic 
inputs of the birth events.

%% file: LDEXAMP1.tex

$\frac {}{}$
\section{Examples} \label{examples}

So far we have largely performed formal calculations, 
not proofs.  In this section we illustrate our results in 
some specific examples and here, unless otherwise 
stated, our results are mathematically rigorous.  
In 
\S\ref{sectslfv}, we present two different approaches to 
the process known as the spatial $\Lambda$-Fleming-Viot process 
(which we shall also define).  The first, based on 
one-for-one replacement, yields, in the high intensity 
limit, the process with levels of \cite{VW15} (under 
somewhat weaker conditions).  The second, based on 
discrete births of Poisson numbers of offspring and death
by independent thinning, corresponds in the prelimit to 
the particle system studied in \cite{BEH09}.  In 
\S\ref{sectslfv2}, we extend this second approach to 
discrete birth mechanisms in which the number of 
offspring is no longer required to be Poisson.  This 
yields a new class of population models, in which the 
replacement mechanism mirrors that of the spatial 
$\Lambda$-Fleming-Viot process, but the population intensity can 
vary with spatial position.  In particular, these models 
provide one approach to combining ecology and genetics 
as described in the introduction.  
In \S\ref{sectbranching}, we revisit branching processes 
and the Dawson-Watanabe superprocess.  In \S\ref{sectMoran}, we 
use one for one replacement, in the special case in 
which just two individuals are involved in each event, 
to recover, in particular, the lookdown construction 
of~\cite{GLW05} for a spatially interacting Moran model.  
 In 
\S\ref{sectspde}, we use the lookdown construction to 
derive a stochastic partial differential equation as the 
limit of rescaled spatially interacting Moran models of 
the type discussed in \S\ref{sectMoran}.  Finally, in 
\S\ref{sectvote}, we give a lookdown construction for a 
class of voter models and use the construction to give a 
heuristic argument for a result of Mueller and Tribe 
\cite{MT95} showing that the rescaled voter model 
converges to a solution of the stochastic partial 
differential equation obtained in \S\ref{sectspde}.  

\subsection{Spatial $\Lambda$-Fleming-Viot 
process}\label{sectslfv} The spatial $\Lambda$-Fleming-Viot 
process was introduced in \cite{Eth08} and rapidly 
developed by a number of authors 
\cite{BEV10,BEH09,VW15}.  The primary motivation is to 
model a spatially distributed population in such a way 
that the distribution of the population is stable in space 
and one can recover the genealogical trees relating 
individuals in a sample from the population in an 
analytically tractable way.  A survey can be found in 
\cite{BEV13}.  The process is driven by spatially 
distributed birth/death events in which a significant 
fraction of the local population is replaced.  The 
location, spatial extent, and `impact' of these events (by 
which we mean the proportion of the local population 
replaced in an event) is determined by a Poisson random 
measure, and stability of the population is maintained by 
ensuring that the numbers of births and deaths balance.  

We now explicitly distinguish between the location of a 
particle $x\in {\Bbb R}^d$ and its type $\kappa\in {\Bbb K}$. 
Let $E={\Bbb R}^d\times {\Bbb K}$,
${\Bbb U}={\Bbb R}^d\times [0,1]\times [0,\infty )$, and 
$\mu =\ell^d\times\nu^1(w,d\zeta )\times\nu^2(dw)$ where $\ell^d$ is
Lebesgue measure on ${\Bbb R}^d$, $\nu^2$ is a $
\sigma$-finite measure 
on $[0,\infty )$ 
and $\nu^1$ is a transition function from $[0,\infty )$ to $[0,1]$.  

If $C\subset {\Bbb R}^d$ is Borel measurable, then $|C|=\ell^d(C)$. If $
C$ is 
a finite or countable set, 
then $|C|$ will denote the number of elements in $C$.  
Which interpretation applies should be clear in context.

Each point in ${\Bbb U}$ specifies a point $y\in {\Bbb R}^d$, $w\in 
[0,\infty )$
and $\zeta\in [0,1]$.  The corresponding reproduction event will affect
the population in the ball $D_{y,w}\subseteq {\Bbb R}^d$ centered at $
y$ with radius $w$, 
and $\zeta$ will determine the impact within the ball.  
The model is driven by a space-time Poisson random 
measure on ${\Bbb U}\times [0,\infty )$ with mean measure $\mu\times
\ell$.  If a 
birth/death event occurs at time $t$ corresponding to 
$(y,\zeta ,w)\in {\Bbb U}$, an
individual located in $D_{y,w}$ is selected at random to be 
the `parent', a fraction $\zeta$ of the individuals in $D_{y,w}$ are 
killed and replaced by individuals of the same type as 
the parent, with the locations of the new individuals 
uniformly distributed over $D_{y,w}$. 

We will give two constructions of processes following 
this recipe which differ 
substantially for finite $\lambda$ but,
under conditions for which both constructions are valid,  
yield the same measure-valued model in the limit. 
The first 
construction follows ideas of V\'eber and Wakolbinger 
\cite{VW15}.

In order to rigorously define the generators of our 
processes, we will need to restrict the domains.  In both 
cases the domains will be subsets of 
\begin{eqnarray}
	{\cal D}_{\lambda}&=&\{f(\eta )=\prod_{(x,\kappa ,u)\in\eta} g(x,\kappa,u):0\leq 
g\leq 1,\exists\mbox{\rm \ compact }K_g\subset {\Bbb R}^d,0<u_g\leq
\lambda ,\non\\
&&\qquad\qquad\qquad\qquad\qquad g(x,\kappa ,u)=1\mbox{\rm \ for }
(x,u)\notin K_g\times [0,u_g)\}\label{scrddef}\\
{\cal D}_{\infty}&=&\cup_{\lambda >0}{\cal D}_{\lambda}.\non\end{eqnarray}
Without loss of generality, we can assume that 
$K_g=D_{0,\rho_g}$ for $0<\rho_g<\infty$. 

Consider $A_{dr,2}$ defined in (\ref{dr2lam}).  Recall that 
with this mechanism, for each replacement event, we 
specify the probability $r(x)$ that an individual of type $x$ 
is replaced and the parent is taken to be the individual 
chosen to be replaced that has the lowest level.  For an 
event corresponding to $z=(y,\zeta ,w)$, let $r(x,z)$ be $\zeta {\bf 1}_{
D_{y,w}}(x)$ 
and for $(x,\kappa )\in E$, the transition function $q$ of 
Section~\ref{sectoneforone} becomes 
\[q(x,\kappa ,z,dx'\times d\kappa')=\upsilon_{y,w}(dx')\delta_{\kappa}
(d\kappa'),\]
where $\upsilon_{y,w}$ is the uniform distribution over the ball 
$D_{y,w}$, that is, the offspring have the same type as the 
parent and are independently and uniformly distributed 
over the ball.  Consequently, $\hat {g}$ in (\ref{ghatdef}) becomes
\[\hat {g}_{y,w}(\kappa ,u)\equiv\int g(x',\kappa ,u)\upsilon_{y,
w}(dx').\]
In addition, recalling that $\bar {g}(x,\kappa )=\lambda^{-1}\int_
0^{\lambda}g(x,\kappa ,u)du$, we define
\begin{equation}\label{gywbar}
\bar {g}_{y,w}(\kappa )\equiv\int\bar {g}(x',\kappa )\upsilon_{
y,w}(dx')=\frac 1{\lambda}\int_0^{\lambda}\hat {g}_{y,w}(\kappa ,
u)du=\int\frac 1{\lambda}\int_0^{\lambda}g(x',\kappa ,u)du\upsilon_{
y,w}(dx').\end{equation}
We postpone giving precise conditions on 
$\nu^1$  and $\nu^2$ until we have formally derived the
generators.

We define 
\begin{equation}\eta_{y,w}=\sum_{(x,\kappa ,u)\in\eta :x\in D_{y,
w}}\delta_{(x,\kappa ,u)}\mbox{\rm \  and  }\eta_{y,w}^g=\sum_{(x
,\kappa ,u)\in\eta :x\in D_{y,w},u\leq u_g}\delta_{(x,\kappa ,u)}.\label{etatest}\end{equation}
That is $\eta^g_{y,w}=\eta (\cdot\cap D_{y,w}\times {\Bbb K}
\times [0,u_g))$ 
is the restriction of $\eta$ to $D_{y,w}\times {\Bbb K}\times [
0,u_g)$. 
From (\ref{dr2inf}) and~(\ref{expn}),   
$A^{\infty}_{dr,2}$ is given by  
\begin{equation}A^{\infty}_{dr,2}f(\eta )=f(\eta )\int_{{\Bbb R}^
d\times [0,1]\times [0,\infty )}\Bigg(\frac {\sum_{S\subset\eta^g_{
y,w}}H(g,\hat {g},S,y,\zeta ,w)}{\prod_{(x,\kappa ,u)\in\eta^g_{y
,w}}g(x,\kappa ,u)}-1\Bigg)dy\nu^1(w,d\zeta )\nu^2(dw),
\label{vwgen}\end{equation}
where 
\[H(g,\hat {g},S,y,\zeta ,w)=\prod_{(x,\kappa ,u)\in S}(\hat {g}_{
y,w}(\kappa^{*}(S),u)\zeta )\prod_{(x,\kappa ,u)\in\eta_{y,w}^g,(
x,\kappa ,u)\notin S}((1-\zeta )g(x,\kappa ,u)),\]
$\kappa^{*}(S)$ being the type of the lowest level particle in $S$. 
 $A^{\infty}_{dr,2}$ is the generator for the lookdown 
construction of \cite{VW15}. Again, for an event 
corresponding to $(y,\zeta ,w)$, a particle in $D_{y,w}$ is 
involved in the event with probability $\zeta$.

The relationship between the martingale problems for finite and infinite 
$\lambda$ is particularly simple in this setting. 
For finite $\lambda$, $A^{\lambda}_{dr,2}f(\eta )=A_{dr,2}^{\infty}
f(\eta )$ provided $u_g\leq\lambda$.  
Consequently, any solution of the martingale problem for 
$A_{dr,2}^{\infty}$ restricted to levels in $[0,\lambda ]$  gives a solution of 
the martingale problem for $A_{dr,2}^{\lambda}$.  In particular, 
existence and 
uniqueness for $A_{dr,2}^{\lambda}$ for all $\lambda >0$ implies existence and 
uniqueness for 
$A^{\infty}_{dr,2}$. 

Setting $\bar{\eta}_{y,w}=\bar{\eta }(\cdot\cap D_{y,w}\times {\Bbb K}
)$, for finite $\lambda$,
\begin{equation}\alpha A^{\lambda}_{dr,2}f(\bar{\eta })=\alpha f(
\bar{\eta })\int_{{\Bbb R}^d\times [0,1]\times [0,\infty )}\Bigg(\frac {
\sum_{S\subset\bar{\eta}_{y,w}}\bar {H}(\bar {g},\hat {g},S,y,\zeta 
,w)}{\prod_{(x,\kappa )\in\bar{\eta}_{y,w}}\bar {g}(x,\kappa )}-1\Bigg
)dy\nu^1(w,d\zeta )\nu^2(dw)\label{lfvdisc}\end{equation}
where, recalling the notation defined in~(\ref{gywbar}), 
\begin{equation}\bar {H}(\bar {g},\hat {g},S,y,\zeta ,w)=\frac 1{
|S|}\sum_{(x,\kappa )\in S}(\bar {g}_{y,w}(\kappa )\zeta )^{|S|}\prod_{
(x,\kappa )\in\bar{\eta}_{y,w},(x,\kappa )\notin S}((1-\zeta )\bar {
g}(x,\kappa )).\label{hbarlam}\end{equation}

Finally, setting $h^{*}_{y,w}(\kappa )=\int_0^{\infty}(1-\hat {g}_{
y,w}(\kappa ,u))du$ (recall 
$h(x,\kappa )=\int_0^{\infty}(1-g(x,\kappa ,u))du$) and 
\[{\Bbb H}_1(h^{*}_{y,w},\Xi ,y,\zeta ,w)=\frac 1{\Xi (D_{y,w}\times 
{\Bbb K})}\int_{D_{y,w}\times {\Bbb K}}e^{-\zeta h^{*}_{y,w}(\kappa 
)\Xi (D_{y,w}\times {\Bbb K})}\Xi (dx\times d\kappa ),\]
we have
\begin{eqnarray}
&&\alpha A^{\infty}_{dr,2}f(\Xi )\label{firstlam}\\
&&\quad =e^{-\int h(x,\kappa )\Xi (dx,d\kappa )}\non\\
&&\qquad\times\int_{{\Bbb R}^d\times [0,1]\times [0,\infty )}({\Bbb H}_
1(h^{*}_{y,w},\Xi ,y,\zeta ,w)e^{\zeta\int_{D_{y,w}\times {\Bbb K}}
h(x,\kappa )\Xi (dx,d\kappa )}-1)dy\nu^1(w,d\zeta )\nu^2(dw).\non\end{eqnarray}
Note that if $\Xi$ is a solution of the martingale problem 
for $\alpha A^{\infty}_{dr,2}$ and $\Xi (0,dx\times {\Bbb K})$ is Lebesgue measure, then 
$\Xi (t,dx\times {\Bbb K})$ is Lebesgue measure for all $t\geq 0$.  (Consider 
the generator with $h$ not depending on $\kappa$.)

Before establishing conditions under which the 
construction above is valid, let us describe an alternative 
lookdown construction of the spatial $\Lambda$-Fleming-Viot 
process employing discrete births (Section~\ref{sectdb}) and
independent thinning (Section~\ref{sectthin})
as in Section \ref{sectevent}.  With $z=(y,\zeta ,w)$ as above, 
the thinning parameter is 
\begin{equation}\rho (x,\kappa ,z)=1+\frac {\zeta}{1-\zeta}{\bf 1}_{
D_{y,w}}(x).\label{thin}\end{equation}
As we saw in Section~\ref{sectthin}, this assumption 
ensures that the probability that an existing  individual 
(other than the parent) dies is zero outside the ball 
$D_{y,w}$ and $\zeta$ within it.  

If there is at least one individual in $D_{y,w}$ (to serve as 
parent), the discrete 
birth event corresponding to $z$ produces a Poisson 
number of offspring with parameter $\lambda\alpha_z$ conditioned to 
be positive, where 
$\alpha_z=\zeta |D_{y,w}|$, $r(x,z)={\bf 1}_{D_{y,w}}(x)$, and 
\[q(x,\kappa ,z,dx',d\kappa')=\upsilon_{y,w}(dx')\delta_{\kappa}(
d\kappa').\]
The finite intensity model is then essentially that 
considered in \cite{BEH09}, differing only in the 
assumptions that the parent is selected before the 
thinning and the offspring distribution is conditioned
to be positive.  Note that the definition of $r$ 
in this construction is different from the definition in 
the previous construction.  There, $r$ determined the 
chance of being involved in the event; here we use it to 
weight the chance of being a parent.  This distinction 
becomes important in modelling different forms of 
natural selection when we would choose $r$ to depend on 
type.  

As in Section \ref{sectevent}, but with a slight change of 
notation, let 
$\eta (y,w)=\eta (D_{y,w}\times {\Bbb K}\times [0,\lambda ])$ and
\begin{equation}A^{\lambda}_{th,db}f(\eta )=\int_{{\Bbb U}}{\bf 1}_{
\{\eta (y,w)>0\}}(H_z^{\lambda}(g,\eta )-f(\eta ))(1-e^{-\alpha_z
\lambda})dy\nu^1(w,d\zeta )\nu^2(dw).\label{lamf12}\end{equation}
We introduce the factor $1-e^{-\alpha_z\lambda}$ in the event measure, 
and then condition on there being at least one offspring. 
If $\eta (D_{y,w}\times {\Bbb K}\times [0,\lambda ])=\eta (y,w)\neq 
0$, we obtain an expression for 
$H_z^{\lambda}(g,\eta )$ by partitioning on the lowest level selected for 
the offspring.   Since the levels $\{v_i\}$ selected for the 
offspring are the jump times in $[0,\lambda ]$ of a Poisson 
process with intensity $\alpha_z$, this yields 
\begin{eqnarray}
H_z^{\lambda}(g,\eta )&=&\prod_{(x,\kappa ,u)\in\eta ,x\notin D_{
y,w}}g(x,\kappa ,u)\label{hzlamdef}\\
&&\qquad\times\frac 1{1-e^{-\alpha_z\lambda}}\int_0^{\lambda}\Big
[\alpha_ze^{-\alpha_zv^{*}}\hat {g}_{y,w}(\kappa^{*},v^{*})e^{-\alpha_
z\int_{v^{*}}^{\lambda}(1-\hat {g}_{y,w}(\kappa^{*},v))dv}\nonumber\\
&&\qquad\qquad\times\prod_{(x,\kappa ,u)\in\eta ,x\in D_{y,w},u\neq 
u^{*}}g(x,\kappa ,\frac 1{1-\zeta}{\cal J}_{y,w}^{\lambda}(x,u,\eta 
,v^{*}))\Big]dv^{*},\nonumber\end{eqnarray}
where $(x^{*},\kappa^{*},u^{*})$ is the point in $\eta$ satisfying 
$x^{*}\in D_{y,w}$ and
\[u^{*}=\mbox{\rm argmax}\{\frac {\lambda -u}{\lambda -v^{*}}:(x,
\kappa ,u)\in\eta ,x\in D_{y,w},u\geq v^{*}\}\cup \{\frac u{v^{*}}
:(x,\kappa ,u)\in\eta ,x\in D_{y,w},u\leq v^{*}\},\]
and ${\cal J}_{y,w}^{\lambda}(x,u,\eta ,v^{*})$ is obtained as in (\ref{levtran}) with 
$r={\bf 1}_{D_{y,w}}$.  Recall that we thin the existing population 
after we select the parent, and the thinning is 
accomplished by multiplying ${\cal J}^{\lambda}_{y,w}$  by $\rho$ defined in 
(\ref{thin}).

Let $\bar{\eta}_{|D_{y,w}}$ denote $\bar{\eta}$ restricted to $D_{
y,w}\times {\Bbb K}$.  Since 
conditional on $\bar{\eta}$ and $v^{*}$, $(x^{*},\kappa^{*})$ is selected uniformly at 
random from $\bar{\eta}_{|D_{y,w}}$ and, for $u\neq u^{*}$ (see 
Lemma~\ref{uniflev}), the ${\cal J}_{y,w}^{\lambda}(x,u,\eta ,v^{
*})$ are independent 
and uniform over $[0,\lambda ]$, partitioning on the level of the 
lowest offspring, define 
\begin{eqnarray*}
{\cal H}^{\lambda}_z(\bar {g},\bar{\eta })&=&\frac 1{|\bar{\eta}_{
|D_{y,w}}|}\sum_{(x^{*},\kappa^{*})\in\bar{\eta}_{|D_{y,w}}}\prod_{
(x,\kappa )\in\bar{\eta}_{|D_{y,w}},(x,\kappa )\neq (x^{*},\kappa^{
*})}((1-\zeta )\bar {g}(x,\kappa )+\zeta )\\
&&\qquad\qquad\times\frac 1{1-e^{-\lambda\alpha_z}}\int_0^{\lambda}
\alpha_ze^{-\alpha_zv^{*}}\hat {g}_{y,w}(\kappa^{*},v^{*})e^{-\alpha_
z\int_{v^{*}}^{\lambda}(1-\hat {g}_{y,w}(\kappa^{*},v))dv}dv^{*}\\
&=&\frac 1{|\bar{\eta}_{|D_{y,w}|}|}\sum_{(x^{*},\kappa^{*})\in\bar{
\eta}_{|D_{y,w}}}\prod_{(x,\kappa )\in\bar{\eta}_{|D_{y,w}},(x,\kappa 
)\neq (x^{*},\kappa^{*})}((1-\zeta )\bar {g}(x,\kappa )+\zeta )\\
&&\qquad\qquad\times\frac 1{1-e^{-\lambda\alpha_z}}(e^{-\alpha_z\int_
0^{\lambda}(1-\hat {g}_{y,w}(\kappa^{*},v))dv}-e^{-\lambda\alpha_
z})\\
&=&\frac 1{|\bar{\eta}_{|D_{y,w}}|}\sum_{(x^{*},\kappa^{*})\in\bar{
\eta}_{|D_{y,w}}}\prod_{(x,\kappa )\in\bar{\eta}_{|D_{y,w}},(x,\kappa 
)\neq (x^{*},\kappa^{*})}((1-\zeta )\bar {g}(x,\kappa )+\zeta )\\
&&\qquad\qquad\times\frac 1{1-e^{-\lambda\alpha_z}}(e^{-\zeta |D_{
y,w}|h^{*}_{y,w}(\kappa^{*})}-e^{-\alpha_z\lambda}),\end{eqnarray*}
where, as before, $h^{*}_{y,w}(\kappa )=\int_0^{\infty}(1-\hat {g}_{
y,w}(\kappa ,u))du$.  To 
understand this quantity, recall first that in our 
discrete births model, the parent is eliminated from the 
population.  Next, 
for points within $D_{y,w}$, they survive with probability 
$(1-\zeta )$, otherwise they are removed (giving the product on 
the right side of the first line).  The final term 
corresponds to the offspring (recalling the notation 
$\bar {g}_{y,w}(\kappa )$ from~(\ref{gywbar}) and that we have
conditioned on there being at least one offspring).  Then 
\[\alpha A^{\lambda}_{th,db}f(\bar{\eta })=\alpha f(\bar{\eta })\int_{
{\Bbb U}}{\bf 1}_{\{\bar{\eta }(y,w)>0\}}(\frac {{\cal H}_z^{\lambda}
(\bar {g},\bar{\eta })}{\prod_{(x,\kappa )\in\bar{\eta}_{|D_{y,w}}}
\bar {g}(x,\kappa )}-1)(1-e^{-\alpha_z})dy\nu^1(w,d\zeta )\nu^2(d
w).\]

Note that $\alpha A^{\lambda}_{th,db}$ constructed here is not the same as 
$\alpha A^{\lambda}_{dr,2}$ given in (\ref{lfvdisc}).  Here, at each 
birth/death event, 
existing particles 
are randomly killed and an independent number of new 
particles are created while in the previous construction, 
the number of births equaled the number of deaths. 
However,
taking $\lambda\rightarrow\infty$, by (\ref{hlim}), 

\begin{equation}A^{\infty}_{th,db}f(\eta )=\int_{{\Bbb U}}{\bf 1}_{
\{\eta (y,w)>0\}}(H_z(g,\eta )-f(\eta ))dy\nu^1(w,d\zeta )\nu^2(d
w),\label{dbth1}\end{equation}
with
\begin{eqnarray*}
H_z(g,\eta )&=&\prod_{(x,\kappa ,u)\in\eta ,x\notin D_{y,w}}g(x,\kappa 
,u)\\
&&\qquad\times\int_0^{\infty}\Big[\alpha_ze^{-\alpha_zv^{*}}\hat {
g}_{y,w}(\kappa^{*},v^{*})e^{-\alpha_z\int_{v^{*}}^{\infty}(1-\hat {
g}_{y,w}(\kappa^{*},v))dv}\\
&&\qquad\qquad\times\prod_{(x,\kappa ,u)\in\eta ,x\in D_{y,w},u>u^{
*}}g(x,\kappa ,\frac 1{1-\zeta}(u-u^{*}+v^{*}))\\
&&\qquad\qquad\times\prod_{(x,\kappa ,u)\in\eta ,x\in D_{y,w},u<u^{
*}}g(x,\kappa ,\frac 1{1-\zeta}u)\Big]dv^{*}.\end{eqnarray*}
Just as in Lemma~\ref{uniflev}, (and using~(\ref{hlim})), 
it is easy to see that $\eta^{*}$ satisfying
\begin{eqnarray*}
\int gd\eta^{*}=\sum_{(x,\kappa ,u)\in\eta ,x\in D_{y,w},u>u^{*}}
g(x,\kappa ,\frac 1{1-\zeta}(u-u^{*}+v^{*}))\\
+\sum_{(x,\kappa ,u)\in\eta ,x\in D_{y,w},u<u^{*}}g(x,\kappa ,\frac 
1{1-\zeta}u)\end{eqnarray*}
is conditionally Poisson with Cox measure 
$(1-\zeta ){\bf 1}_{D_{y,w}}(x)\Xi (dx,d\kappa)$ 
and recalling the definition of $h^*_{y,w}(\kappa)$ from just below
equation~(\ref{hbarlam}) an integration by parts gives
\begin{eqnarray*}
\int_0^{\infty}\alpha_ze^{-\alpha_zv^{*}}\hat {g}_{y,w}(\kappa^{*}
,v^{*})e^{-\alpha_z\int_{v^{*}}^{\infty}(1-\hat {g}_{y,w}(\kappa^{
*},v))dv}dv^{*}
=e^{-\alpha_zh_{y,w}^{*}(\kappa^{*})}.\end{eqnarray*}
Averaging~(\ref{dbth1}) gives 
\begin{eqnarray*}
&&\alpha A^{\infty}_{th,db}f(\Xi )=e^{-\int_{{\Bbb R}^d\times {\Bbb K}}
h(x,\kappa )\Xi (dx,d\kappa )}\\
&&\qquad\qquad\qquad\qquad\times\int_{{\Bbb U}}({\Bbb H}_2(h^{*}_{
y,w},\Xi ,z)e^{\zeta\int_{D_{y,w}\times {\Bbb K}}h(x,\kappa )\Xi 
(dx,d\kappa )}-1)dy\nu^1(w,d\zeta )\nu^2(dw),\end{eqnarray*}
where
\[{\Bbb H}_2(h^{*}_{y,w},\Xi ,z)=\frac 1{\Xi (D_{y,w}\times {\Bbb K}
)}\int_{D_{y,w}\times {\Bbb K}}e^{-\zeta |D_{yw}|h_{y,w}^{*}(\kappa 
)}\Xi (dx\times d\kappa ),\]
(defined to be $1$ if $\Xi(D_{y,w}\times {\Bbb K})=0$)
which, in general, differs from ${\Bbb H}_1$.
However, if $\Xi$ is a solution of the martingale problem 
for $\alpha A_{th,db}^{\infty}$ with
 $\Xi (0,dx\times {\Bbb K})$  Lebesgue measure, then $\Xi (t,dx\times 
{\Bbb K})$ is 
Lebesgue measure for  all $t\geq 0$ and
${\Bbb H}_2(h_{y,w}^{*},\Xi ,z)={\Bbb H}_1(h^{*}_{y,w},\Xi ,z)$.  Consequently, in this case, $
\Xi$ is also
a solution of the martingale problem for
 $\alpha A^{\infty}_{dr,2}$ in the previous construction.

Our calculations so far in this subsection have been 
entirely formal.  We now turn to actually constructing 
the processes that correspond to the generators 
described above.  

\subsubsection{First construction of spatial $\Lambda$-Fleming-Viot with levels}
\label{firstslfv}
The process corresponding to $A^{\infty}_{dr,2}$ appears already in 
\cite{VW15}, but the strategy of our construction, based 
on writing down stochastic equations for the type of the 
particle at the $i$th level for each $i$, is somewhat 
different, and we obtain our process under somewhat 
weaker conditions.  In particular, for existence of our 
construction, we require 
\begin{equation}\int_{[0,1]\times (1,\infty )}\zeta w^d\nu^1(w,d\zeta 
)\nu^2(dw)<\infty\label{flreq1}\end{equation}
and 
\begin{equation}\left\{\begin{array}{ll}
\displaystyle\int_{[0,1]\times [0,1]}\zeta |w|^2\nu^1(w,d\zeta )\nu^
2(dw)<\infty&\mbox{\rm if }d=1,\\
&\\
\displaystyle\int_{[0,1]\times [0,1]}\zeta |w|^{2+d}\nu^1(w,d\zeta 
)\nu^2(dw)<\infty&\mbox{\rm if }d\geq 2;\end{array}
\right.\label{flreq1b}\end{equation}
while \cite{VW15} assumes
\begin{equation}\int_{[0,1]\times (0,\infty )}\zeta w^d\nu^1(w,d\zeta 
)\nu^2(dw)<\infty .\label{flreq2}\end{equation}
We should point out, however, that up to now, we do 
not have a proof of uniqueness for the system of 
stochastic equations under the weaker conditions, except in 
the case $d=1$ when uniqueness is proved in \cite{XZ17}.  The 
solution is unique under (\ref{flreq2}).

To rigorously cover the more general conditions, we need 
to be more careful in the description of  the generators, which for 
simplicity we will call $A^{\lambda}$ and $A^{\infty}$.  In particular, we appeal to 
the construction in Appendix \ref{sumbnd}.  With reference to 
(\ref{scrddef}), we 
restrict the domain to 
\[{\cal D}(A^{\infty})=\{f(\eta )=\prod_{(x,\kappa ,u)\in\eta}g(x
,\kappa ,u)\in {\cal D}_{\infty}:g(\cdot ,\kappa ,u)\in C^2({\Bbb R}^
d)\}.\]

To avoid additional complication of notation, we will also 
assume that for each $k=1,2,\ldots$,
\begin{equation}\nu^2(2^{-k},2^k)<\infty .\label{intbnd}\end{equation}
With the results of Appendix \ref{sumbnd}\ in mind,
define $\Gamma_0=\emptyset$ and for $k=1,2,\ldots$,
\begin{equation}\Gamma_k=D_{0,k}\times [0,1]\times [2^{-k},2^k].\label{gamk}\end{equation}
Set 
\[B_kf(\eta )=\int_{\Gamma_k-\Gamma_{k-1}}f(\eta )\Bigg(\frac {\sum_{
S\subset\eta^g_{y,w}}H(g,\hat {g},S,y,\zeta ,w)}{\prod_{(x,\kappa 
,u)\in\eta^g_{y,w}}g(x,\kappa ,u)}-1\Bigg)dy\nu^1(w,d\zeta )\nu^2
(dw),\]
where as before 
\[H(g,\hat {g},S,y,\zeta ,w)=\prod_{(x,\kappa ,u)\in S}(\hat {g}_{
y,w}(\kappa^{*}(S),u)\zeta )\prod_{(x,\kappa ,u)\in\eta_{y,w}^g,(
x,\kappa ,u)\notin S}((1-\zeta )g(x,\kappa ,u)).\]
Note that, writing $v_d$ for the volume of the unit ball,
$\lambda_k$ in (\ref{bk}) is 
\begin{equation}\lambda_k=v_dk^d\int_{2^{-k}}^{2^k}\nu^1(w,[0,1])
\nu^2(dw)-v_d(k-1)^d\int_{2^{-(k-1)}}^{2^{k-1}}\nu^1(w,[0,1])\nu^
2(dw).\label{lamklamfv}\end{equation}

The definition of $H_k$ is somewhat more complicated than 
the form used in Appendix \ref{sumbnd}, but arguments 
used there carry over immediately.   Let 
${\Bbb U}_k=(\Gamma_k-\Gamma_{k-1})\times ([0,1]\times D_{0,1})^{
\eta}$, 
and 
\[\nu_k(dy,d\zeta ,dw,\ldots ,dz_u,dv_u,\ldots )=\frac 1{\lambda_
k}dy\nu^1(w,d\zeta )\nu^2(dw)\prod_{(x,\kappa ,u)\in\eta}dz_u\upsilon_{
0,1}(dv_u),\]
that is, for each $k$,
we associate a pair of random variables $(Z_{k,u},V_{k,u})$ 
with each element of $\eta$, where $Z_{k,u}$ is uniformly  
distributed on $[0,1]$ and $V_{k,u}$ is uniformly distributed on 
$D_{0,1}$. We can index these random variables by $u$ since in 
our model the levels $u$ will be distinct. Then
\begin{equation}H_k(\eta ,y,\zeta ,w,(z,v)^{\eta})=\sum_{(x,\kappa 
,u)\in\eta}\left((1-{\bf 1}_{D_{y,w}}(x){\bf 1}_{[0,\zeta ]}(z))\delta_{
(x,\kappa ,u)}+{\bf 1}_{D_{y,w}}(x){\bf 1}_{[0,\zeta ]}(z)\delta_{
(y+wv,\kappa_{u^{*}},u)}\right),\label{hklamfv}\end{equation}
where $u^{*}=\min\{u:(x,\kappa ,u)\in\eta ,x\in D_{y,w},z_u\leq\zeta 
\}$.  Note that if 
$V$ is uniformly distributed on $D_{0,1}$, then $y+wV$ is 
uniformly distributed on $D_{y,w}$.

To verify Condition 
\ref{sumcnd2}.  We split $\Gamma_k$, setting
\begin{eqnarray*}
\Gamma_k&=&\Gamma_k^1\cup\Gamma_k^2\equiv D_{0,k}\times [0,1]\times 
[2^{-k},1]\cup D_{0,k}\times [0,1]\times (1,2^k],\\
\Gamma_{\infty}&=&\Gamma_\infty^1\cup\Gamma_\infty^2\equiv {\Bbb R}^d\times 
[0,1]\times (0,1]\cup {\Bbb R}^d\times [0,1]\times (1,\infty ),\end{eqnarray*}
and for $i=1,2$, define
\[B_k^if(\eta )=\int_{\Gamma^i_k-\Gamma^i_{k-1}}f(\eta )\Bigg(\frac {
\sum_{S\subset\eta^g_{y,w}}H(g,\hat {g},S,y,\zeta ,w)}{\prod_{(x,
\kappa ,u)\in\eta^g_{y,w}}g(x,\kappa ,u)}-1\Bigg)dy\nu^1(w,d\zeta 
)\nu^2(dw).\]
Recall the definition of $\eta_{y,w}^g$ from~(\ref{etatest}). For $i=2$, 
as in   
(\ref{dr2lam}), let $\{\xi^{\zeta}_{x,\kappa ,u}\}$ be independent with 
$P\{\xi^{\zeta}_{x,\kappa ,u}=1\}=1-P\{\xi_{x,\kappa ,u}^{\zeta}=
0\}=\zeta$.  Then
\begin{eqnarray}
&&|\sum_{k=m+1}^{\infty}B_k^2f(\eta )|=\Big|\int_{\Gamma_{\infty}^
2-\Gamma_m^2}\prod_{(x,\kappa ,u)\in\eta -\eta^g_{y,w}}g(x,\kappa 
,u)
\nonumber
\\
&&\qquad\times\left(E[\prod_{(x,\kappa ,u)\in\eta_{y,w}^g}(\xi^{\zeta}_{
x,\kappa ,u}\hat g_{y,w}(\kappa^{*},u)+(1-\xi^{\zeta}_{x,\kappa ,
u})g(x,\kappa ,u))]-f(\eta_{y,w}^g)\right)dy\nu^1(w,d\zeta )\nu^2
(dw)\Big|
\nonumber
\\
&&\leq\int_{\Gamma^2_{\infty}-\Gamma_m^2}\sum_{(x,\kappa ,u)\in\eta_{
y,w}^g}E[\xi^{\zeta}_{x,\kappa ,u}|\hat {g}_{y,w}(\kappa^{*},u)-1
|]dy\nu^1(w,d\zeta )\nu^2(dw)\\
&&\qquad +\int_{\Gamma^2_{\infty}-\Gamma_m^2}\sum_{(x,\kappa ,u)\in
\eta_{y,w}^g\cap D_{0,\rho_g}\times {\Bbb K}\times [0,u_g)}E[\xi^{
\zeta}_{x,\kappa ,u}|1-g(x,\kappa ,u)|]dy\nu^1(w,d\zeta )\nu^2(dw
)
\nonumber
\\
&&\leq\int_{\Gamma^2_{\infty}}|\eta_{y,w}^g|\frac {|D_{0,\rho_g}\cap 
D_{y,w}|}{|D_{y,w}|}\zeta dy\nu^1(w,d\zeta )\nu^2(dw)\\
&&\qquad +\int_{\Gamma^2_{\infty}}\eta (D_{0,\rho_g}\times {\Bbb K}
\times [0,u_g)){\bf 1}_{\{D_{y,w}\cap D_{0,\rho_g}\neq\emptyset \}}
\zeta dy\nu^1(w,d\zeta )\nu^2(dw)
\nonumber
\\
&&\leq |D_{0,\rho_g}|\int_{[0,1]\times (1,\infty )}\frac {|\eta_{
y,w}^g|}{|D_{y,w}|}{\bf 1}_{\{D_{y,w}\cap D_{0,\rho_g}\neq\emptyset 
\}}\zeta\nu^1(w,d\zeta )\nu^2(dw)
\nonumber
\\
&&\qquad +\eta (D_{0,\rho_g}\times {\Bbb K}\times [0,u_g))\int_{[
0,1]\times (1,\infty )}v_d(\rho_g+w)^d\zeta\nu^1(w,d\zeta )\nu^2(
dw),\label{bound on B2}\end{eqnarray}
where to obtain the first inequality we have used the 
identity
\[\prod_{k=1}^ma_k-\prod_{k=1}^mb_k=\sum_{k=1}^m(\prod_{1\leq l<k}
a_l)(a_k-b_k)\prod_{k<l\leq m}b_l,\]
observing that, in our case, all factors are less than  or equal 
to one and so we can estimate the right hand side by
$\sum_{k=1}^m|a_k-b_k|$, and the differences are
\begin{eqnarray*}
&&\xi^{\zeta}_{x,\kappa ,u}\hat {g}_{y,w}(\kappa^{*},u)+(1-\xi^{\zeta}_{
x,\kappa ,u})g(x,\kappa ,u))-g(x,\kappa ,u)\\
&&\qquad =\xi_{x,\kappa ,u}(\hat {g}_{y,w}(\kappa^{*},u)-g(x,\kappa 
,u))\\
&&\qquad =\xi_{x,\kappa ,u}^{\zeta}(\hat {g}_{y,w}(\kappa^{*},u)-
1)+\xi_{x,\kappa ,u}^{\zeta}(1-g(x,\kappa ,u)).\end{eqnarray*}
Recalling that $g$ vanishes outside $D_{0,\rho_g}$, in the second inequality
we have then used that for $u<u_g$
\begin{equation}
\label{ghatbd}
|1-\hat{g}_{y,w}(\kappa,u)|\leq \frac{|D_{0,\rho_g}\cap D_{y,w}|}{|D_{y,w}|}.
\end{equation}
The sums in the two integrals are over the $(x,\kappa ,u)$ for 
which the term is nonzero.  If there exists $0<c<\infty$ such 
that $E[\eta (D_{y,w})\times {\Bbb K}\times [0,r])]\leq cr|D_{y,w}
|$ for all $y,w,r$, as 
would be the case if $\eta (dx\times {\Bbb K}\times du)$ were a Poisson 
random measure with Lebesgue mean measure,
then the expectation of the right side of (\ref{bound on B2}) is bounded by
\[c|D_{0,\rho_g}|(1+u_g)\int_{[0,1]\times (1,\infty )}v_d(\rho_g+
w)^d\zeta dy\nu^1(w,d\zeta )\nu^2(dw),\]
which is finite under (\ref{flreq1}).

If $m>\rho_g+1$, then
\begin{eqnarray*}
&&\sum_{k=1}^mB_k^1f(\eta )\\
&=&\int_{\Gamma_m^1}f(\eta )\Bigg(\frac {\sum_{S\subset\eta^g_{y,
w}}H(g,\hat {g},S,y,\zeta ,w)}{\prod_{(x,\kappa ,u)\in\eta^g_{y,w}}
g(x,\kappa ,u)}-1\Bigg)dy\nu^1(w,d\zeta )\nu^2(dw)\\
&=&\int_{\Gamma_m^1}f(\eta )\Bigg(\sum_{S\subset\eta^g_{y,w}}\int_{
D_{y,w}^{|S|}}\left(\frac {\prod_{(x,\kappa ,u)\in S}g(x_u,\kappa^{
*}(S),u)}{\prod_{(x,\kappa ,u)\in S}g(x,\kappa ,u)}-1\right)\prod
\upsilon_{y,w}(dx_u)\Bigg)\\
&&\qquad\qquad\qquad\qquad\qquad\qquad\qquad\times\zeta^{|S|}(1-\zeta 
)^{|\eta^g_{y,w}|-|S|}dy\nu^1(w,d\zeta )\nu^2(dw)\\
\\
&=&\int_{\Gamma_m^1}f(\eta )\Bigg(\sum_{(x^{*},\kappa^{*},u^{*})\in
\eta^g_{y,w}}\int_{D_{y,w}}\left(\frac {g(x_{u^{*}},\kappa^{*},u^{
*})}{g(x^{*},\kappa^{*},u^{*})}-1-\frac {(x_{u^{*}}-x^{*})\cdot\nabla 
g(x^{*},\kappa^{*},u^{*})}{g(x^{*},\kappa^{*},u^{*})}\right)\upsilon_{
y,w}(dx_{u^{*}})\Bigg)\\
&&\qquad\qquad\qquad\qquad\qquad\qquad\qquad\times\zeta (1-\zeta 
)^{|\eta^g_{y,w}|-1}dy\nu^1(w,d\zeta )\nu^2(dw)\\
\\
&&\quad +\int_{\Gamma_m^1}f(\eta )\Bigg(\sum_{S\subset\eta^g_{y,w}
,|S|\geq 2}\int_{D_{y,w}^{|S|}}\left(\frac {\prod_{(x,\kappa ,u)\in 
S}g(x_u,\kappa^{*}(S),u)}{\prod_{(x,\kappa ,u)\in S}g(x,\kappa ,u
)}-1\right)\prod\upsilon_{y,w}(dx_u)\Bigg)\\
&&\qquad\qquad\qquad\qquad\qquad\qquad\qquad\times\zeta^{|S|}(1-\zeta 
)^{|\eta^g_{y,w}|-|S|}dy\nu^1(w,d\zeta )\nu^2(dw),\end{eqnarray*}
where in the first term on the right, we are summing 
over $S\subset\eta^g_{y,w}$ with $|S|=1$ and in the second term, we 
are summing over $S\subset\eta^g_{y,w}$ with $|S|\geq 2$. 
If we assume that 
$m>\rho_g+1$, then since in the integral over $\Gamma_m^1$ we have
$w<1$, for each $x^*$ for which $\nabla g(x^*,\kappa^*,u^*)$
is non-trivial we have
\begin{equation}\int_{D_{0,m}\times [2^{-m},1]}{\bf 1}_{D_{y,w}}(
x^{*})\int_{D_{y,w}}(x'-x^{*})\upsilon_{y,w}(dx')dy\nu^2(dw)=0,\label{mzero}\end{equation}
and so including the gradient term has no effect.
Also, observe that the $\nabla g$ term plays the same role here as it 
does in the generator of a L\'evy process (in fact, the 
location of the particle at a fixed level $u$ is a L\'evy 
process). 

Define 
\begin{equation}C_{y,w}g(x,\kappa ,u)=\int_{D_{y,w}}(g(x',\kappa 
,u)-g(x,\kappa ,u)-(x'-x)\cdot\nabla g(s,\kappa ,u))\upsilon_{y,w}
(dx').\label{cyw}\end{equation}
Then, for $m>\rho_g+1$,
\begin{eqnarray*}
|\sum_{k=m+1}^{\infty}B_k^1f(\eta )|&\leq&\int_{\Gamma^1_{\infty}
-\Gamma_m^1}\sum_{S\subset\eta_{y,w}^g,|S|\geq 2}|\prod_{(x,\kappa 
,u)\in S}\hat {g}_{y,w}(\kappa^{*}(S),u)-\prod_{(x,\kappa ,u)\in 
S}g(x,\kappa ,u))|\\
&&\qquad\qquad\qquad\qquad\qquad\qquad\times\zeta^{|S|}(1-\zeta )^{
|\eta^g_{y,w}|-|S|}dy\nu^1(w,d\zeta )\nu^2(dw)\\
&&\quad +\int_{\Gamma^1_{\infty}-\Gamma_m^1}\sum_{(x^{*},\kappa^{
*},u^{*})\in\eta_{y,w}^g}|C_{y,w}g(x^{*},\kappa^{*},u^{*})|\\
&&\qquad\qquad\qquad\qquad\qquad\qquad\times\zeta (1-\zeta )^{|\eta^
g_{y,w}|-1}dy\nu^1(w,d\zeta )\nu^2(dw)\\
&\leq&\int_{\Gamma^1_{\infty}-\Gamma_m^1}(1-(1-\zeta )^{|\eta^g_{
y,w}|}-|\eta^g_{y,w}|\zeta (1-\zeta )^{|\eta^g_{y,w}|-1}){\bf 1}_{
\{D_{0,\rho_g}\cap D_{y,w}\neq\emptyset \}}dy\nu^1(w,d\zeta )\nu^
2(dw)\\
&&\quad +\int_{\Gamma^1_{\infty}-\Gamma_m^1}\Vert\partial^2g\Vert 
|\eta^g_{y,w}|w^2\zeta {\bf 1}_{\{D_{0,\rho_g}\cap D_{y,w}\neq\emptyset 
\}}dy\nu^1(w,d\zeta )\nu^2(dw)\\
&\leq&\int_{\Gamma^1_{\infty}}|\eta_{y,w}^g|(|\eta_{y,w}^g|-1){\bf 1}_{
\{D_{0,\rho_g}\cap D_{y,w}\neq\emptyset \}}\zeta^2dy\nu^1(w,d\zeta 
)\nu^2(dw)\\
&&\quad +\int_{\Gamma^1_{\infty}}\Vert\partial^2g\Vert |\eta^g_{y
,w}|w^2\zeta {\bf 1}_{\{D_{0,\rho_g}\cap D_{y,w}\neq\emptyset \}}
dy\nu^1(w,d\zeta )\nu^2(dw).\end{eqnarray*}

We are primarily interested in solutions $\{\eta_t\}$ of the martingale 
problem for $A^{\infty}$ such that at each time $t$, $\eta_t(\cdot
\times {\Bbb K}\times\cdot )$ is a 
Poisson point process  on 
${\Bbb R}^d\times [0,\infty )$ with mean measure $\ell^{d+1}$.  
Consequently, if, as 
in the discussion of 
$\sum B_k^2$, we 
require that
\begin{equation}E[\eta (D_{y,w}\times {\Bbb K}\times [0,r])]\leq 
cr|D_{y,w}|,\quad\mbox{\rm for all }y\in {\Bbb R}^d,w>0,\label{ord1}\end{equation}
and in addition require
\begin{equation}E[\eta (D_{y,w}\times {\Bbb K}\times [0,r])(\eta 
(D_{y,w}\times {\Bbb K}\times [0,r])-1)]\leq cr^2|D_{y,w}|^2\quad\mbox{\rm for all }
0<w\leq 1,\label{ord2}\end{equation}
the solution of  primary interest will meet these 
requirements.  Under these assumptions
\begin{eqnarray*}
E[|\sum_{k=m+1}^{\infty}B_k^1f(\eta )|]&\leq&cu_g^2v_d(\rho_g+1)^
d\int_{[0,1]\times [0,1]}v_d^2\zeta^2w^{2d}\nu^1(w,d\zeta )\nu^2(
dw)\\
&&\qquad +\Vert\partial^2g\Vert u_g(\rho_g+1)^d\int_{[0,1]\times 
[0,1]}v_dw^{d+2}\zeta\nu^1(w,d\zeta )\nu^2(dw).\end{eqnarray*}
Note that by (\ref{flreq1}) and (\ref{flreq1b}), the right
side is finite.  Comparing the two terms, for 
$d=1$, the first term dominates, while for $d\geq 2$, the 
second term dominates (explaining the need for 
alternative conditions in (\ref{flreq1b})).

From this point on, our approach is reminiscent 
of that in Section~\ref{pure death process}.  For 
$l=1,2,\ldots$,  define $\eta^l_{y,w}=\eta (\cdot\cap D_{y,w}\times 
{\Bbb K}\times [0,l))$.  Set 
\begin{eqnarray*}
\psi_l(\eta )&=&\eta (D_{0,l}\times {\Bbb K}\times [0,l])\int_{[0
,1]\times (1,\infty )}v_d(l+w)^d\zeta\nu^1(w,d\zeta )\nu^2(dw)\\
&&+v_dl^d\int_{{\Bbb R}^d\times [0,1]\times (1,\infty )}\frac {|\eta^
l_{y,w}|}{v_dw^d}{\bf 1}_{\{D_{0,l}\cap D_{y,w}\neq\emptyset \}}\zeta 
dy\nu^1(w,d\zeta )\nu^2(dw)\\
&&+\int_{{\Bbb R}^d\times [0,1]\times [0,1]}|\eta_{y,w}^l|(|\eta^
l_{y,w}|-1){\bf 1}_{\{D_{0,l}\cap D_{y,w}\neq\emptyset \}}dy\nu^1
(w,d\zeta )\nu^2(dw)\\
&&+l\int_{{\Bbb R}^d\times [0,1]\times [0,1]}|\eta^l_{y,w}|w^2{\bf 1}_{
\{D_{0,l}\cap D_{y,w}\neq\emptyset \}}\zeta dy\nu^1(w,d\zeta )\nu^
2(dw).\end{eqnarray*}
Select $\delta_l>0$ so that if (\ref{ord1}) and (\ref{ord2}) are 
satisfied, then $\sum_l\delta_lE[\psi_l(\eta_t)]<\infty$, 
and define
\[\psi (\eta )=1+\sum_l\delta_l\psi_l(\eta ).\]
Then for each $g$ such that 
$f(\eta )=\prod_{(x,\kappa ,u)\in\eta}g(x,\kappa ,u)\in {\cal D}(
A^{\infty})$, there exists $l$ such that 
$\rho_g\leq l$, $u_g\leq l$ and $\Vert\partial^2g\Vert\leq l$, and hence for $
m\geq 0$,
\[|\sum_{k=m+1}^{\infty}B_k^{}|\leq\frac 1{\delta_l}\psi (\eta ).\]
Consequently, we can take $c_f$ in Theorems \ref{mf} and 
Condition \ref{sumcnd2} to be $\delta_l^{-1}$ and $m_f$ in Condition 
\ref{sumcnd2} to be $[\rho_g+2]$.  We have the following.

\begin{theorem}\label{lamfvthm}
Assume that (\ref{flreq1}) and (\ref{flreq1b}) hold and
that $\eta$ is  a solution of the martingale problem 
for $A^{\infty}_{dr,2}$ given by (\ref{vwgen}) satisfying (\ref{ord1})  and (\ref{ord2}).  
Then with the $H_k$ given by (\ref{hklamfv}) and $\lambda_k$ given 
by (\ref{lamklamfv}), the conclusion of Theorem 
\ref{sumeq}\ holds.

Let $\Xi$ be a solution of the martingale problem for $\alpha A_{
dr.2}^{\infty}$ 
given by  (\ref{firstlam}) satisfying
\[E[\int_0^tE[\psi (\eta_s)|\Xi (s)]ds]<\infty ,\]
for all $t>0$, where for each $s$, $\eta_s$ is a conditionally 
Poisson process with Cox measure $\Xi (s)\times\ell$.  Then $\Xi$ 
can be obtained 
from a solution of the martingale problem for $A^{\infty}_{dr,2}$.  In 
particular, the conclusion holds for any solution with
$\Xi (0,dx\times {\Bbb K})$ equal to Lebesgue measure.
\end{theorem}\

\begin{remark}
With the above formulation of  the generator, for finite $\lambda$,
\begin{eqnarray*}
\alpha A^{\lambda}_{dr,2}f(\bar{\eta })=\alpha f(\bar{\eta })\int_{
{\Bbb R}^d\times [0,1]\times [0,\infty )}\sum_{(x,\kappa )\in\bar{
\eta}_{y,w}}\frac {C_{y,w}\bar {g}(x,\kappa )}{\bar {g}(x,\kappa 
)}\zeta (1-\zeta )^{|\bar{\eta}_{y,w}|-1}dy\nu^1(w,d\zeta )\nu^2(
dw)\\
+\alpha f(\bar{\eta })\int_{{\Bbb R}^d\times [0,1]\times [0,\infty 
)}\Bigg(\frac {\sum_{S\subset\bar{\eta}_{y,w},|S|\geq 2}\bar {H}(
\bar {g},\hat {g},S,y,\zeta ,w)}{\prod_{(x,\kappa )\in\bar{\eta}_{
y,w}}\bar {g}(x,\kappa )}-1\Bigg)dy\nu^1(w,d\zeta )\nu^2(dw)\end{eqnarray*}
where $\bar {H}$ is as in (\ref{hbarlam}) for $|S|\geq 2$.

For $\lambda =\infty$, $\alpha A^{\infty}_{dr,2}$ is as in (\ref{firstlam}), that is,
\begin{eqnarray*}
&&\alpha A^{\infty}_{dr,2}f(\Xi )\\
&&\quad =e^{-\int h(x,\kappa )\Xi (dx,d\kappa )}\\
&&\qquad\times\int_{{\Bbb R}^d\times [0,1]\times [0,\infty )}({\Bbb H}_
1(h^{*}_{y,w},\Xi ,y,\zeta ,w)e^{\zeta\int_{D_{y,w}\times {\Bbb K}}
h(x,\kappa )\Xi (dx,d\kappa )}-1)dy\nu^1(w,d\zeta )\nu^2(dw),\end{eqnarray*}
where $h^{*}_{y,w}(\kappa )=\int_0^{\infty}(1-\hat {g}_{y,w}(\kappa 
,u))du$.
\end{remark}

\subsubsection{Stochastic equations for locations and 
types} 
Recall 
\[H_k(\eta ,y,\zeta ,w,(z,v)^{\eta})=\sum_{(x,\kappa ,u)\in\eta}\left
((1-{\bf 1}_{D_{y,w}}(x){\bf 1}_{[0,\zeta ]}(z))\delta_{(x,\kappa 
,u)}+{\bf 1}_{D_{y,w}}(x){\bf 1}_{[0,\zeta ]}(z)\delta_{(y+wv,\kappa_{
u^{*}},u)}\right).\]
Write 
\[\eta (t)=\sum\delta_{(X_u(t),\kappa_u(t),u)}.\]
Under the conditions of Theorem \ref{lamfvthm}, we have 
\begin{eqnarray*}
&&f(\eta (t))=f(\eta (0))\\
&&\quad +\lim_{m\rightarrow\infty}\sum_{k=1}^m\int_0^t\Big(f(H_k(
\eta (s-),Y_k(s-),\zeta_k(s-),W_k(s-),\{(Z_{k,u}(s-),V_{k,u}(s-))
\})\\
&&\qquad\qquad\qquad\qquad\qquad\qquad\qquad\qquad\qquad\qquad -f
(\eta (s-))\Big)dN_k(s),\end{eqnarray*}
where, for each $k$, $N_k$ is a Poisson process with parameter
$\lambda_k$ (as defined in Theorem \ref{lamfvthm}) and
at each jump $\tau$ of $N_k$, if $X_u(\tau -)\in D_{Y_k(\tau 
-),W_k(\tau -)}$, 
then $(Z_{k,u}(\tau -),V_{k,u}(\tau -))$ is replaced by an independent pair 
of random variables $(Z_{k,u}(\tau ),V_{k,u}(\tau ))\in [0,1]\times 
D_{0,1}$  with 
distribution $dz\times\upsilon_{0,1}(dv)$.  Consequently, the location of 
the particle with level $u$ will satisfy 
\begin{eqnarray}
&&X_u(t)=X_u(0)\label{xueq1}\\
&&\quad +\lim_{m\rightarrow\infty}\sum_{k=1}^m\int_0^t{\bf 1}_{D_{
Y_k(s-),W_k(s-)}}(X_u(s-)){\bf 1}_{[0,\zeta_k(s-)]}(Z_{k,u}(s-))(
Y_k(s-)+W_k(s-)V_{k,u}(s-)\non\\
&&\qquad\qquad\qquad\qquad\qquad\qquad\qquad\qquad\qquad\qquad -X_
u(s-))dN_k(s).\non\end{eqnarray}
Note that $\xi_0$ defined by
\[\int_{[0,t]\times {\Bbb R}^d\times [0,1]\times [0,\infty )}f(y,
\zeta ,w)\xi_0(ds,dy,d\zeta ,dw)=\sum_{k=1}^{\infty}\int_0^tf(Y_k
(s-),\zeta_k(s-),W_k(s-))dN_k(s)\]
is a Poisson random measure with mean measure 
$dsdy\nu^1(w,d\zeta )\nu^2(dw)$, and $\xi_u$ defined by 
\begin{eqnarray*}
&&\int_{[0,t]\times \{0,1\}\times D_{0,1}\times {\Bbb R}^d\times 
[0,1]\times [0,\infty )}f(\theta ,v,y,\zeta ,w)\xi_u(ds,d\theta ,
dv,dy,d\zeta ,dw)\\
&&\qquad\qquad\qquad =\sum_{k=1}^{\infty}\int_0^tf({\bf 1}_{[0,\zeta_
k(s-)]}(Z_{k,u}(s-)),V_{k,u}(s-),Y_k(s-),\zeta_k(s-),W_k(s-))dN_k
(s)\end{eqnarray*}
is a Poisson random measure with mean measure
\[ds((1-\zeta )\delta_0(\theta )+\zeta\delta_1(\theta ))\upsilon_{
0,1}(dv)dy\nu^1(w,d\zeta )\nu^2(dw).\]
Then letting $z=(\theta ,v,y,\zeta ,w)$, so 
$\xi_u(ds,d\theta ,dv,dy,d\zeta ,dw)=\xi_u(ds,dz)$, 
(\ref{xueq1}) becomes  
\begin{eqnarray}
X_u(t)&=&X_u(0)\label{loc1}\\
&&\qquad +\lim_{k\rightarrow\infty}\int_{[0,t]\times \{0,1\}\times 
D_{0,1}\times\Gamma_k}{\bf 1}_{D_{y,w}}(X_u(s-))\theta (y+wv-X_u(
s-))\xi_u(ds,dz)\non\\
&=&X_u(0)+\int_{[0,t]\times \{0,1\}\times D_{0,1}\times {\Bbb R}^
d\times [0,1]\times [0,\infty )}{\bf 1}_{D_{y,w}}(X_u(s-))\theta 
(y+wv-X_u(s-))\tilde{\xi}_u(ds,dz),\non\\
\end{eqnarray}
where $\tilde{\xi}_u$ is $\xi_u$ centered by its mean measure.  The 
centering has no effect on the right side of the first 
equality once $k$ is large enough that $X_u(s)\in D_{0,k-1}$ for 
$0\leq s\leq t$.  In particular, if $X_u(s)\in D_{0,k-1}$,
\begin{eqnarray*}
&&\int_{\{0,1\}\times D_{0,1}\times D_{0,k}\times [0,1]\times [2^{
-k},1]}{\bf 1}_{D_{y,w}}(X_u(s))\theta (y+wv-X_u(s))\\
&&\qquad\qquad\qquad\qquad\qquad\qquad ((1-\zeta )\delta_0(\theta 
)+\zeta\delta_1(\theta ))\upsilon_{0,1}(dv)dy\nu^1(w,d\zeta )\nu^
2(dw)\\
&&\qquad =\int_{D_{0,k}\times [0,1]\times [0,1]}{\bf 1}_{D_{0,w}}
(X_u(s)-y)\zeta (y-X_u(s))dy\nu^1(w,d\zeta )\nu^2(dw)\\
&&\qquad =\int_{{\Bbb R}^d\times [0,1]\times [0,1]}{\bf 1}_{D_{0,
w}}(y)\zeta ydy\nu^1(w,d\zeta )\nu^2(dw)\\
&&\qquad =0,\end{eqnarray*}
where in the first equality we have used that $v$ is uniformly distributed on
$D_{0,1}$ and so has mean zero.
Furthermore, (\ref{flreq1}) and (\ref{flreq1b}) imply the 
existence of the stochastic integrals in the limiting 
equation.

Set $R_{y,w,v}(x)=(y+wv-x)(y+wv-x)^T$.
Assuming existence of a solution, the 
centered integral in (\ref{loc1}) is a square integrable 
martingale $M_u$ with covariation matrix
\[[M_u]_t=\int_{[0,t]\times \{0,1\}\times D_{0,1}\times {\Bbb R}^
d\times [0,1]\times [0,1]}{\bf 1}_{D_{y,w}}(X_u(s-))\theta R_{y,w
,v}(X_u(s-))\xi_u(ds,d\theta ,dv,dy,d\zeta ,dw)\]
and, by translation invariance,
\begin{eqnarray*}
E[[M_u]_t]&=&t\int_0^1\int_0^1\int_{D_{0,w}}\int_{D_{0,1}}\zeta (
y+wv)(y+wv)^T\upsilon_{0,1}(dv)dy\nu^1(w,d\zeta )\nu^2(dw)\\
&=&t\int_0^1\int_0^1\int_{D_{0,w}}\zeta (yy^T+|w|^2c_dI)dy\nu^1(w
,d\zeta )\nu^2(dw)\\
&=&tC_d\int_0^1\int_0^1\zeta |w|^{2+d}\nu^1(w,d\zeta )\nu^2(dw)I,\end{eqnarray*}
for appropriate choices of $c_d$ and $C_d$, which is finite by
 (\ref{flreq1b}).  

\begin{lemma}\label{lamexun}
Assume (\ref{flreq1}) and (\ref{flreq1b}).  Then weak 
(distributional)
existence holds for the system (\ref{loc1}).  If, in 
addition, (\ref{flreq2}) holds, then strong uniqueness (and 
hence strong existence) holds.
\end{lemma}

\begin{remark}
Weak uniqueness (uniqueness in distribution) 
for a single $X_u$ follows by uniqueness 
of the corresponding martingale problem ($X_u$ is a L\'evy 
process).  Unfortunately, weak uniqueness for a single 
$X_u$ does not imply weak uniqueness for the system.  If 
we consider the joint distribution of $X_u$ and $X_{u'}$, weak 
uniqueness only implies uniqueness of the marginal 
distributions.  Strong existence means that $X_u$ can be 
written as a function of the stochastic inputs, and 
strong uniqueness implies there is only one such 
function (up to modification on events of probability 
zero).  Strong uniqueness for a single $X_u$ would give 
strong uniqueness for the system.  In the case of $d=1$, 
strong uniqueness is proved in \cite{XZ17} under the 
more general conditions (\ref{flreq1}) and (\ref{flreq1b}) 
\end{remark}

\begin{proof}
It is enough to consider an arbitrary 
but finite subsystem  $\{X_{u_i},1\leq i\leq m\}$. With reference to 
(\ref{cyw}), let $f(x)=\prod_{i=1}^mg(x_i)$, $g\in C^2({\Bbb R}^d
)$, $0\leq g\leq 1$, 
$g(x)=1$ for $x$ outside $D_{0,\rho_g}$,  
and for
$S\subset \{i:x_i\in D_{y,w},1\leq i\leq m\}$, let $f_S(x)=\prod_{
i\in S}g(x_i)$ and
\begin{eqnarray}
B_{y,w}^{|S|}f_S(x)&=&\int_{D_{y,w}^{|S|}}\bigg(\prod_{i\in S}g(x_
i')-\prod_{i\in S}g(x_i)\label{byw2}\\
&&\qquad -{\bf 1}_{\{w\leq 1\}}f_S(x)\sum_{i\in S}\frac {(x_i'-x_
i)\cdot\nabla g(x_i)}{g(x_i)}\bigg)\prod_{i\in S}\upsilon_{y,w}(d
x_i').\non\end{eqnarray}
Then setting $S_{y,w}(x)=\{i:x_i\in D_{y,w}\}$, 
the generator for  the subsystem becomes 
\begin{eqnarray}
A^mf(x)&=&f(x)\int_{{\Bbb R}^d\times [0,1]\times [0,\infty )}\sum_{
S\subset S_{y,w}(x)}\frac {B_{y,w}^{|S|}f_S(x)}{\prod_{i\in S}g(x_
i))}\label{fingen}\\
&&\qquad\qquad\qquad\qquad\qquad\qquad\times\zeta^{|S|}(1-\zeta )^{
|S_{y,w}(x)|-|S|}dy\nu^1(w,d\zeta )\nu^2(dw).\non\end{eqnarray}
Note that $A^mf(x)$ is a continuous function of $x$.

Existence of solutions of the martingale problem for 
(\ref{fingen}) follows 
by approximation. 
To obtain an approximation $X^{\epsilon}=(X_1^{\epsilon},\cdots ,
X_m^{\epsilon})$, 
consider the system obtained 
by replacing $\nu^2$ by $\nu^2_{\epsilon}$ given by $\nu^2_{\epsilon}
(C)=\nu^2(C\cap [\epsilon ,\infty ))$.  
 The generator $A^{m,\epsilon}$ is then a bounded operator (the 
gradient term integrates to zero), and existence and 
uniqueness for the martingale problem is immediate. 
For each $i$, $X^{\epsilon}_i$ is  a L\'evy process with L\'evy measure
\begin{eqnarray*}
\nu_{\epsilon}(C)&=&\int_{{\Bbb R}^d\times [0,1]\times [\epsilon 
,\infty )}\int_{D_{y,w}}{\bf 1}_C(x'-x){\bf 1}_{D_{y,w}}(x)\upsilon_{
y,w}(dx')\zeta dy\nu^1(w,d\zeta )\nu^2(dw)\\
&=&\int_{{\Bbb R}^d\times [0,1]\times [\epsilon ,\infty )}\int_{D_{
0,w}}{\bf 1}_C(x'+y-x){\bf 1}_{D_{0,w}}(x-y)\upsilon_{0,w}(dx')\zeta 
dy\nu^1(w,d\zeta )\nu^2(dw)\\
&=&\int_{{\Bbb R}^d\times [0,1]\times [\epsilon ,\infty )}\int_{D_{
0,w}}{\bf 1}_C(x'-z){\bf 1}_{D_{0,w}}(z)\upsilon_{0,w}(dx')\zeta 
dz\nu^1(w,d\zeta )\nu^2(dw),\end{eqnarray*}
and convergence in distribution of $\{X^{\epsilon}_i\}$ follows from 
convergence of the L\'evy measures.  Convergence for 
each component implies relative compactness of $\{X^{\epsilon}\}$ at 
least in $D_{{\Bbb R}^d}[0,\infty )\times\cdots\times D_{{\Bbb R}^
d}[0,\infty )$ (\cite{EK86}, 
Proposition 3.2.4), if not in $D_{({\Bbb R}^d)^m}[0,\infty )$. For a 
convergent subsequence, convergence in 
the product topology still implies convergence of the 
integrals
\[\int_0^tA^{m,\epsilon}f(X^{\epsilon}(s))ds\Rightarrow\int_0^tA^
mf(X(s))ds,\]
which in turn ensures that the limit is a solution of the 
martingale problem for $A^m$.  The fact that the limit is a 
weak solution of the stochastic differential equation 
follows by Theorem 2.3 of \cite{Kur11}.

If (\ref{flreq2}) holds, then a solution of (\ref{loc1}) jumps 
only finitely often in a finite time interval, that is, 
$\{(s,z)\in\xi_u:s\leq t,X_u(s-)\in D_{y,w},\theta =1\}$ is finite for each $
t>0$. 
Consequently, the equation is uniquely solved by moving 
from one such $(s,z)$ to the next, and this solution 
depends only on the stochastic inputs, that is, it is a 
strong solution.
\end{proof}

We still need to consider the evolution of the type 
of each particle.  Note that the particle with index $u$ changes 
type only if it is involved in a birth/death event with a 
particle having a lower level.  The number of times that 
particle $u_1$ and particle $u_2$ are involved in the same 
birth/death 
event up to time $t$ can be written as
\[N_{u_1u_2}(t)=\int_{[0,t]\times {\Bbb R}^d\times [0,1]\times [0
,\infty )}{\bf 1}_{D_{y,w}}(X_{u_2}(s-)){\bf 1}_{D_{y,w}}(X_{u_1}
(s-))\theta_{u_1}(s)\theta_{u_2}(s)\xi_0(ds,dy,d\zeta ,dw)\]
and  since $\theta_{u_1}$ and $\theta_{u_2}$ are conditionally independent given 
$\xi_0$,
\[E[N_{u_1u_2}(t)]=\int_{[0,t]\times {\Bbb R}^d\times [0,1]\times 
[0,\infty )}E[{\bf 1}_{D_{y,w}}(X_{u_2}(s)){\bf 1}_{D_{y,w}}(X_{u_
1}(s))]\zeta^2dsdy\nu^1(w,d\zeta )\nu^2(dw).\]
Let $C\subset {\Bbb R}^d$ be bounded and $u>0$, and let $N_{C,u}(
t)$ be the 
number of times by time $t$ that two particles with 
levels below $u$  and locations in 
$C$ are involved in the same birth/death event.  Then, 
assuming (\ref{ord1}) and (\ref{ord2}), 
\begin{eqnarray*}
E[N_{C,u}(t)]&=&\int_{[0,t]\times {\Bbb R}^d\times [0,1]\times [0
,\infty )}E[\sum_{u_1<u_2\leq u}{\bf 1}_C(X_{u_1}(s)){\bf 1}_C(X_{
u_2}(s)){\bf 1}_{D_{y,w}}(X_{u_2}(s)){\bf 1}_{D_{y,w}}(X_{u_1}(s)
)]\\
&&\qquad\qquad\qquad\qquad\qquad\qquad\qquad\qquad\qquad\qquad\times
\zeta^2dsdy\nu^1(w,d\zeta )\nu^2(dw)\\
&=&\int_{[0,t]\times {\Bbb R}^d\times [0,1]\times [0,\infty )}E[{{
\eta (s,D_{y,w}\cap C\times [0,u])}\choose 2}]\zeta^2dsdy\nu^1(w,
d\zeta )\nu^2(dw)\\
&\leq&\int_{[0,t]\times {\Bbb R}^d\times [0,1]\times [0,\infty )}
u^2|D_{y,w}\cap C|^2\zeta^2dsdy\nu^1(w,d\zeta )\nu^2(dw)\\
&\leq&u^2\int_{[0,t]\times {\Bbb R}^d\times [0,1]\times [0,\infty 
)}(v_d^2w^{2d}\wedge |C|^2)\zeta^2dsdy\nu^1(w,d\zeta )\nu^2(dw)\\
&<&\infty .\end{eqnarray*}
It follows that no single particle will change type more 
than finitely often by time $t$. 

It is now straightforward to write down an equation for 
the way in which individuals' types change with time.  
For $u_1<u_2$, define 
\[L_{u_1u_2}(t)=\#\{s\leq t:N_{u_1u_2}(s)-N_{u_1u_2}(s-)=1,N_{u_3
u_2}(s)-N_{u_3u_2}(s-)=0,\;\forall u_3<u_1\}.\]
Then, writing $\kappa_u$ for the type of the individual with level 
$u$,
\begin{equation}\kappa_{u_2}(t)=\kappa_{u_2}(0)+\sum_{u_1<u_2}\int_
0^t(\kappa_{u_1}(s-)-\kappa_{u_2}(s-))dL_{u_1u_2}(s).\label{typ1}\end{equation}

As in Section 5 of  \cite{DK99a}, the genealogy of the 
particles alive at time $t$ is determined by the $L_{u_1u_2}$.  In 
particular, the index of the ancestor at time $r<t$ of the 
particle at level $u_2$ at time $t$ satisfies 
\begin{equation}J_{u_2}(t,r)=u_2-\sum_{u_3<u_1\leq u_2}\int_r^t(u_
1{\bf 1}_{\{J_{u_2}(t,s)=u_1\}}-u_3)dL_{u_3u_1}(s).\label{anlev}\end{equation}
Since the lookdown construction for the
discrete population model is simply the 
restriction of the lookdown construction of the infinite 
density population model, the genealogies of the discrete 
model converge to those of the infinite 
density model. To be precise:

\begin{theorem}
For any solution of the infinite system (\ref{loc1}) 
(regardless of the uniqueness question), the counting 
processes $N_{u_1u_2}$ and $L_{u_1u_2}$, the type processes $\kappa_
u$, and 
the ancestral index $J_u$ are uniquely determined, and the 
genealogies of the $\lambda <\infty$ model converge to those of 
the $\lambda =\infty$ model.  
\end{theorem}

\begin{remark}\label{geneaeq}
As noted in Section \ref{sectransam}, one can model the 
selection of a random sample of size $n$ from a region $C$ 
satisfying $0<\Xi (t,C\times {\Bbb K})<\infty$ simply by selecting the particles 
located in $C$ with the $n$ lowest levels.  The genealogy of 
the sample can then be obtained using equation 
(\ref{anlev}).
\end{remark}

\subsubsection{Second construction of spatial 
$\Lambda$-Fleming-Viot with levels} \label{sectfv2} The particle 
dynamics for $A^{\lambda}_{dr,2}$ given by (\ref{loc1}) and (\ref{typ1}) 
are very different from the particle dynamics that are 
natural for $A^{\lambda}_{th,db}$ defined in (\ref{lamf12}).  The event 
measures $\mu (dz)\equiv dy\nu^1(w,d\zeta )\nu^2(dw)$ are of the same form, 
but what happens at each event $z=(y,\zeta ,w)$ is very 
different.  In particular, for a birth/death event in the 
ball $D_{y,w}$, the total population size in $D_{y,w}$ does not 
change for $A_{dr,2}^{\lambda}$, but typically it will change for $
A_{th,db}^{\lambda}$.  
As will become apparent when we analyze the behavior 
of the levels, we will need to assume stronger 
conditions on the event measures than were used in the 
previous construction.  

With ${\cal D}_{\lambda}$ defined in (\ref{scrddef}), 
we take the domain of $A^{\lambda}\equiv A^{\lambda}_{th,db}$ to be
\[\{f(\eta )=\prod_{(x,\kappa ,u)\in\eta}g(x,\kappa ,u)\in {\cal D}_{
\lambda}:g(\cdot ,\kappa ,\cdot )\in C^2,\Vert\partial_ug\Vert\equiv\sup_{
x,\kappa ,u}|\partial_ug(x,\kappa ,u)|<\infty \},\]
and ${\cal D}(A^{\infty})=\cup_{\lambda}{\cal D}(A^{\lambda})$.
In a birth/death event determined by $z=(y,\zeta ,w)$, the 
parent is killed and, with 
probability $1$, for $\lambda <\infty$
all other particles in the event region $D_{y,w}$ 
will change levels and for $\lambda =\infty$, all particles with 
levels above that of the parent will change levels.  

For finite $\lambda$, $v^{*}$ has density 
$(1-e^{-\lambda\alpha_z})^{-1}\alpha_ze^{-\alpha_zv}$ on $[0,\lambda 
]$, where $\alpha_z=\zeta |D_{y,w}|$.  Let 
$u_1^{*}=\min\{u:u>v^{*},(x,\kappa ,u)\in\eta_{|D_{y,w}}\}$ and 
$u_2^{*}=\max\{u:u<v^{*},(x,\kappa ,u)\in\eta_{|D_{y,w}}\}$ and define
\[\tau_1^{*}=\log\frac {\lambda -v^{*}}{\lambda -u_1^{*}},\quad\tau_
2^{*}=\log\frac {\lambda}{u^{*}_2},\mbox{\rm \ and }\tau^{*}=\tau_
1^{*}\wedge\tau_2^{*}.\]
Setting
\[u^{*}=u_1^{*}{\bf 1}_{\{\tau_1^{*}<\tau_2^{*}\}}+u_2^{*}{\bf 1}_{
\{\tau_1^{*}>\tau_2^{*}\}},\]
for $(x,\kappa ,u)\in\eta_{|D_{y,w}}$, we have
\begin{equation}{\cal J}_{y,w}^{\lambda}(x,u,\eta ,v^{*})={\bf 1}_{
\{u>v^{*}\}}(ue^{\tau^{*}}-\lambda (e^{\tau^{*}}-1))+{\bf 1}_{\{u
<v^{*}\}}ue^{\tau^{*}},\label{levtran2}\end{equation}
and for $\lambda =\infty$, 
\begin{equation}{\cal J}^{\infty}_{y,w}(x,u,\eta ,v^{*})=u-{\bf 1}_{
[v^{*},\infty )}(u)(u^{*}-v^{*}).\label{levtran3}\end{equation}
If $(x,\kappa ,u)\in\eta_{|D_{y,w}}$ is not the parent, that is, $
u\neq u^{*}$, then 
$(x,\kappa ,u)$ jumps to $(x,\kappa ,\frac {{\cal J}^{\lambda}_{y
,w}(x,u,\eta ,v^{*})}{1-\zeta})$. 

For reasons that will become clear below, we also require the 
stronger condition (\ref{flreq2}), that is, 
\begin{equation}\int_{(0,\infty )\times [0,1]}\zeta w^d\nu^1(w,d\zeta 
)\nu^2(dw)<\infty .\label{flreq22}\end{equation}

Recall (\ref{ghatbd}), and note that at an event 
$z=(y,\zeta ,w)$, 
the expected number of new particles with level 
below $u_g$ is bounded by 
\[\frac 1{1-e^{-\lambda\alpha_z}}\zeta |D_{y,w}|u_g.\]
Setting ${\Bbb U}={\Bbb R}^d\times [0,1]\times [0,\infty )$ and assuming (\ref{intbnd}), 
define $\Gamma_k$ as in (\ref{gamk}).  With reference to Appendix 
\ref{sumbnd}, define
\[B_k^{\lambda}f(\eta )=\int_{\Gamma_k-\Gamma_{k-1}}(H_z^{\lambda}
(g,\eta )-f(\eta ))(1-e^{-\lambda\alpha_z})dy\nu^1(w,d\zeta )\nu^
2(dw),\]
where $H_z^{\lambda}(g,\eta )$ is defined in (\ref{hzlamdef}).
Note that $u^{*}$ and $\kappa^{*}$ are
determined by  $v^{*}$ and $\eta_{y,w}$.   Then

\begin{eqnarray}
|B_k^{\lambda}f(\eta )|&\leq&\int_{\Gamma_k-\Gamma_{k-1}}|H_z^{\lambda}
(g,\eta )-f(\eta )|(1-e^{-\alpha_z\lambda})dy\nu^1(w,d\zeta )\nu^
2(dw)\label{genbnd3}\\
&\leq&\int_{\Gamma_k-\Gamma_{k-1}}\int_0^{\lambda}\alpha_ze^{-\alpha_
zv^{*}}\left(1-\hat g_{y,w}(\kappa^{*},v^{*})e^{-\alpha_z\int_{v^{
*}}^{\lambda}(1-\hat {g}_{y,w}(\kappa^{*},v))dv}\right)dv^{*}\mu 
(dz)\non\\
&&\qquad +\int_{\Gamma_k-\Gamma_{k-1}}\int_0^{\lambda}\alpha_ze^{
-\alpha_zv^{*}}\bigg|\prod_{(x,\kappa ,u)\in\eta ,x\in D_{y,w},u\neq 
u^{*}}g(x,\kappa ,\frac 1{1-\zeta}{\cal J}_{y,w}^{\lambda}(x,u,\eta 
,v^{*}))\non\\
&&\qquad\qquad\qquad\qquad\qquad\qquad\qquad\qquad\qquad -\prod_{
(x,\kappa ,u)\in\eta ,x\in D_{y,w}}g(x,\kappa ,u)\bigg|dv^{*}\mu 
(dz).\non\end{eqnarray}
Note that the integrand in the first term on the right is 
zero if $v^{*}\geq u_g$ and
\[\hat {g}_{y,w}(\kappa^{*},v)\geq 1-\frac {|D_{y,w}\cap D_{0,\rho_
g}|}{|D_{y,w}|}\equiv\underline {\hat {g}}_{y,w}\]
Then, bounding the two terms on the right of 
(\ref{genbnd3}),
\begin{eqnarray*}
|B_k^{\lambda}f(\eta )|&\leq&\int_{\Gamma_k-\Gamma_{k-1}}\int_0^{
u_g}\alpha_ze^{-\alpha_zv^{*}}\left(1-\underline {\hat {g}}_{y,w}
e^{-\alpha_z(u_g-v^{*})(1-\underline {\hat {g}}_{y,w})}\right)dv^{
*}\mu (dz)\non\\
&&\qquad +\int_{\Gamma_k-\Gamma_{k-1}}(1-e^{-\alpha_z\lambda})\bar{
\eta }(D_{y,w}\cap D_{0,\rho_g}\times {\Bbb K})\mu (dz)\non\\
&\leq&\int_{\Gamma_k-\Gamma_{k-1}}\left((1-\underline {\hat {g}}_{
y,w})(1-e^{-\alpha_zu_g})+\int_0^{u_g}\alpha_ze^{-\alpha_zv^{*}}\underline {
\hat {g}}_{y,w}\left(1-e^{-\alpha_z(u_g-v^{*})(1-\underline {\hat {
g}}_{y,w})}\right)dv^{*}\right)\mu (dz)\non\qquad\\
&&\qquad +\int_{\Gamma_k-\Gamma_{k-1}}(1-e^{-\alpha_z\lambda})\bar{
\eta }(D_{y,w}\cap D_{0,\rho_g}\times {\Bbb K})\mu (dz)\non\\
&\leq&\int_{\Gamma_k-\Gamma_{k-1}}2\alpha_zu_g(1-\underline {\hat {
g}}_{y,w})\mu (dz)\non\qquad\\
&&\qquad +\int_{\Gamma_k-\Gamma_{k-1}}(1-e^{-\alpha_z\lambda})\bar{
\eta }(D_{y,w}\cap D_{0,\rho_g}\times {\Bbb K})\mu (dz)\non\\
&\leq&\int_{\Gamma_k-\Gamma_{k-1}}2u_g\zeta |D_{y,w}|\frac {|D_{y
,w}\cap D_{0,\rho_g}|}{|D_{y,w}|}\mu (dz)\non\\
&&\qquad +\int_{\Gamma_k-\Gamma_{k-1}}\zeta |D_{y,w}|\lambda\bar{
\eta }(D_{y,w}\cap D_{0,\rho_g}\times {\Bbb K})\mu (dz).\non\end{eqnarray*}

The construction of the $\psi$ needed to apply Theorem 
\ref{mf} and Theorem \ref{sumeq}
is similar to the construction in the previous 
section.   Bounding the parameters in the estimates 
above that depend on $g$ by a positive integer $l$, we have 
\begin{eqnarray}
\sum_{k=1}^{\infty}|B_k^{\lambda}f(\eta )|&\leq&\int_{[0,1]\times 
[0,\infty )}2v_d\zeta l(w\wedge l)^d\nu^1(w,d\zeta )\nu^2(dw)\non\\
&&\qquad +\int_{{\Bbb U}}\zeta |D_{y,w}|\lambda\bar{\eta }(D_{y,w}
\cap D_{0,l}\times {\Bbb K})\mu (dz)\non\\
&\equiv&\psi_l(\eta ),\end{eqnarray}
provided $u_g$ and $\rho_g$ are less than $l$.

We are primarily interested in solutions of the 
martingale problem for which $\eta (\cdot\times {\Bbb K}\times\cdot 
)$ will be 
dominated by a 
Poisson random measure on ${\Bbb R}^d\times [0,\lambda ]$ with Lebesgue 
mean measure, so  restricting our attention to 
solutions of the martingale problem satisfying
\begin{equation}E[\bar{\eta }(D_{y,w}\times {\Bbb K})]\leq c|D_{y
,w}|=cv_dw^d,\label{mombnd2}\end{equation}
we have 
\[E[\psi_l(\eta )]\leq\int_{[0,1]\times [0,\infty )}(v_d\zeta l(w
\wedge l)^d+c\lambda v_d^2\zeta w^d(w\wedge l)^d)\nu^1(w,d\zeta )\nu^
2(dw),\]
which is finite under (\ref{flreq22}). Then, as before, we 
set 
 $\psi (\eta )=1+\sum_{l=1}^{\infty}\delta_l\psi_l(\eta )$, where we select $
\delta_l>0$  satisfying
\[\sum_{l=1}^{\infty}\delta_l\int_{[0,1]\times [0,\infty )}(v_d\zeta 
l(w\wedge l)^d+c\lambda v_d^2\zeta w^d(w\wedge l)^d)\nu^1(w,d\zeta 
)\nu^2(dw)<\infty .\]

The $\lambda =\infty$ case takes a little more care.  Note that we 
exploit the fact that if $u<v^{*}$ and $x\in D_{y,w}$, then
\[{\cal J}^{\infty}_{y,w}(x,u,\eta ,v^{*})=u,\]
so after a birth-death event, the new level is $\frac 1{1-\zeta}u$.
Let 
${\Bbb U}_1=\{(y,\zeta ,w)\in {\Bbb U}:\zeta\leq\frac 12\}$ and $
{\Bbb U}_2=\{(y,\zeta ,w)\in {\Bbb U}:\zeta >\frac 12\}$.
\begin{eqnarray}
|B_kf(\eta )|&\leq&\int_{\Gamma_k-\Gamma_{k-1}}|H_z(g,\eta )-f(\eta 
)|dy\nu^1(w,d\zeta )\nu^2(dw)\label{genbnd4}\\
&\leq&\int_{\Gamma_k-\Gamma_{k-1}}\int_0^{\infty}\alpha_ze^{-\alpha_
zv^{*}}\left(1-\hat g_{y,w}(\kappa^{*},v^{*})e^{-\alpha_z\int_{v^{
*}}^{\infty}(1-\hat {g}_{y,w}(\kappa^{*},v))dv}\right)dv^{*}\mu (
dz)\non\\
&&\qquad +\int_{\Gamma_k-\Gamma_{k-1}}\int_0^{\infty}\alpha_ze^{-
\alpha_zv^{*}}\bigg|\prod_{(x,\kappa ,u)\in\eta ,x\in D_{y,w},u\neq 
u^{*}}g(x,\kappa ,\frac 1{1-\zeta}{\cal J}^{\infty}_{y,w}(x,u,\eta 
,v^{*}))\non\\
&&\qquad\qquad\qquad\qquad\qquad\qquad\qquad\qquad\qquad -\prod_{
(x,\kappa ,u)\in\eta ,x\in D_{y,w}}g(x,\kappa ,u)\bigg|dv^{*}\mu 
(dz)\non\\
&\leq&\int_{\Gamma_k-\Gamma_{k-1}}2u_g\zeta |D_{y,w}|\frac {|D_{y
,w}\cap D_{0,\rho_g}|}{|D_{y,w}|}\mu (dz)\non\qquad\qquad\\
&&\qquad +\int_{(\Gamma_k-\Gamma_{k-1})}\int_0^{u_g}\alpha_ze^{-\alpha_
zv^{*}}\eta (D_{y,w}\cap D_{0,\rho_g}\times {\Bbb K}\times [0,u_g
+u_1^{*}-v^{*}))dv^{*}\mu (dz)\non\\
&&\qquad +\int_{(\Gamma_k-\Gamma_{k-1})\cap {\Bbb U}_1}e^{-\alpha_
zu_g}\eta (D_{y,w}\cap D_{0,\rho_g}\times {\Bbb K}\times [0,u_g])
\Vert\partial_ug\Vert\frac {\zeta}{1-\zeta}u_g\mu (dz)\non\\
&&\qquad +\int_{(\Gamma_k-\Gamma_{k-1})\cap {\Bbb U}_2}e^{-\alpha_
zu_g}\eta (D_{y,w}\cap D_{0,\rho_g}\times {\Bbb K}\times [0,u_g])
\mu (dz).\non\end{eqnarray}
The first term corresponds to offspring of the event, the second 
accounts for the change in levels of individuals already present in
the population in
the case $v^*<u_g$ and the final two terms 
to the corresponding changes when $v^*>u_g$. As in Section~\ref{firstslfv},
we are bounding the difference of two products in which all the factors
are less than or equal to one, by a sum of differences
of factors.  

As before, for $u_g$, $\Vert\partial_ug\Vert$, and $\rho_g$ less than $
l$, 
\begin{eqnarray}
\sum_{k=1}^{\infty}|B_kf(\eta )|&\leq&\int_{{\Bbb U}}2l\zeta |D_{
y,w}|\frac {|D_{y,w}\cap D_{0,l}|}{|D_{y,w}|}\mu (dz)\non\qquad\qquad\\
&&\qquad +\int_{{\Bbb U}}\int_0^l\alpha_ze^{-\alpha_zv^{*}}\eta (
D_{y,w}\cap D_{0,l}\times {\Bbb K}\times [0,l+u_1^{*}-v^{*}))dv^{
*}\mu (dz)\non\\
&&\qquad +\int_{{\Bbb U}_1}e^{-\alpha_zl}\eta (D_{y,w}\cap D_{0,l}
\times {\Bbb K}\times [0,l])\frac {\zeta}{1-\zeta}l^2\mu (dz)\non\\
&&\qquad +\int_{{\Bbb U}_2}e^{-\alpha_zl}\eta (D_{y,w}\cap D_{0,l}
\times {\Bbb K}\times [0,l])\mu (dz)\non\\
&\equiv&\psi_l(\eta )\non\end{eqnarray}

Note that in the second term on the right, $v^{*}\leq l$, and 
$u_1^{*}-v^{*}<u_1^{*}\leq 2l$, if $\eta (D_{y,w}\times {\Bbb K}\times 
(l,2l])>0$.  In 
general we have
\begin{eqnarray*}
&&\eta (D_{y,w}\cap D_{0,l}\times {\Bbb K}\times [0,l+u_1^{*}-v^{*}
))\leq\eta (D_{y,w}\cap D_{0,l}\times {\Bbb K}\times [0,3l))\\
&&\quad+\sum_{k=2}^{\infty}{\bf 1}_{\{\eta (D_{y,w}\times {\Bbb K}\times 
(l,kl])=0,\eta (D_{y,w}\times {\Bbb K}\times (kl,(k+1)l])>0\}}\eta 
(D_{y,w}\cap D_{0,l}\times {\Bbb K}\times ((k+1)l,(k+2)l)),\end{eqnarray*}
and assuming $\eta$ is conditionally Poisson with Cox 
measure $\Xi (dx,d\kappa )du$, the conditional independence of 
$\eta$ on disjoint sets gives
\begin{eqnarray*}
&&E[\eta (D_{y,w}\cap D_{0,l}\times {\Bbb K}\times [0,l+u_1^{*}-v^{
*}))]\\
&&\quad\leq E[\Xi (D_{y,w}\cap D_{0,l}\times {\Bbb K})]3l\\
&&\qquad\qquad +\sum_{k=2}^{\infty}E[e^{-\Xi (D_{y,w}\times {\Bbb K}
)(k-1)l}(1-e^{-\Xi (D_{y,w}\times {\Bbb K})l})\Xi (D_{y,w}\cap D_{
0,l}\times {\Bbb K})]2l\\
&&\quad\leq 5lE[\Xi (D_{y,w}\cap D_{0,l}\times {\Bbb K})].\end{eqnarray*}
Consequently, if there exists $c>0$ such that 
\begin{equation}E[\Xi (D_{y,w}\times {\Bbb K})]\leq c|D_{y,w}|,\label{cmbnd}\end{equation}
 then $E[\psi_l(\eta )]<\infty$, and the conclusions of Theorem 
\ref{mf} and Theorem \ref{sumeq} hold. 

For solutions of the martingale problem for $A^{\lambda}_{th,db}$ or 
$A^{\infty}_{th,db}$, the initial level of each particle will be distinct, 
and we will index particles by their initial level.  Each 
particle has birth time $b_u$, 
which we will take to be $0$ for the particles in the 
population at time $0$, and an initial location $x_u=X_u(b_u)$  
and 
a type $\kappa_u$ which do not change with time.  

Let ${\cal N}={\cal N}({\Bbb R}^d\times [0,\infty ))$ be the space of counting measures on 
${\Bbb R}^d\times [0,\infty )$.
The evolution of the process is determined by a Poisson 
random measure $\xi$ on $[0,\infty )\times {\cal N}\times {\Bbb R}^
d\times [0,1]\times [0,\infty )$ 
with mean measure 
\[ds\nu^3(y,\zeta ,w,d\gamma )dy\nu^1(w,d\zeta )\nu^2(dw),\]
where $\nu^3(y,\zeta ,w,d\gamma )$ is the distribution of the Poisson random 
measure on ${\Bbb R}^d\times [0,\infty )$ with mean measure 
\[\zeta {\bf 1}_{D_{y,w}}(x)dxdv=\zeta |D_{y,w}|\upsilon_{y,w}(dx
)dv.\]
Note that a ``point'' in $\xi$ is of the form 
$\beta =(s,\{(x_k,v_k),k\geq 1\},y,\zeta ,w)$, where we will assume that the 
$\{(x_k,v_k)\}$ are indexed in increasing order of the $v_k$. 
Then, the 
birth times, locations, and levels of ``new'' particles are given by 
\[{\cal B}=\cup_{\beta\in\xi}\{(s,x_k,v_k),k\geq 1,v_k<\lambda \}
.\]

Then $v^{*}\equiv v(\beta )=v_1$ and $x(\beta )=x_1$. Then
$(x^{*},\kappa^{*},u^{*})\equiv (x^{*}(\beta ,\eta ),\kappa^{*}(\beta 
,\eta ),u^{*}(\beta ,\eta ))$ is the point in $\eta$ 
satisfying $x^{*}\in D_{y,w}$ and 
\[u^{*}=\mbox{\rm argmax}\{\frac {\lambda -u}{\lambda -v^{*}}:(x,
\kappa ,u)\in\eta ,x\in D_{y,w},u\geq v^{*}\}\cup \{\frac u{v^{*}}
:(x,\kappa ,u)\in\eta ,x\in D_{y,w},u\leq v^{*}\}.\]

Set $G={\cal N}\times {\Bbb R}^d\times [0,1]\times [0,\infty )$.  By Theorem \ref{sumeq}, 
we have 

\begin{theorem}\label{eqvar1} 
For $0<\lambda\leq\infty$, any solution of the martingale problem 
for $A^{\lambda}=A^{\lambda}_{th,db}$ given in (\ref{lamf12}) that satisfies
\[E[\int_0^t\psi (\eta (s))ds]<\infty ,\quad\mbox{\rm for all }t\geq 
0,\]
can be obtained as a solution of the stochastic equation
\begin{eqnarray*}
f(\eta (t))&=&f(\eta (0))\\
&&\quad +\int_{[0,t]\times G}(\frac {f(\eta (s-))}{f(\eta_{y,w}(s
-))}\prod_{(x,u)\in\gamma (\beta )}g(x,\kappa^{*}(\beta ,\eta (s-
)),u)\\
&&\qquad\qquad\qquad\times\prod_{(x,\kappa ,u)\in\eta_{y,w}(s-),u
\neq u^{*}(\beta ,\eta (s-))}g(x,\kappa ,\frac {{\cal J}^{\lambda}_{
y,w}(x,u,\eta (s-),v(\beta ))}{1-\zeta})\\
&&\qquad\qquad\qquad\qquad\qquad\qquad\qquad\qquad\qquad\qquad -f
(\eta (s-))){\bf 1}_{\{v(\beta )<\lambda \}}\xi (ds,d\beta ).\end{eqnarray*}

\end{theorem}

To construct a more useful system of equations, 
if $(x,\kappa ,u)\equiv (x_u,\kappa_u,u)\in\eta (0)$, the level evolves by 
\begin{eqnarray}
U_u(t)&=&u+\int_{(0,t]\times G}{\bf 1}_{D_{y,w}}(x_u)(\frac {{\cal J}_{
y,w}^{\lambda}(x_u,U_u(s-),\eta (s-),v(\beta ))}{1-\zeta}-U_u(s-)
)\nonumber\\
&&\qquad\qquad\qquad\qquad\qquad\qquad\qquad {\bf 1}_{\{v(\beta )
<\lambda \}}\xi (ds,d\beta ),\label{equ1}\end{eqnarray}
and the particle dies at time  
\begin{equation}d_u=\inf\{t>0:U_u(t)>\lambda\mbox{\rm \ or }U_u(t
-)=u^{*}(\beta ,\eta (t-)),(t,\beta )\in\xi \}.\label{dthu}\end{equation}
If there is 
a birth/death event at time $s$,
\[(s,\beta )=(s,\{(x_k,v_k),k\geq 1\},y,\zeta ,w)\in\xi ,\]
then for $u=v_k$, we set $x_u=x_k$ and $b_u=s$.
The levels for the new 
particles satisfy
\begin{eqnarray}
U_u(t)&=&u+\int_{(b_u,t]\times G}{\bf 1}_{D_{y,w}}(x_u)(\frac {{\cal J}_{
y,w}^{\lambda}(x_u,U_u(s-),\eta (s-),v(\beta ))}{1-\zeta}-U_u(s-)
)\nonumber\\
&&\qquad\qquad\qquad\qquad\qquad\qquad\qquad {\bf 1}_{\{v(\beta )
<\lambda \}}\xi (ds,d\beta ),\label{equ2}\end{eqnarray}
for $t\geq b_u$, and the type is given by $\kappa_u=\kappa^{*}(\beta 
,\eta (s-))$.  
Again, the particle dies at time $d_u$ given by (\ref{dthu}), 
so
\[\eta (t)=\sum {\bf 1}_{[b_u,d_u)}(t)\delta_{(x_u,\kappa_u,U_u(t
))}.\]

With reference to (\ref{hlim}),
passing to the limit as $\lambda\rightarrow\infty$, the equations become 
\begin{eqnarray*}
U_u(t)&=&u+\int_{(b_u,t]\times G}{\bf 1}_{D_{y,w}}(x_u)(\frac {U_
u(s-)-{\bf 1}_{\{U_u(s-)\geq v(\beta )\}}(u^{*}(\beta ,\eta (s-))
-v(\beta ))}{1-\zeta}-U_u(s-))\xi (ds,d\beta ),\end{eqnarray*}
for $t\geq b_u$, and defining 
\begin{equation}\tau_u=\lim_{k\rightarrow\infty}\inf\{t:U_u(t)>k\}
,\label{hitinf}\end{equation}
the particle dies at time  
\begin{equation}d_u=\tau_u\wedge\inf\{t>0:U_u(t-)=u^{*}(\beta ,\eta 
(t-)),(t,\beta )\in\xi \}.\label{dthu2}\end{equation}

Since the downward jumps in $U_u$, when they occur, will 
typically be $O(1)$, we can only allow finitely many per 
unit time.  Conditional on $U_u$, the intensity of downward jumps is
\[\int_0^{\infty}\int_0^1v_dw^d(1-e^{-\zeta v_dw^dU_u(t)})\nu^1(w
,d\zeta )\nu^2(dw),\]
which is finite by (\ref{flreq22}).  (Recall that $v_d$ is the 
volume of the unit ball.)  The cumulative effect of the 
upward jumps on $\log U_u$ is bounded by 
\[-\int_{G_t}{\bf 1}_{D_{y,w}}(x_u)\log(1-\zeta )\xi (ds,d\beta ),\]
which has expectation
\[-\int_0^t\int_0^{\infty}\int_0^1v_dw^d\log(1-\zeta )\nu^1(w,d\zeta 
)\nu^2(dw),\]
which is again finite by (\ref{flreq22}), that is, assuming 
(\ref{flreq22}), $\tau_u$ defined in (\ref{hitinf}) is infinite.

We are going 
to prove existence by a tightness and weak convergence 
argument, so we need to view $\xi$ as a random variable 
in an appropriate metric space.  Let $\varphi\in C_b([0,\infty )\times 
{\Bbb R}^d)$ be 
strictly positive and satisfy $\int_{{\Bbb R}^d}\int_0^{\infty}\varphi 
(s,y)dsdy<\infty$. Let 
${\Bbb M}$ be the space of measures on 
${\Bbb S}=[0,\infty )\times {\cal N}\times {\Bbb R}^d\times [0,1]
\times [0,\infty )$ and define convergence in 
${\Bbb M}$ by the requirement that $\mu_n\rightarrow\mu$ if and only if
\[\int_{{\Bbb S}}\varphi (s,y)f(s,\gamma ,y,\zeta ,w)\mu_n(ds,d\gamma 
,dy,d\zeta ,dw)\rightarrow\int_{{\Bbb S}}\varphi (s,y)f(s,\gamma 
,y,\zeta ,w)\mu (ds,d\gamma ,dy,d\zeta ,dw),\]
for all $f\in C_b({\Bbb S})$.  Then ${\Bbb M}$ is metrizable and complete.

\begin{theorem}
For $\lambda <\infty$, assume that with probability one, 
$\bar{\eta}^{\lambda}(0,K\times {\Bbb K})<\infty$ for every compact $
K\subset {\Bbb R}^d$, and that conditioned on 
$\bar{\eta}^{\lambda}(0)$, the levels in $\eta^{\lambda}(0)$ are independent and uniform on $
[0,\lambda ]$.  
Then existence holds for the solution of the system 
of stochastic equations (\ref{equ1}) and (\ref{equ2})
and hence for the corresponding 
martingale problem.  

For $\lambda =\infty$, assume that $\eta (0)$ is conditionally Poisson with 
Cox measure $\Xi (0)\times\ell$ on ${\Bbb R}^d$ and 
$\sup_{y\in {\Bbb R}^d}E[\Xi (0,D_{y,1}\times {\Bbb K})]<\infty$.  For  $
\lambda <\infty$, let $U^{\lambda}$ be a 
solution of the system (\ref{equ1}) and (\ref{equ2}) with 
$\eta^{\lambda}(0)$ the restriction of $\eta (0)$ to $u\in [0,\lambda 
]$.  Then $\{(U^{\lambda},\xi )\}$ is 
relatively compact in $D_{{\Bbb R}}[0,\infty )^{\infty}\times {\Bbb M}$ and any limit point 
is a solution for the system with $\lambda =\infty$. Consequently,
existence holds for the $\lambda =\infty$ system of 
stochastic equations and hence for the corresponding 
martingale problem, and along the convergent 
subsequence, the genealogies corresponding to 
$U^{\lambda}$ converge to the genealogies of the limit.

\end{theorem}

\begin{proof}
Assume $\lambda <\infty$.  There are only countably 
many particles that ever live, and the levels must satisfy the 
countable system of equations 
\begin{eqnarray*}
U_u(t)&=&u+\int_{(b_u,t]\times G}{\bf 1}_{D_{y,w}}(x_u)(\frac {{\cal J}_{
y,w}^{\lambda}(x_u,U_u(s-),\eta (s-),v(\beta ))}{1-\zeta}-U_u(s-)
){\bf 1}_{\{v(\beta )<\lambda \}}\xi (ds,d\beta ),\end{eqnarray*}
including the initial particles with $b_u=0$.

Let $U_u^{\varepsilon}$  satisfy 
\begin{eqnarray*}
U_u^{\varepsilon}(t)&=&u+\int_{(b_u,t]\times G}{\bf 1}_{D_{y,w}}(
x_u)(\frac {{\cal J}_{y,w}^{\lambda}(x_u,U^{\varepsilon}_u(s-),\eta^{
\varepsilon}([s/\varepsilon ]\varepsilon ),v(\beta ))}{1-\zeta}-U_
u^{\varepsilon}(s-)){\bf 1}_{\{v(\beta )<\lambda \}}\xi (ds,d\beta 
).\end{eqnarray*}
With probability one, no jump in $\xi$ occurs at times of 
the form $[s/\varepsilon ]\varepsilon$, and it follows that $U^{\varepsilon}$ is uniquely 
determined.  On any bounded time interval, each particle 
is involved in only finitely many events, that is, $U^{\varepsilon}_
u$ 
jumps only finitely often, and the jumps are bounded.   
Consequently, $\{(U^{\varepsilon},\xi )\}$ is relatively compact in 
$D_{{\Bbb R}}[0,\infty )^{\infty}\times {\Bbb M}$ in the sense 
of convergence in distribution.  Selecting a convergent 
subsequence with limit $(U,\xi )$,  the only issue is the 
continuity of ${\cal J}_{y,w}^{\lambda}$.  Suppose $(\beta ,t)\in
\xi$.  Then since ${\cal J}^{\lambda}_{y,w}$ 
only depends on finitely many of the $U_u$, and, with 
probability one, no particle locations are on the boundary 
of $D_{y,w}$, the necessary continuity will be satisfied if 
$U_{u_1}(t-)\neq U_{u_2}(t-)$ for all $u_1$ and $u_2$ with $x_{u_
1},x_{u_2}\in D_{y,w}$ 
and there are no ties in the determination of $u^{*}(\beta ,\eta 
)$.  
But the first requirement holds since $U_{u_1}(t-)$ and $U_{u_2}(
t-)$ 
will be independent and uniform and the second holds 
since $v(\beta )$ will be independent of $U(t-)$.  

Essentially the same argument works 
for the relative compactness of $\{(U^{\lambda},\xi )\}$ and taking a
convergent subsequence, we obtain existence for 
$\lambda =\infty$ and convergence of the genealogies.
\end{proof}

\begin{remark}\label{unqlack}
At this point, we do not have a uniqueness result for 
the martingale problem or the stochastic equations.  
This question will be pursued elsewhere.
\end{remark}

\subsection{Spatial $\Lambda$-Fleming-Viot 
process with general offspring distribution}\label{sectslfv2}
In the discrete birth/independent thinning model 
described in the previous section, the offspring 
distribution was Poisson and the model was constructed 
so that for $\lambda =\infty$, the locations and levels of the particles 
form a spatial Poisson process that is stationary in 
time.  We now drop the Poisson assumption and allow an 
offspring distribution restricted only by the requirement 
that the expected number of offspring for an event 
$z=(y,\zeta ,w)$ in the
ball $D_{y,w}$ with thinning probability $\zeta$ is 
\[\sum_{k=0}^{\infty}kp(k,z)=\lambda\zeta |D_{y,w}|.\]

To avoid the uniqueness problem mentioned in Remark 
\ref{unqlack}, we replace ${\Bbb R}^d$ by a torus ${\Bbb T}$.  Taking 
${\Bbb U}={\Bbb T}\times [0,1]\times [0,\infty )$ and setting 
$\mu (dy,d\zeta ,dw)=dy\nu^1(w,d\zeta )\nu^2(dw)$, we assume $\mu 
({\Bbb U})<\infty$ and
 define
\[A^{\lambda}f(\eta )=\int_{{\Bbb U}}\sum_{k=1}^{\infty}p(k,z)(H_{
k,z}^{\lambda}(g,\eta )-f(\eta ))dy\nu^1(w,d\zeta )\nu^2(dw),\]
where as before, if $\eta (D_{y,w}\times {\Bbb K})=0$, $H_{k,z}^{
\lambda}(g,\eta )=f(\eta )$, and if 
$\eta (D_{y,w}\times {\Bbb K})\neq 0$, 
\begin{eqnarray*}
H_{k,z}^{\lambda}(g,\eta )&=&\prod_{(x,\kappa ,u)\in\eta ,x\notin 
D_{y,w}}g(x,\kappa ,u)\\
&&\qquad\times\int_0^{\lambda}\Big[\frac k{\lambda}(1-\frac {v^{*}}{
\lambda})^{k-1}\hat {g}_{y,w}(\kappa^{*},v^{*})\left(\frac 1{\lambda 
-v^{*}}\int_{v^{*}}^{\lambda}\hat g_{y,w}(\kappa^{*},v)dv\right)^{
k-1}\\
&&\qquad\qquad\times\prod_{(x,\kappa ,u)\in\eta ,x\in D_{y,w},u\neq 
u^{*}}g(x,\kappa ,\frac 1{1-\zeta}{\cal J}_{y,w}^{\lambda}(x,u,\eta 
,v^{*})\Big]dv^{*},\end{eqnarray*}
where as before 
$\hat {g}_{y,w}(\kappa ,u)\equiv\int g(x',\kappa ,u)\upsilon_{y,w}
(dx')$.

Again, $(x^{*},\kappa^{*},u^{*})$ is the point in $\eta$ satisfying 
$x^{*}\in D_{y,w}$ and
\[u^{*}=\mbox{\rm argmax}\{\frac {\lambda -u}{\lambda -v^{*}}:(x,
\kappa ,u)\in\eta ,x\in D_{y,w},u\geq v^{*}\}\cup \{\frac u{v^{*}}
:(x,\kappa ,u)\in\eta ,x\in D_{y,w},u\geq v^{*}\},\]
and ${\cal J}_{y,w}^{\lambda}(x,u,\eta ,v^{*})$ is obtained as in (\ref{levtran2}).

Recalling that $\bar {g}_{y,w}(\kappa )=\lambda^{-1}\int_0^{\lambda}
\int g(x,\kappa ,u)\upsilon_{y,w}(dx)du$ and 
averaging, we define 
\begin{eqnarray*}
{\cal H}^{\lambda}_{k,z}(\bar {g},\bar{\eta })&=&\frac 1{|\bar{\eta}_{
|D_{y,w}}|}\sum_{(x^{*},\kappa^{*})\in\bar{\eta}_{|D_{y,w}}}\bar {
g}_{y,w}(\kappa^{*})^k\frac 1{\bar {g}(x^{*},\kappa^{*})}\\
&&\qquad\qquad\times\prod_{(x,\kappa )\in\bar{\eta}_{|D_{y,w}},(x
,\kappa )\neq (x^{*},\kappa^{*})}((1-\zeta )+\zeta\frac 1{\bar {g}
(x,\kappa )})\end{eqnarray*}
and obtain
\begin{eqnarray*}
\alpha A^{\lambda}f(\bar{\eta })=\alpha f(\bar{\eta })\int_{{\Bbb U}}
\sum_{k=1}^{\infty}p(k,z)({\cal H}_{k,z}^{\lambda}(\bar {g},\bar{
\eta })-1)dy\nu^1(w,d\zeta )\nu^2(dw).\end{eqnarray*}

To obtain a limit as $\lambda\rightarrow\infty$, for each $z$, let $
\mu (dq,z)$ be a 
probability distribution on $[0,\infty )$ satisfying
\[\int_0^{\infty}q\mu (dq,z)=\alpha_z\equiv\zeta |D_{y,w}|,\]
and assume that as
$\lambda\rightarrow\infty$, for each $\varphi\in C_b({\Bbb R})$,
\[\sum_k\varphi (\frac k{\lambda})p^{\lambda}(k,z)\rightarrow\int_
0^{\infty}\varphi (q)\mu (dq,z).\]

These conditions imply
\[\sum_kp^{\lambda}(k,z)\int_0^{\lambda}\frac k{\lambda}(1-\frac {
v^{*}}{\lambda})^{\lambda\frac {k-1}{\lambda}}f(v^{*})dv^{*}\rightarrow
\int_0^{\infty}\int_0^{\infty}qe^{-qv^{*}}f(v^{*})dv^{*}\mu (dq,z
).\]
Observing that $\frac k{\lambda}\rightarrow q$ implies 
\begin{eqnarray*}
\left(\frac 1{\lambda -v^{*}}\int_{v^{*}}^{\lambda}\hat g_{y,w}(\kappa^{
*},\frac 1{1-\zeta}v)dv\right)^{k-1}&\rightarrow&\exp\{-q\int_{v^{
*}}^{\infty}(1-\hat {g}_{y,w}(\kappa^{*},v)dv)\}\\
&=&\exp\{-q(\hat {h}_{y,w}(\kappa^{*},v^{*})\},\end{eqnarray*}
where $\hat {h}_{y,w}(\kappa ,u)=\int_u^{\infty}(1-\hat {g}_{y,w}
(\kappa ,v))dv$,
it follows that $\sum_kp^{\lambda}(k,z)H_{k,z}^{\lambda}(g,\eta )$ converges to 
\begin{eqnarray*}
H_z(g,\eta )&=&\prod_{(x,\kappa ,u)\in\eta ,x\notin D_{y,w}}g(x,\kappa 
,u)\\
&&\qquad\times\int_0^{\infty}\int_0^{\infty}\Big[qe^{-qv^{*}}\hat {
g}_{y,w}(\kappa^{*},v^{*})\exp\{-q(1-\zeta )\hat {h}_{y,w}(\kappa^{
*},v^{*})\}\\
&&\qquad\qquad\times\prod_{(x,\kappa ,u)\in\eta ,x\in D_{y,w},u>u^{
*}}g(x,\kappa ,\frac 1{1-\zeta}(u-(u^{*}-v^{*}))\\
&&\qquad\qquad\times\prod_{(x,\kappa ,u)\in\eta ,x\in D_{y,w},u<u^{
*}}g(x,\kappa ,\frac 1{1-\zeta}u)\Big]dv^{*}\mu (dq,z),\end{eqnarray*}
and 
\begin{equation}A^{\infty}f(\eta )=\int_{{\Bbb U}}(H_z(g,\eta )-f
(\eta ))dy\nu^1(w,d\zeta )\nu^2(dw).\label{dbth2}\end{equation}

As before, setting $h^{*}_{y,w}(\kappa )=\int_0^{\infty}(1-\hat {
g}_{y,w}(\kappa ,u))du$ and 
\[{\Bbb H}_3(h^{*}_{y,w},q,\Xi ,y,\zeta ,w)=\frac 1{\Xi (D_{y,w}\times 
{\Bbb K})}\int_{D_{y,w}\times {\Bbb K}}e^{-q(1-\zeta )h^{*}_{y,w}
(\kappa )}\Xi (dx\times d\kappa ),\]
for $f(\Xi )=e^{-\int h(x,\kappa )\Xi (dx,d\kappa )}$, we have
\begin{eqnarray*}
\alpha A^{\infty}f(\Xi )=e^{-\int h(x,\kappa )\Xi (dx,d\kappa )}\int_{
[0,\infty )\times {\Bbb R}^d\times [0,1]\times [0,\infty )}({\Bbb H}_
3(h^{*}_{y,w},q,\Xi ,y,\zeta ,w)e^{\zeta\int_{D_{y,w}\times {\Bbb K}}
h(x,\kappa )\Xi (dx,d\kappa )}-1)\\
\times\mu (dq,z)dy\nu^1(w,d\zeta )\nu^2(dw).\end{eqnarray*}

Since we are assuming that $\mu ({\Bbb U})<\infty$, the martingale 
problems for the $A^{\lambda}$ and $A^{\infty}$ are well posed, and we 
have the following.

\begin{theorem}If $\eta^{\lambda}$, 
$0<\lambda <\infty$, is a solution of the martingale problems for $
A^{\lambda}$, 
and $\eta^{\lambda}(0)\Rightarrow\eta (0)$, then $\eta^{\lambda}$ converges in 
distribution to 
the unique solution of the martingale problem 
for  $A^{\infty}$ with initial distribution the distribution of $
\eta (0)$.
\end{theorem}

If $\mu (dq,z)$ is degenerate for every $z$, that is, 
$\mu (dq,z)=\delta_{\alpha_z}$, then (\ref{dbth2}) is the same as 
(\ref{dbth1}).  Of course, if $\{p^{\lambda}(k,z)\}$ is the Poisson 
distribution with mean $\lambda\zeta |D_{y,w}|$, then degeneracy holds.  
However, we can also construct non-degenerate examples, 
for example, by choosing a geometric offspring 
distribution, in which case $\mu (dq,z)$ is exponential. 

For $\lambda <\infty$, let ${\cal N}^{\lambda}$ be the collection of counting measures 
on ${\Bbb T}\times [0,\infty )$ and let
$\xi^{\lambda}$ be a Poisson random measure on 
$[0,\infty )\times {\cal N}^{\lambda}\times {\Bbb T}\times [0,1]\times 
[0,\infty )$ with mean measure
\[ds\nu^3(\{p(k,z)\},y,w,d\gamma )dy\nu^1(w,d\zeta )\nu^2(dw),\]
where $\nu^3(\{p(k,z)\},y,w,d\gamma )$ is the probability distribution 
on ${\cal N}^{\lambda}$ of the point process
\[\sum_{i=1}^K\delta_{(X_i,V_i)},\]
where $K$ is integer-valued with distribution $\{p(k,z)\}$ and 
the $(X_i,V_i)$ are independent and uniformly distributed 
over $D_{y,w}\times [0,\lambda ]$.

For $\lambda =\infty$, let $\xi$ be a Poisson random measure on 
$[0,\infty )\times {\cal N}\times [0,\infty )\times {\Bbb T}\times 
[0,1]\times [0,\infty )$ with mean measure
\[ds\nu^3(q,y,w,d\gamma )\mu (dq,z)dy\nu^1(w,d\zeta )\nu^2(dw),\]
where $\nu^3(q,y,w,d\gamma )$ is the probability distribution of the 
Poisson random measure on ${\Bbb T}\times [0,\infty )$ with mean 
measure 
\[q{\bf 1}_{D_{y,w}}(x)dxdv=q|D_{y,w}|\upsilon_{y,w}(dx)dv.\]
Under our boundedness assumption, we can take $\psi\equiv 1$ in 
Theorem \ref{mf} and in Theorem \ref{sumeq}.
The form of the stochastic equation is the same as 
in the previous section.   

Set  
$G^{\lambda}={\cal N}^{\lambda}\times {\Bbb T}\times [0,1]\times 
[0,\infty )$.

\begin{lemma}\label{eqvar2}
Any solution of the martingale problem for $A^{\lambda}$ 
can be obtained as a solution of the stochastic equation
\begin{eqnarray*}
f(\eta (t))&=&f(\eta (0))\\
&&\quad +\int_{[0,t]\times G^{\lambda}}\bigg(\frac {f(\eta (s-))}{
f(\eta_{y,w}(s-))}\prod_{(x,u)\in\gamma (\beta )}g(x,\kappa^{*}(\beta 
,\eta (s-),u)\\
&&\qquad\qquad\qquad\times\prod_{(x,\kappa ,u)\in\eta (s-),u\neq 
u^{*}(\eta (s-),v^{*}(\beta ))}g(x,\kappa ,\frac {{\cal J}^{\lambda}_{
y,w}(x,u,\eta (s-),v^{*})}{1-\zeta})\\
&&\qquad\qquad\qquad\qquad\qquad\qquad\qquad\qquad\qquad\qquad -f
(\eta (s-))\bigg)\xi^{\lambda}(ds,d\beta ).\end{eqnarray*}
For $\lambda =\infty$, the equation is the same with $G^{\lambda}$ replaced by 
$G$,  $\xi^{\lambda}$ replaced by 
$\xi$ and ${\cal J}^{\lambda}_{y,w}$ replaced by ${\cal J}^{\infty}_{
y,w}$.  
\end{lemma}

As before, the level processes satisfy 
\begin{eqnarray}
U_u(t)&=&u+\int_{(b_u,t]\times G}{\bf 1}_{D_{y,w}}(x_u)(\frac {{\cal J}_{
y,w}^{\lambda}(x_u,U_u(s-),\eta (s-),v(\beta ))}{1-\zeta}-U_u(s-)
)\nonumber\\
&&\qquad\qquad\qquad\qquad\qquad\qquad\qquad\xi^{\lambda}(ds,d\beta 
),\label{equ22}\end{eqnarray}
where $b_u=0$ if $(x,\kappa ,u)\in\eta^{\lambda}(0)$, and the death time of a 
particle satisfies 
\begin{equation}d_u=\inf\{t>0:U_u(t)>\lambda\mbox{\rm \ or }U_u(t
-)=u^{*}(\beta ,\eta (t-)),(t,\beta )\in\xi \}.\label{dthu22}\end{equation}

Passing to the $\lambda =\infty$ limit, we can derive the equation 
for the population distribution.  Let 
 $\Xi (t,dx,d\kappa )du$ be the Cox measure for $\eta (t)$.
Define
\[P(t,C)=\Xi (t,C\times {\Bbb K}),\quad C\in {\cal B}({\Bbb T}).\]
If $P(0,dx)=P(0,x)dx$, that is, $P(0,\cdot )$ is absolutely 
continuous with respect to Lebesgue measure, then since 
locations of new points are uniformly distributed over 
disks, $P(t,dx)=P(t,x)dx$ for all $t\geq 0$.  Since ${\Bbb T}\times 
{\Bbb K}$ is a 
complete, separable metric space, we can write
\[\Xi (t,dx,d\kappa )=P(t,x)\Xi_x(t,d\kappa )dx,\]
where $\Xi_x(t,\cdot )\in {\cal P}({\Bbb K})$.

\begin{theorem}\  
For $\lambda =\infty$, if $\eta (0)=\sum_{(x,\kappa ,u)\in\eta (0
)}\delta_{(x,\kappa ,u)}$ is 
conditionally Poisson with Cox measure 
$\Xi (0,dx,d\kappa )du=P(0,x)\Xi_x(0,d\kappa )dxdu$, then $\eta (
t)$ is 
conditionally Poisson with Cox measure $\Xi (t,dx,d\kappa )du$,
$\Xi (t,dx,d\kappa )=P(t,x)\Xi_x(t,d\kappa )dx$, where $\Xi_x(t,{\Bbb K}
)\equiv 1$. Then 
(dropping the $\gamma$ coordinate from $\xi$),
\[P(t,x)=P(0,x)+\int_{[0,t]\times [0,\infty )\times {\Bbb T}\times 
[0,1]\times [0,\infty )}\left(\frac q{|D_{y,w}|}-\zeta P(s-,x)\right
){\bf 1}_{D_{y,w}}(x)\xi (ds,dq,dy,d\zeta ,dw).\]
\end{theorem}

\begin{remark}

Note that in the degenerate case, $q\equiv\alpha_z=\zeta |D_{y,w}
|$, and 
$P(t,x)\equiv 1$ is a solution of this equation. 
\end{remark}\ 

To write an equation including $\Xi_x$, we need to enrich $\xi$ 
so that each point includes a coordinate that is 
independent and uniformly distributed over $[0,1]$, that is, 
for $\hat {G}=[0,1]\times [0,\infty )\times {\Bbb T}\times [0,1]\times 
[0,\infty )$,
we let $\xi$ be the Poisson random measure on 
$[0,\infty )\times\hat {G}$  with mean measure 
$dsdr\mu (dq,z)dy\nu^1(w,d\zeta )\nu^2(dw)$.  Let $K:[0,1]\times 
{\cal P}({\Bbb K})\rightarrow {\Bbb K}$ be a 
measurable function such that if $R$ is uniformly 
distributed over $[0,1]$  and $\rho\in {\cal P}({\Bbb K})$, then $
K(R,\rho )$ has 
distribution $\rho$.  Note that if an event $z=(y,\zeta ,w)$ occurs 
at time $t$, 
then the distribution of the type of the parent will be
\[\int_{D_{y,w}}\Xi_{x'}(t-,\cdot )\upsilon_{y,w}(dx').\]

\begin{theorem}
For $\varphi\in C_c({\Bbb T}\times {\Bbb K})$,
\[\begin{array}{rcl}
\langle\Xi (t),\varphi\rangle =\langle\Xi (0),\varphi\rangle +\int_{
[0,t]\times\hat {G}}\bigg[q\int_{D_{y,w}}\varphi (x,K\left(r,\int_{
D_{y,w}}\Xi_{x'}(s-,\cdot )\upsilon_{y,w}(dx')\right))\upsilon_{y
,w}(dx)\\
-\zeta\langle\Xi (s-),{\bf 1}_{D_{y,w}}\varphi\rangle\bigg]\xi (d
s,dr,dq,dy,d\zeta ,dw).\end{array}
\]
\end{theorem}

\begin{remark}
The above construction is more than complicated enough 
at least for a first reading, but still keep in mind that 
the parameters of the this model, as well as other kinds 
of population models, 
could be taken to be 
functions of $\bar{\eta}$ for $\lambda <\infty$ or $\Xi$ for $\lambda 
=\infty$.  
For example, $\mu (dq,z)$ could be replaced 
by $\mu (dq,z,\Xi (t))$, or in a genealogical construction of the 
Bolker-Pacala model \cite{BP99}, the death rate would be 
$d_0(x,\bar{\eta })=\int d(x-y)\bar{\eta }(dy)$.  
Equally, we could consider frequency dependent selection, in which
the strength of selection in favour of a particular genetic type at a
specific location depends on the current freqency of types there. 
For example \cite{FP17} consider the spatial $\Lambda$-Fleming-Viot model for a 
haploid population with
general frequency dependent selection.
Variations like this lead to a rich 
class of models in which we can combine the forces of 
ecology and genetics.  

\end{remark}

%% file: LDEXAMP2.tex


\subsection{Branching processes}
\label{sectbranching}

Next, we recover a lookdown construction for the 
Dawson-Watanabe superprocess.
Let $A_{cb,k}$ be given by (\ref{cbdef}), and let $A_{pd,k}$ be the 
pure death generator with
$d_0(x)=r(x)k$.  Let ${\cal D}_{\lambda}$ be defined as in (\ref{scrddef}) with 
${\Bbb R}^d$ replaced by $E$, and let ${\cal D}(A^{\lambda})=\{f\in 
{\cal D}_{\lambda}:\partial_ug\mbox{\rm \ is continuous}\}$.
Then, recalling the definition of $G_k^{\lambda}(u)$
from~(\ref{gklu}),
\begin{eqnarray*}
A^{\lambda}f(\eta )&=&\lambda (A_{cb,k}f(\eta )+A_{pd,k}f(\eta ))\\
&=&f(\eta )\sum_{(x,u)\in\eta}\lambda r(x)\bigg[\frac {(k+1)}{\lambda^
k}\int_u^{\lambda}\cdots\int_u^{\lambda}\left(\prod_{i=1}^kg(x,v_i)-1\right)
dv_1\cdots dv_k\\
&&\qquad\qquad\qquad\qquad\qquad\qquad\qquad\qquad +(G_k^{\lambda}
(u)+ku)\frac {\partial_ug(x,u)}{g(x,u)}\bigg]\\
&\rightarrow&f(\eta )\sum_{(x,u)\in\eta}r(x)(k+1)k\left(\int_u^{\infty}
\Big(g(x,v)-1\Big)dv+\frac 12u^2\frac {\partial_ug(x,u)}{g(x,u)}\right)\\
&=&A^{\infty}f(\eta ),\end{eqnarray*}
and
\[\alpha A^{\infty}f(\Xi )=e^{-\int_Eh(x)\Xi (dx)}\int_Er(x)\frac {
k(k+1)}2h^2(x)\Xi (dx),\]
which is the generator of a Dawson-Watanabe process 
without any spatial motion.  
(See Section 1.5 of \cite{Eth00} or Section 3.4 of 
\cite{KR11}.)  Note that for finite $\lambda$, each birth event 
produces $k$ offspring.

For more general offspring distribution, one can take 
\begin{eqnarray*}
A^{\lambda}f(\eta )&=&\lambda\int_{{\Bbb U}}(A_{cb,k(z)}f(\eta )+
A_{pd,k(z)}f(\eta ))\mu (dz)\\
&=&f(\eta )\int_{{\Bbb U}}\sum_{(x,u)\in\eta}\lambda r(x,z)\bigg[\frac {
(k(z)+1)}{\lambda^{k(z)}}\int_u^{\lambda}\cdots\int_u^{\lambda}(\prod_{
i=1}^{k(z)}g(x,v_i)-1)dv_1\cdots dv_{k(z)}\\
&&\qquad\qquad\qquad\qquad\qquad\qquad\qquad +(G_{k(z)}^{\lambda}
(u)+k(z)u)\frac {\partial_ug(x,u)}{g(x,u)}\bigg]\mu (dz)\\
&\rightarrow&f(\eta )\int_{{\Bbb U}}\sum_{(x,u)\in\eta}r(x,z)(k(z
)+1)k(z)\left(\int_u^{\infty}(g(x,v)-1)dv+\frac 12u^2\frac {\partial_
ug(x,u)}{g(x,u)}\right)\mu (dz)\\
&=&A^{\infty}f(\eta ),\end{eqnarray*}
assuming $\sup_{x\in E}\int_{{\Bbb U}}r(x,z)(k(z)+1)k(z)\mu (dz)<
\infty$.  We can 
take $\psi$ in Theorem \ref{mf}\ to be of the form 
$\sum_l\delta_l\eta (K_l\times [0,l])$ for appropriately selected $
\delta_l$.

This construction is a special case of the results in 
\cite{KR11} which considers more general offspring 
distributions (for example, offspring distributions 
without second moments), and other variants of 
branching processes including random environments and 
processes conditioned on extinction and nonextinction.

\subsection{Spatially interacting Moran model}\label{sectMoran}

Consider $A_{dr,3}$, as defined in \S\ref{sectoneforone}, 
in the special case in which the sum is over all subsets 
with $|S|=2$.  In other words, each replacement event 
involves just two individuals.  Specifically, we take 
$r(S,z)=r(x,x')$ for $S=\{x,x'\}$.  We include independent 
motion with generator $B\subset C_b(E)\times C_b(E)$, 
set $q(x,z,dy)=\delta_x(dy)$, and 
assume $r(x,x')=r(x',x)$.  (Note that this symmetry is 
needed for $\alpha Af$ to be a generator applied to $\alpha f$.)  The 
generator becomes 
\begin{equation}Af(\eta )=f(\eta )\sum_{(x,u)\in\eta}\frac {Bg(x,
u)}{g(x,u)}+f(\eta )\sum_{(x,u)\neq (x',u')\in\eta}r(x,x'){\bf 1}_{
\{u'<u\}}(\frac {g(x',u)}{g(x,u)}-1)\label{mmgen}\end{equation}
for 
\[f\in {\cal D}(A)=\{f\in {\cal D}_{\lambda}:g\in {\cal D}(B)\}\]
and 
\begin{equation}\alpha Af(\bar{\eta })=\alpha f(\bar{\eta })\sum_{
x\in\bar{\eta}}\frac {B\bar {g}(x)}{\bar {g}(x)}+\alpha f(\bar{\eta }
)\sum_{\{x,x'\}\subset\bar{\eta}}r(x,x')(\frac 12\frac {\bar {g}(
x')}{\bar {g}(x)}+\frac 12\frac {\bar {g}(x)}{\bar {g}(x')}-1),\label{avmmgen}\end{equation}
that is, at rate $r(x,x')$ one of the pair is killed and 
replaced by a copy of the other.  

Since either particles move or a particle of one 
type is replaced by a particle of another type, if the 
initial number of particles is finite, then, as in the 
classical Moran model, the total 
number of  particles is preserved. Consequently, 
if $r(x,x')$ is bounded, we can apply Theorem 
\ref{mf} with $\psi (\eta )=1+|\eta |^2$.  If the number of particles 
is infinite, the following condition is useful.

\begin{condition}\label{domstr}
Let ${\cal K}=\{K_1,K_2,\ldots \}$, $K_k\subset E$.
For each $f\in {\cal D}(A)$, $f(\eta )=\prod_{(x,u)\in\eta}g(x,u)$,
there exists $K_g\in {\cal K}$ such that $g(x,u)=1$ 
and $Bg(x,u)=0$ for all $x\notin K_g$ and $r(x,x')=0$ for $x\notin 
K_g$ 
and $x'$ in the support of $1-g$.
\end{condition}

\begin{lemma}
Assume Condition \ref{domstr}. Then for each $f\in {\cal D}(A)$, 
there exists $c_f$ such that
\[|Af(\eta )|\leq c_f(\bar{\eta }(K_g)+\int_{K_g\times K_g}r(x,x'
)\bar{\eta }(dx)\bar{\eta }(dx')).\]
Then for  $\delta_k>0$, $k=1,2,\ldots$, $\psi$ of the form 
\[\psi (\eta )=\sum_k\delta_k(\bar{\eta }(K_k)+\int_{K_k\times K_
k}r(x,x')\bar{\eta }(dx)\bar{\eta }(dx'))\]
satisfies (\ref{opest}).
\end{lemma}

\begin{remark}
Of course, to 
apply Theorem \ref{mf} one must verify that
\begin{equation}\int_0^tE[\tilde{\psi }(\bar{\eta }(s))]ds<\infty 
,\quad t\geq 0\label{mfcnd2}\end{equation}
for the solution of interest.
For example, in the spatially interacting Moran model in 
\cite{GLW05}, particles have a location and type 
($(x,\kappa )\in E=G\times {\Bbb K}$ rather than $x$) for a countable set $
G$,
 $r((x,\kappa ),(x',\kappa'))=\gamma {\bf 1}_{\{x=x'\}}$,
the locations evolve independently according to a Markov 
chain with transition intensities $q(x,y)$, that is, 
\[Bg(x,\kappa )=\sum_{y\in G}q(x,y)(g(y,\kappa )-g(x,\kappa ))+Cg
(x,\kappa ),\]
where $C$ is a mutation operator that acts only on the 
type.  The location Markov chain is assumed to satisfy 
estimates that imply $E[\eta (t,\{x\}\times {\Bbb K})^2]<\infty$ provided $
\eta (0)$ 
satisfies specified conditions. 
Consequently, if we take 
$K_k=G_k\times {\Bbb K}$ for finite subsets $G_k$, we can select $
\delta_k$ so 
that (\ref{mfcnd2}) is satisfied.  
\end{remark}

Note that $\lambda$  does not appear in the formula for the 
generator (\ref{mmgen}).  Consequently, the same formula 
gives the limiting generator as $\lambda\rightarrow\infty$, and with reference 
to (\ref{prws}),
\begin{eqnarray*}
\alpha A^{\infty}f(\Xi )&=&e^{-\int_Eh(x)\Xi (dx)}\Big[-\int_EBh(
x)\Xi (dx)\\
&&\qquad\qquad\qquad +\int_{E\times E}r(x,x')\left(\frac 12h^2(x)
+\frac 12h^2(x')-h(x')h(x)\right)\Xi (dx)\Xi (dx')\Big].\end{eqnarray*}

For $\lambda =\infty$, if the number of particles below any level is 
finite, we can take $\psi (\eta )=\sum_{l=1}^{\infty}\delta_l(1+\eta 
(E\times [0,l])^2)$.  If the 
number of particles below a level is infinite, then $\psi$ of 
the form 
\[\psi (\eta )=\sum_{k,l}\delta_{k,l}(\eta (K_k\times [0,l])+\int_{
K_k\times [0,l]\times K_k\times [0,l]}r(x,x')\eta (dx,du)\eta (dx'
,du'))\]
meets the requirements of Theorem \ref{mf}.
  
  For the limiting process, one 
can also see that mass is preserved directly from the 
limiting generator.  Suppose $\Xi$ is a solution of the 
martingale problem with $\Xi (0,E)<\infty$.
 Take $h(x)\equiv c>0$, and observe that $e^{-c\Xi (t,E)}$ 
is a martingale. But, in general, if $M$ and $M^2$ are both martingales, 
then $M$ must be constant, so consider $e^{-c\Xi (t,E)}$ and 
$e^{-2c\Xi (t,E)}$.

If $r(x,x')\equiv\gamma$ and $\Xi (0,E)=1$, then $\Xi$ is a neutral 
Fleming-Viot process.  Since the set of levels is fixed, 
in this case, the lookdown construction is equivalent to 
the construction given in \cite{DK96}.  If as above, 
$r((x,\kappa ),(x',\kappa'))=\gamma {\bf 1}_{\{x=x'\}}$, then the lookdown construction 
for $\lambda =\infty$ is just the lookdown construction for the 
interacting Fisher-Wright diffusions discussed in 
\cite{GLW05}.

\subsection{A stochastic partial differential 
equation}\label{sectspde}
Consider a spatially interacting Moran model with both 
location $x\in\lambda^{-1}{\Bbb Z}$ and type $\kappa\in {\Bbb K}$.  Assume that the particle 
locations follow a simple symmetric random walk, and 
for simplicity, assume that the types of the particles do 
not change.  Killing and replacement of the previous section
now takes place locally 
at each site.  The generator then becomes 
\begin{eqnarray*}
Af(\eta )=f(\eta )\sum_{(x,\kappa ,u)\in\eta}\lambda^2\frac {g(x+
\lambda^{-1},\kappa ,u)+g(x-\lambda^{-1},\kappa ,u)-2g(x,\kappa ,
u)}{2g(x,\kappa ,u)}\\
+f(\eta )\sum_{(x,\kappa ,u)\neq (x',\kappa',u')\in\eta}\lambda {\bf 1}_{
\{x=x'\}}{\bf 1}_{\{u'<u\}}(\frac {g(x,\kappa',u)}{g(x,\kappa ,u)}
-1).\end{eqnarray*}
Note that particles move independently, so that the 
number of particles at a site will fluctuate; however, if 
the initial site occupancies are i.i.d.~Poisson, then they 
will remain i.i.d.~Poisson. The averaged generator becomes 
\begin{eqnarray*}
\alpha Af(\eta )=\alpha f(\bar{\eta })\sum_{(x,\kappa )\in\bar{\eta}}
\lambda^2\frac {\bar {g}(x+\lambda^{-1},\kappa )+\bar {g}(x-\lambda^{
-1},\kappa )-2\bar {g}(x,\kappa )}{2\bar {g}(x,\kappa )}\\
+\alpha f(\bar{\eta })\sum_{(x,\kappa )\neq (x',\kappa')\in\bar{\eta}}\frac {
\lambda}2{\bf 1}_{\{x=x'\}}(\frac {\bar {g}(x,\kappa')}{\bar {g}(
x,\kappa )}-1)\end{eqnarray*}
(c.f.~(\ref{mmgen}) and (\ref{avmmgen})).

Let $(X_u^{\lambda}(t),\kappa_u(t))$ denote the position and type of a 
particle at level $u$, Assume that $\{(X_u^{\lambda}(0),\kappa_u(
0),u)\}$ 
determines a conditionally Poisson random measure with 
Cox measure $\lambda^{-1}\times\ell^{\lambda}(dx)\times\nu_0(x,d\kappa 
)\times du$ on $(\lambda^{-1}{\Bbb Z}\times {\Bbb K}\times [0,\lambda 
])$, 
where $\ell^{\lambda}$ is counting measure on $\lambda^{-1}{\Bbb Z}$ and $
\nu_0$ is a 
random mapping $\nu_0:x\in {\Bbb R}\rightarrow\nu_0(x,\cdot )\in 
{\cal P}({\Bbb K})$.  Note that as 
$\lambda\rightarrow\infty$, the $\{X_u^{\lambda}-X_u^{\lambda}(0)
\}$ converge to independent standard 
Brownian motions $\{W_u\}$.  

For $u'<u$, let $L_{u'u}^{\lambda}(t)$ be the number of times by time $
t$ 
that there has been a `lookdown' from $u$ to $u'$.  Then 
$L_{u'u}^{\lambda}$ is a counting process with integrated intensity 
\[\Lambda_{u'u}^{\lambda}(t)=\lambda\int_0^t{\bf 1}_{\{X_u^{\lambda}
(s)=X_{u'}^{\lambda}(s)\}}ds,\]
and we can write
\[L_{u'u}(t)=Y_{u'u}(\Lambda_{u'u}^{\lambda}(t)),\]
where the $Y_{u'u}$ are independent unit Poisson processes 
and are independent of $X_{u'u}^{\lambda}(t)\equiv X_{u'}^{\lambda}
(t)-X_u^{\lambda}(t)$.
To identify the limit of $\Lambda^{\lambda}_{u'u}$ as $\lambda\rightarrow
\infty$, define  
\[N_{u'u}^{\lambda}(t)=\#\{s\leq t:X_{u'u}^{\lambda}(s-)=0,X_{u'u}^{
\lambda}(s)\neq 0\}.\]
Then $N_{u'u}^{\lambda}$ is a counting process with intensity 
$\lambda^2{\bf 1}_{\{X_{u'u}^{\lambda}(t)=0\}}$.  Define 
\[\tilde {N}_{u'u}^{\lambda}(t)=N_{u'u}^{\lambda}(t)-\int_0^t\lambda^
2{\bf 1}_{\{X_{u'u}^{\lambda}(s)=0\}}ds.\]
Then 
\begin{eqnarray*}
|X^{\lambda}_{u'u}(t)|=|X_{u'u}^{\lambda}(0)|+\int_0^t\mbox{\rm sign}
(X_{u'u}^{\lambda}(s-))dX_{u'u}^{\lambda}(s)+\frac 1{\lambda}\tilde {
N}_{u'u}^{\lambda}(t)\\
+\lambda\int_0^t{\bf 1}_{\{X_{u'u}^{\lambda}(s)=0\}}ds.\end{eqnarray*}
Since $X_{u'u}^{\lambda}\Rightarrow X_{u'u}=X_{u'}-X_u$ and $\lambda^{
-1}\tilde {N}_{u'u}^{\lambda}\Rightarrow 0$, it follows 
that $X_{u'u}^{\lambda}$ and $\Lambda_{u'u}^{\lambda}=\lambda\int_
0^t{\bf 1}_{\{X_{u'u}^{\lambda}(s)=0\}}ds$ converge to $X_{u'u}$ and 
$\Lambda_{u'u}$ respectively satisfying Tanaka's formula
\begin{equation}|X_{u'u}(t)|=|X_{u'u}(0)|+\int_0^t\mbox{\rm sign}
(X_{u'u}(s-))dX_{u'u}(s)+\Lambda_{u'u}(t).\label{tan}\end{equation}
An application of It\^o's formula gives
\begin{equation}\Lambda_{u'u}(t)=\lim_{\varepsilon\rightarrow 0}\frac 
1{\varepsilon}\int_0^t{\bf 1}_{(-\varepsilon ,\varepsilon )}(X_{u'}
(s)-X_u(s))ds.\label{loctim}\end{equation}

To summarize, $\{(X_u(0),\kappa_u(0),u)\}$ determines a conditionally 
Poisson random measure with Cox measure 
$dx\times\nu_0(x,d\kappa )\times du$ and $X_u(t)=X_u(0)+W_u(t)$, where the $
W_u$ 
are independent, standard Brownian motions.  $L_{u'u}$ is 
determined by (\ref{tan}) and 
\[L_{u'u}(t)=Y_{u'u}(\Lambda_{u'u}(t)),\]
where the $Y_{u'u}$ are independent unit Poisson processes 
that are independent of $\{(X_u(0),\kappa_u,u)\}$ and $\{W_u\}$. The 
particle types satisfy
\[\kappa_u(t)=\kappa_u(0)+\sum_{u'<u}\int_0^t(\kappa_{u'}(s-)-\kappa_
u(s-))dL_{u'u}(s).\]

Then $\{(X_u(t),\kappa_u(t),u)\}$ determines a conditionally Poisson 
random measure with Cox measure 
\[\Xi_t(dx,d\kappa )\times du=dx\times\nu_t(x,d\kappa )\times du.\]
For details and related results 
see Buhr \cite{Buh02}.  In particular, for
$\varphi (x,\kappa )$ bounded, $C^2$ in $x$, and having compact support in 
$x$,
\[M_{\varphi}(t)=\langle\Xi_t,\varphi\rangle -\int_0^t\langle\Xi_
s,\frac 12\partial_x^2\varphi\rangle ds\]
is a $\{{\cal F}_t^{\Xi}\}$-martingale with quadratic variation
\[[M_{\varphi}]_t=\int_0^t\int_{{\Bbb R}}\int_{{\Bbb K}\times {\Bbb K}}
(\varphi (x,\kappa')-\varphi (x,\kappa ))^2\nu_s(x,d\kappa')\nu_s
(x,d\kappa )dxds,\]
identifying $\Xi$ as a solution of a martingale problem.

Suppose ${\Bbb K}=\{0,1\}$ and $\nu_s(x)\equiv\nu_s(x,\{1\})$.  Then taking 
$\varphi (x,\kappa )=\kappa\psi (x)$,
\[M_{\psi}(t)=\int_{{\Bbb R}}\psi (x)\nu_t(x)dx-\int_0^t\int_{{\Bbb R}}\frac 
12\psi^{\prime\prime}(x)\nu_s(x)dxds\]
is a martingale with quadratic variation
\[[M_{\psi}]_t=\int_0^t\int_{{\Bbb R}}2\psi^2(x)\nu_s(x)(1-\nu_s(
x))dxds,\]
which implies $\nu_t$ is a weak solution of the stochastic 
partial differential equation
\begin{eqnarray}
\int_{{\Bbb R}}\psi (x)\nu_t(x)dx=\int_{{\Bbb R}}\psi (x)\nu_0(x)
dx+\int_0^t\int_{{\Bbb R}}\frac 12\psi^{\prime\prime}(x)\nu_s(x)d
xds\label{mtspde}\\
+\int_{[0,t]\times {\Bbb R}}\psi (x)\sqrt {2\nu_s(x)(1-\nu_s(x)}W
(ds,dx),\nonumber\end{eqnarray}
where $W$ is Gaussian white noise on $[0,\infty )\times {\Bbb R}$ with 
$E[W(A)W(B)]=\ell (A\cap B)$ for Lebesgue measure $\ell$ on 
$[0,\infty )\times {\Bbb R}$.

\subsection{Voter model}\label{sectvote}\ The stochastic partial 
differential equation (\ref{mtspde}) is a special case of 
the equation that arises as the limit of rescaled voter 
models in the work of Mueller and Tribe  \cite{MT95}. 
To see the relationship of their work to our current 
approach, we give a construction of a class of voter 
models.

Let $E={\Bbb Z}\times {\Bbb K}$, where ${\Bbb Z}$ is the space of locations and $
{\Bbb K}$ 
the space of types.  We assume that there is one 
particle at each location, and consider
\[\begin{array}{rcl}
A_{dr,3}f(\eta )=f(\eta )\sum_{i\neq j}r(|x_i-x_j|){\bf 1}_{\{u_i
<u_j\}}\left(\frac {g(x_i,\kappa_i,u_i)g(x_j,\kappa_i,u_j)+g(x_i,
\kappa_i,u_j)g(x_j,\kappa_i,u_i)}{2g(x_i,\kappa_i,u_i)g(x_j,\kappa_
j,u_j)}-1\right)\end{array}
\]
where
$$\sigma^2\equiv\frac{1}{2}\sum_ll^2r(l)<\infty.$$
Then
\[\alpha A_{dr,3}f(\bar{\eta })=\alpha f(\bar{\eta })\sum_{i<j}r(
|x_i-x_j|)\left(\frac 12\frac {\bar {g}(x_j,\kappa_i)}{\bar {g}(x_
j,\kappa_j)}+\frac 12\frac {\bar {g}(x_i,\kappa_j)}{\bar {g}(x_i,
\kappa_i)}-1\right)\]
which is the generator for a voter model.  Particle 
motion involves two particles exchanging places, so in 
this model, the occupancy at each site is preserved.

Note that the 
collection of levels does not change, and the location of the 
particle associated with level $u$ will satisfy a 
stochastic equation of the form
\[X_u(t)=X_u(0)+\sum_{k<l}\int_{[0,t]\times \{0,1\}}\theta ({\bf 1}_{
\{X_u(s-)=l\}}(k-l)+{\bf 1}_{\{X_u(s-)=k\}}(l-k))\xi_{kl}(ds,d\theta 
),\]
where the $\xi_{kl}$ are independent Poisson random measures 
with mean measures 
\[r(|k-l|)(\frac 12\delta_1(d\theta )+\frac 12\delta_0(d\theta ))
ds.\]
For $k>l$, assume $\xi_{kl}\equiv\xi_{lk}$.  Let $U_l(t)$ and $\hat {
K}_l(t)$ denote 
the level and type of the particle with location $l$.  Then 
the type for the particle with level $u$ satisfies 
\[K_u(t)=K_u(0)+\sum_{l\neq k}\int_{[0,t]\times \{0,1\}}{\bf 1}_{
\{U_l(s-)<u\}}{\bf 1}_{\{X_u(s-)=k\}}(\hat {K}_l(s-)-K_u(s-))\xi_{
kl}(ds,d\theta ).\]

Now, as $\lambda\rightarrow\infty$, assume that $\{(\lambda^{-1}X_
u(0),K_u(0),u)\}$ converges 
to a conditionally Poisson point process on ${\Bbb R}\times {\Bbb K}
\times [0,\infty )$ 
with Cox measure $dx\times\nu_0(x,d\kappa )\times du$.  Set 
$X_u^{\lambda}(t)=\frac 1{\lambda}X_u(\lambda^2t)$ and $K_u^{\lambda}
(t)=K_u(\lambda^2t)$.  Then $X_u^{\lambda}$ is a 
martingale with quadratic variation 
\[[X_u^{\lambda}]_t=\sum_{k<l}\frac 1{\lambda^2}\int_{[0,\lambda^
2t]\times \{0,1\}}\theta ({\bf 1}_{\{X_u(s-)=l\}}(k-l)^2+{\bf 1}_{
\{X_u(s-)=k\}}(l-k)^2)\xi_{kl}(ds,d\theta )\]
and
\[[X_u^{\lambda}]_t\rightarrow\frac 12\sum_{k<l}(k
-l)^2r(|k-l|)t =\sigma^2t.\]
In addition, for $u\neq u'$,
\[[X_u^{\lambda},X^{\lambda}_{u'}]_t\rightarrow 0,\]
so the $X_u^{\lambda}$ converge to a collection of independent 
Brownian motions $X_u$.

For $u'<u$, let
\[N_{u',u}^{\lambda}(t)=\sum_{l\neq k}\int_{[0,\lambda^2t]\times 
\{0,1\}}{\bf 1}_{\{X_{u'}(s-)=l\}}{\bf 1}_{\{X_u(s-)=k\}}\xi_{kl}
(ds,d\theta ).\]
Then $N^{\lambda}_{u',u}$ is a counting process with integrated intensity
\[\int_0^t\lambda^2r(\lambda |X^{\lambda}_{u'}(s)-X_u^{\lambda}(s
)|)ds.\]
Under appropriate time-scaling conditions, this integral should 
converge to a constant times the intersection local time given 
in (\ref{loctim}).  Then, up to changes in parameters, the 
limit of the lookdown construction would be the same as 
in \S\ref{sectspde}.


%% file: POPAPP2.tex

\renewcommand {\theequation}{A.\arabic{equation}}
\appendix

\setcounter{equation}{0}

\section{Appendix}
\subsection{Poisson identities} \label{poisson identities}
\begin{lemma}\label{app1} 
If $\xi$ is a Poisson random measure on $S$ with $\sigma$-finite 
mean measure $\nu$ and $f\in L^1(\nu )$, then 
\begin{equation}E[e^{\int_Sf(z)\xi (dz)}]=e^{\int_S(e^f-1)d\nu},\label{lapfunc}\end{equation}
\begin{equation}E[\int_Sf(z)\xi (dz)]=\int_Sfd\nu ,\quad Var(\int_
Sf(z)\xi (dz))=\int_Sf^2d\nu ,\label{poisvar}\end{equation}
allowing $\infty =\infty$.

Letting $\xi =\sum_i\delta_{Z_i}$, for $g\geq 0$ with 
$\log g\in L^1(\nu )$,
\[E[\prod_ig(Z_i)]=e^{\int_S(g-1)d\nu}.\]
Similarly, if $hg,g-1\in L^1(\nu )$, then
\begin{equation}E[\sum_jh(Z_j)\prod_ig(Z_i)]=\int_Shgd\nu e^{\int_
S(g-1)d\nu},\label{pois1fac}\end{equation}
\begin{equation}E[\sum_{i\neq j}h(Z_i)h(Z_j)\prod_kg(Z_k)]=(\int_
Shgd\nu )^2e^{\int_S(g-1)d\nu},\label{pois2fac}\end{equation}
and more generally, if $\nu$ has no atoms and $r\in M(S\times S)$, 
$r\geq 0$,
\begin{equation}E[\sum_{i\neq j}r(Z_i,Z_j)\prod_{k\neq i,j}g(Z_k)
]=\int_{S\times S}r(x,y)\nu (dx)\nu (dy)e^{\int (g-1)d\nu},\label{prws}\end{equation}
allowing $\infty =\infty$.
\end{lemma}

\begin{proof}
The independence properties of $\xi$ imply (\ref{lapfunc}) 
and (\ref{poisvar}) for 
simple functions.  The general case follows by 
approximation.  

To prove (\ref{prws}), it is enough to consider a finite 
measure $\nu$ and bounded continuous $r$ and $g$ and extend by 
approximation.
Let $\{B_k^n\}$ be a partition of $S$ with $\mbox{\rm diam}(B_k^n
)\leq n^{-1}$, and let 
$x_k^n\in B_k^n$.  Define
\[\xi_n=\sum_k\delta_{x_k^n}{\bf 1}_{\{\xi (B_k^n)>0\}}.\]
Then $\xi_n\rightarrow\xi$ in the sense that $\int fd\xi_n\rightarrow
\int fd\xi$ for every 
bounded
continuous $f$, and 
\begin{eqnarray*}
\sum_{i\neq j}r(x_i^n,x_j^n){\bf 1}_{\{\xi (B_i^n)>0\}}{\bf 1}_{\{
\xi (B_j^n)>0\}}\prod_{k\neq i,j}(g(x_k^n){\bf 1}_{\{\xi (B_k^n)>
0\}}+{\bf 1}_{\{\xi (B_k^n)=0\}})\\
\rightarrow\sum_{i\neq j}r(Z_i,Z_j)\prod_{k\neq i,j}g(Z_k).\end{eqnarray*}
By independence, the expectation of the left side is
\begin{eqnarray*}
&&\sum_{i\neq j}r(x_i^n,x_j^n)(1-e^{-\nu (B_i^n)})(1-e^{-\nu (B_j^
n)})\prod_{k\neq i,j}(g(x_k^n)(1-e^{-\nu (B_k^n)})+e^{-\nu (B_k^n
)})\\
&&\qquad\qquad\qquad\approx\sum_{i\neq j}r(x_i^n,x_j^n)\nu (B_i^n
)\nu (B_j^n)\exp\{\sum_{k\neq i,j}(g(x_k^n)-1)\nu (B_k^n)\}\\
&&\qquad\qquad\qquad\rightarrow\int_{S\times S}r(x,y)\nu (dx)\nu 
(dy)e^{\int (g-1)d\nu},\end{eqnarray*}
where the convergence follows from the assumed 
continuity of $r$ and $g$ and the fact that
$\sum_i\nu (B_i^n)^2\rightarrow 0$.

The other identities follow in a similar 
manner.  Note that the integrability of the random 
variables in the expectations above can be verified by 
replacing $g$ by $(g\vee (-a))\wedge a{\bf 1}_A+{\bf 1}_{A^c}$ and $
h$ by $(h\vee (-a))\wedge a{\bf 1}_A$ for 
$0<a<\infty$ and $\nu (A)<\infty$ and passing to the limit as 
$a\rightarrow\infty$ and $A\nearrow E$.
\end{proof}


\subsection{Markov mapping theorem}\label{mpsect} 
The 
following theorem (extending Corollary 3.5 from 
\cite{Kur98}) plays an essential role in justifying the 
particle representations and can also be used to prove 
uniqueness for the corresponding measure-valued 
processes.   Let $(S,d)$ and $(S_0,d_0)$ be complete, separable 
metric spaces, $B(S)\subset M(S)$ be the Banach space of 
bounded measurable functions on $S$, with $\|f\|=\sup_{x\in S}|f(
x)|$, and 
$C_b(S)\subset B(S)$ be the subspace of bounded continuous 
functions.  An operator $A\subset B(S)\times B(S)$ is {\em dissipative\/} if 
$\Vert f_1-f_2-\epsilon (g_1-g_2)\Vert\geq\Vert f_1-f_2\Vert$ for all $
(f_1,g_1),(f_2,g_2)\in A$ and 
$\epsilon >0$; $A$ is a {\em pre-generator\/} if $A$ is dissipative and there 
are sequences of functions $\mu_n:S\rightarrow {\cal P}(S)$ and $
\lambda_n:S\rightarrow [0,\infty )$ 
such that for each $(f,g)\in A$ 
\begin{equation}g(x)=\lim_{n\rightarrow\infty}\lambda_n(x)\int_S(
f(y)-f(x))\mu_n(x,dy)\label{gencomp}\end{equation}
for each $x\in S$. $A$ is {\em countably determined }
if there exists a countable subset 
$\{g_k\}\subset {\cal D}(A)\cap\bar {C}(S)$ 
such that 
every solution of the martingale problem for $\{(g_k,Ag_k)\}$ 
is a solution of the martingale problem for $A$.  (For 
example, $A$ is  countably determined if it is {\em graph }
{\em separable\/} in the sense that
there exists $\{(g_k,h_k)\}\subset A\cap\bar {C}(S)\times B(S)$ 
such that $A$ is contained in the bounded pointwise 
closure of $\{(g_k,h_k)\}$.)  These  
conditions are satisfied by essentially all operators $A$ 
that might reasonably be thought to be generators of 
Markov processes.  Note that $A$ is graph separable if 
$A\subset L\times L$, where $L\subset B(S)$ is separable in the sup norm 
topology, for example, if $S$ is locally compact and $L$ is 
the space of continuous functions vanishing at infinity.

A collection of functions $D\subset\bar {C}(S)$ is {\em separating\/} if 
$\nu ,\mu\in {\cal P}(S)$ and $\int_Sfd\nu =\int_Sfd\mu$ for all $
f\in D$ imply 
$\mu =\nu$.

For an $S_0$-valued, measurable process $Y$, $\hat {{\cal F}}^Y_t$ will denote the 
completion of the $\sigma$-algebra $\sigma (Y(0),\int_0^rh(Y(s))d
s,r\leq t,h\in B(S_0))$.  
For almost every $t$, $Y(t)$ will be $\hat {{\cal F}}_t^Y$-measurable, but in 
general, $\hat {{\cal F}}^Y_t$  does not contain ${\cal F}^Y_t=\sigma 
(Y(s):s\leq t)$.  Let 
${\bf T}^Y=\{t:Y(t)\mbox{\rm \ is }\hat {{\cal F}}_t^Y\mbox{\rm \ measurable}
\}$.  If $Y$ is c\`adl\`ag and has no 
fixed points of discontinuity (that is, for every $t$, 
$Y(t)=Y(t-)$ a.s.), then ${\bf T}^Y=[0,\infty )$. Let $D_S[0,\infty 
)$ denote the 
space of c\`adl\`ag, $S$-valued functions with the Skorohod 
topology, and $M_S[0,\infty )$ denotes the space of Borel 
measurable functions, $x:[0,\infty )\rightarrow S$, topologized by 
convergence in Lebesgue measure.

\begin{theorem}\label{mf}Let 
$(S,d)$ and $(S_0,d_0)$ be complete, separable metric spaces.  
Let $A\subset\bar {C}(S)\times C(S)$ and $\psi\in C(S)$, $\psi\geq 
1$.  Suppose that for 
each $f\in {\cal D}(A)$ there exists $c_f>0$ such that 
\begin{equation}|Af(x)|\leq c_f\psi (x),\quad x\in A,\label{opest}\end{equation}
and define $A_0f(x)=Af(x)/\psi (x)$.

Suppose that $A_0$ is a countably determined pre-generator, 
and suppose that ${\cal D}(A)={\cal D}(A_0)$ is closed under multiplication and 
is separating.  Let $\gamma :S\rightarrow S_0$ be Borel measurable, and let 
$\alpha$ be a transition function from $S_0$ into $S$ 
($y\in S_0\rightarrow\alpha (y,\cdot )\in {\cal P}(S)$ is Borel measurable) satisfying 
$\int h\circ\gamma (z)\alpha (y,dz)=h(y)$, $y\in S_0$, $h\in B(S_
0)$, that is, 
$\alpha (y,\gamma^{-1}(y))=1$.  Assume that $\tilde{\psi }(y)\equiv
\int_S\psi (z)\alpha (y,dz)<\infty$ for each 
$y\in S_0$ and define 
\[C=\{(\int_Sf(z)\alpha (\cdot ,dz),\int_SAf(z)\alpha (\cdot ,dz)
):f\in {\cal D}(A)\}\;.\]
Let $\mu_0\in {\cal P}(S_0)$, and define $\nu_0=\int\alpha (y,\cdot 
)\mu_0(dy)$. 

\begin{itemize}

\item[a)]  If 
$\tilde {Y}$ satisfies $\int_0^tE[\tilde{\psi }(\tilde {Y}(s))]ds
<\infty$ for all $t\geq 0$ and $\tilde {Y}$ 
is a solution of the martingale problem for $(C,\mu_0)$, 
then there exists a solution $X$ of the martingale problem 
for $(A,\nu_0)$ such that $\tilde {Y}$ has the same distribution on 
$M_{S_0}[0,\infty )$ as $Y=\gamma\circ X$.  If $Y$ and $\tilde {Y}$ are c\`adl\`ag, then $
Y$ 
and $\tilde {Y}$ have the same distribution on $D_{S_0}[0,\infty 
)$.

\item[b)]For  $t\in {\bf T}^Y$,
\begin{equation}P\{X(t)\in\Gamma |\hat {{\cal F}}^Y_t\}=\alpha (Y
(t),\Gamma ),\quad\Gamma\in {\cal B}(S).\label{rpid}\end{equation}

\item[c)]If, in addition, uniqueness holds for the martingale 
problem for $(A,\nu_0)$, then
uniqueness holds for the $M_{S_0}[0,\infty )$-martingale 
problem for $(C,\mu_0)$.  If $\tilde {Y}$ has sample paths in $D_{
S_0}[0,\infty )$, then 
uniqueness holds for the $D_{S_0}[0,\infty )$-martingale problem for 
$(C,\mu_0)$.  

\item[d)]If uniqueness holds for the martingale problem for 
$(A,\nu_0)$,
then $Y$ restricted to ${\bf T}^Y$ is a Markov process.

\end{itemize}

\end{theorem}

\begin{remark}\label{discgen}
Theorem \ref{mf} can be extended to cover a large class 
of generators whose range contains discontinuous 
functions.  (See \cite{Kur98}, Corollary 3.5 and Theorem 
2.7.)  In particular, suppose $A_1,\ldots ,A_m$ satisfy the 
conditions of Theorem \ref{mf} for a common domain 
${\cal D}={\cal D}(A_1)=\cdots ={\cal D}(A_m)$ and $\beta_1,\ldots 
,\beta_m$ are nonnegative 
functions in $B(S)$.  Then the conclusions of Theorem 
\ref{mf} hold for 
\[Af=\beta_1A_1f+\cdots +\beta_mA_mf.\]

By (\ref{rpid}), $X$ and $Y$ are ``intertwined'' in the sense of 
\cite{RP81}.
\end{remark}

\begin{proof}
Theorem 3.2 of \cite{Kur98} can be extended to operators 
satisfying (\ref{opest}) by applying Corollary 1.12 of 
\cite{KS01} (with the operator $B$ in that corollary set 
equal zero) in place of Theorem 2.6 of \cite{Kur98}.  
Alternatively, see Corollary 3.2 of \cite{KN11}
\end{proof}

\subsection{Stochastic equations for processes built from 
bounded generators}\label{sumbnd}
We are primarily interested in generators of the form
\begin{equation}Af(x)=\int_{{\Bbb U}}(P_zf(x)-f(x))\mu (dz),\label{modgen}\end{equation}
where for each $z\in {\Bbb U}$, $P_z$ is a transition operator on a 
complete, separable metric space $E$, appropriately 
measurable as a function of $z\in {\Bbb U}$, and $\mu$ is a $\sigma$-finite 
measure on ${\Bbb U}$.  To illustrate the type of stochastic 
equation we have in mind, let 
\[A_0f(x)=\lambda_0\int_E(f(y)-f(x))\eta (x,dy),\]
where $0<\lambda_0<\infty$ and $\eta$ is a transition function on $
E$.  We 
can always find a probability measure $\nu_0$ on a 
measurable space ${\Bbb U}_0$  and a measurable function 
$H_0(x,u):E\times {\Bbb U}_0\rightarrow E$ 
satisfying $\eta (x,C)=\int_{{\Bbb U}_0}{\bf 1}_C(H_0(x,u))\nu (d
u)$, $C\in {\cal B}(E)$, so that 
\[\lambda_0\int_E(f(y)-f(x))\eta (x,dy)=\lambda_0\int_{{\Bbb U}_0}
(f(H_0(x,u))-f(x))\nu_0(du).\]
See, for example, the construction in \cite{BD83}.

If $N$ is a Poisson process with parameter $\lambda_0$, 
$U_0,U_1,\ldots$ are independent ${\Bbb U}_0$-valued random variables 
with distribution $\nu$$_0$, and $X(0)$ is an $E$-valued random 
variable, $N$, $\{U_i\}$, and $X(0)$ independent, then there is a 
unique, $E$-valued process $X$ satisfying
\begin{equation}f(X(t))=f(X(0))+\int_0^t(f(H_0(X(s-),U_{N(s-)}))-
f(X(s-))dN(s),\label{beq}\end{equation}
for all $f\in B(E)$, and $X$ will be a solution of the 
martingale problem for $A_0$.  Since in this case, $A_0$ is  a 
bounded operator and the martingale problem is 
well-posed, it follows that the martingale problem and 
the stochastic equation are equivalent in the sense that 
every solution of  the stochastic equation is a solution 
of the martingale problem and every solution of  the 
martingale problem is a weak solution of the stochastic 
equation.

In general, we are interested in situations where 
uniqueness is not necessarily known for either the 
martingale problem or the stochastic equation, but we 
still want to know that the two are equivalent.  We 
will obtain our result by application of  the Markov 
mapping theorem using arguments similar to those used 
in \cite{Kur11}.  Let us illustrate these arguments by 
proving what we already know regarding the martingale 
problem for  
$A_0$ and (\ref{beq}).

Let $\hat {B}_0$ be the generator for a process in $S=E\times {\Bbb U}_
0\times \{-1,1\}$ 
given by 
\[\hat {B}_0\hat {f}(x,u,\theta )=\lambda_0\int_{{\Bbb U}_0}(\hat {
f}(H_0(x,u),u',-\theta )-\hat {f}(x,u,\theta ))\nu_0(du'),\quad\hat {
f}\in B(S),\]
and setting 
\[f(x)=\frac 12\int_{{\Bbb U}_0}\hat {f}_0(x,u,1)\nu_0(du)+\frac 
12\int_{{\Bbb U}_0}\hat {f}_0(x,u,-1)\nu_0(du),\]
observe that 
\[A_0f(x)=\frac 12\int_{{\Bbb U}_0}\hat {B}_0\hat {f}(x,u,1)\nu_0
(du)+\frac 12\int_{{\Bbb U}_0}\hat {B}_0\hat {f}(x,u,-1)\nu_0(du)
.\]
The Markov mapping theorem implies that if $\hat {X}$ is a 
solution of the martingale problem for  $A_0$, there exists 
a solution $Z=(X,U,\Theta )$ of the martingale problem for $\hat {
B}_0$ such 
that $X$ has the same distribution as $\hat {X}$.

Let $N(t)$ be the counting process satisfying 
$\Theta (t)=\Theta (0)(-1)^{N(t)}$.
 Note that setting $\hat {f}(x,u,\theta )=\theta$, 
\[M_{\theta}(t)=\Theta (t)-\int_0^t\hat {B}_0f(Z(s))ds=\Theta (t)
+2\int_0^t\lambda_0\Theta (s)ds\]
is a martingale and
\[N(t)=-\frac 12\int_0^t\Theta (s-)d\Theta (s)=-\frac 12\int_0^t\Theta 
(s-)dM_{\theta}(s)+\lambda_0t.\]
Consequently, $N(t)-\lambda_0t$ is a martingale, and hence
 $N$ is a Poisson process with 
intensity $\lambda_0$.

\begin{lemma}
For any bounded function $f$ on $E$,
\begin{equation}f(X(t))=f(X(0))+\int_0^t(f(H_0(X(s-,U(s-)))-f(X(s
-)))dN(s).\label{steq1}\end{equation}
\end{lemma}

\begin{proof}
To see that this identity holds, let
\[M_f(t)=f(X(t))-f(X(0))-\int_0^t\hat {B}_0f(X(s),U(s),\Theta (s)
)ds.\]
We have the following Meyer processes (see Lemma 
5.1 of \cite{Kur11}).
\begin{eqnarray*}
\langle M_f\rangle_t&=&\int_0^t(\lambda_0(f^2(H_0(X(s),U(s))-f^2(
X(s))\\
&&\qquad\qquad\qquad -2f(X(s))\lambda_0(f(H_0(X(s),U(s)))-f(X(s))
)ds\\
&=&\int_0^t\lambda_0(f(H_0(X(s),U(s)))-f(X(s)))^2ds,\\
\langle M_f,M_{\theta}\rangle_t&=&\int_0^t\Big(\lambda_0(f(H_0(X(
s),U(s)))(-1)\Theta (s)-f(X(s))\Theta (s))\\
&&\qquad\qquad +2\lambda_0f(X(s))\Theta (s)-\Theta (s)\lambda_0(f
(H_0(X(s),U(s))-f(X(s)))\Big)ds\\
&=&-\int_0^t2\Theta (s)\lambda_0(f(H_0(X(s),U(s))-f(X(s)))ds,\\
\langle M_{\theta}\rangle_t&=&\int_0^t2\Theta (s)2\lambda_0(s)\Theta 
(s)=4\lambda_0t.\end{eqnarray*}
Then
\begin{eqnarray*}
M(t)&=&f(X(t))-f(X(0))-\int_0^t(f(H_0(X(s-),U(s-)))-f(X(s-)))dN(s
)\\
&=&M_f(t)+\frac 12\int_0^t(f(H_0(X(s-),U(s-)))-f(X(s-)))\Theta (s
-)dM_{\theta}(s)\end{eqnarray*}
is a martingale and
\begin{eqnarray*}
\langle M\rangle&=&\langle M_f\rangle +\int_0^t(f(H_0(X(s-),U(s-)
))-f(X(s-)))\Theta (s-)d\langle M_f,M_{\theta}\rangle\\
&&\qquad +\frac 14\int_0^t(f(H_0(X(s),U(s)))-f(X(s)))^2d\langle M_{
\theta}\rangle_s\\
&=&0,\end{eqnarray*}
so $M=0$ and (\ref{steq1}) holds.
\end{proof}

We now assume that $\mu$ is in (\ref{modgen}) is infinite, 
but $\sigma$-finite.  Writing 
${\Bbb U}=\cup_{k=1}^{\infty}{\Bbb U}_k$ as a disjoint union of sets of finite measure, 
we can write 
\begin{equation}Af(x)=\sum_{k=1}^{\infty}\int_{{\Bbb U}_k}(P_zf(x
)-f(x))\mu (dz)\equiv\sum_{k=1}^{\infty}B_kf(x),\label{Asum}\end{equation}
where each $B_k$ is a bounded generator, and hence can be 
written as 
\begin{equation}B_kf(x)=\lambda_k\int_E(f(y)-f(x))\eta_k(x,dy)=\lambda_
k\int_{{\Bbb U}_k}(f(H_k(x,u))-f(x))\nu_k(du),\label{bk}\end{equation}
for  $\lambda_k=\mu ({\Bbb U}_k)$, and some $H_k:E\times {\Bbb U}_
k\rightarrow E$, and $\nu_k\in {\cal P}({\Bbb U}_k)$. 
We are implicitly assuming that ${\Bbb U}_k$ is rich enough to 
support a measure $\nu_k$ for which the desired $H_k$ will 
exist.  One can always replace ${\Bbb U}$ by ${\Bbb U}\times [0,1
]$ and $\mu$ by 
$\mu\times\ell$.

To be specific, we will simply assume that $B_k$ is given 
by the right side of (\ref{bk}).
To make the definition of $A$ as the sum of the $B_k$ 
precise, let ${\cal D}\subset C_b(E)$, and
 assume the following conditions.

\begin{condition}\label{sumcnd2}

\begin{itemize}
\item[a)] ${\cal D}$ is  closed under multiplication and separates 
points in $E$.

\item[b)] For each $f\in {\cal D}$,
\[Af(x)\equiv\lim_{m\rightarrow\infty}\sum_{k=1}^mB_kf(x)\]
exists pointwise in $E$.  

\item[c)]\label{sumcndc}\ There exists $\psi\in M(E)$ such that $
\psi\geq 1$ 
and for each $f\in {\cal D}$, 
there exists $c_f$ and $m_f$ such that for $m\geq m_f$, 
\[|\sum_{k=m+1}^{\infty}B_kf(x)|\equiv |Af(x)-\sum_{k=1}^mB_kf(x)
|\leq c_f\psi (x),\quad x\in E.\]
\end{itemize}
\end{condition}

Let ${\Bbb E}_m=E\times {\Bbb U}_1\times\cdots\times {\Bbb U}_m\times 
\{-1,1\}^m$,
\[{\cal D}(\hat {A}_m)=\{\hat {f}(x,u,\theta )=f(x)\prod_{k=1}^mg_
k(u_k,\theta_k):f\in {\cal D},g_k\in C_b({\Bbb U}_k\times \{-1,1\}
),1\leq k\leq m\},\]
and define a generator $\hat {A}_m$ for a 
process in ${\Bbb E}_m$ by 
\begin{eqnarray*}
&&\hat {A}_m\hat {f}(x,u_1,\ldots ,u_m,\theta_1,\ldots ,\theta_m)\\
&&\qquad =\sum_{k=1}^m\lambda_k\int_{{\Bbb U}_k}(\hat {f}(H_k(x,u_
k),\eta_k(u|u'_k),\eta_k(\theta |-\theta_k))-\hat {f}(x,u,\theta 
))\nu_k(du'_k)\\
&&\qquad\qquad\qquad\qquad\qquad\qquad\qquad\qquad\qquad +\prod_{
k=1}^mg(u_k,\theta_k)\sum_{l=m+1}^{\infty}B_lf(x),\end{eqnarray*}
where for an arbitrary set $S$, 
for $z\in S^{\infty}$ and $z'_k\in S$, $\eta_k(z|z'_k)$ is the element of  
$S^{\infty}$ obtained from $z$  by replacing $z_k$ by $z_k'$.
If $\hat {X}$ is a solution of the martingale problem for 
$A$ satisfying 
\[E[\int_0^t\psi (\hat {X}(s))ds]<\infty ,\quad t\geq 0,\]
the Markov mapping theorem implies that for each $m$,
there exists a solution $(X^{(m)},U^{(m)},\Theta^{(m)})$ of the martingale 
problem for $\hat {A}_m$ such that $X^{(m)}$ and $\hat {X}$  have the same 
distribution.  By induction, the sequence of processes 
can be constructed so that the restriction of 
$(X^{(m+1)},U^{(m+1)},\Theta^{(m+1)})$ to ${\Bbb E}_m$ has the same distribution as  
$(X^{(m)},U^{(m)},\Theta^{(m)})$, and it follows that there exists a process 
$(X,{\Bbb U},\Theta )$ in ${\Bbb E}=E\times {\Bbb U}_1\times {\Bbb U}_
2\times\cdots\times \{-1,1\}^{\infty}$ 
so that the restriction of 
$(X,U,\Theta )$ to ${\Bbb E}_m$ has the same distribution as  
$(X^{(m)},U^{(m)},\Theta^{(m)})$.

Consequently,
\begin{eqnarray*}
\hat {M}_f^m(t)&=&f(X(t))-f(X(0))-\sum_{k=1}^m\int_0^t(f(H_k(X(s-
),U_k(s-)))-f(X(s-)))dN_k(s)\\
&&\qquad -\int_0^t\sum_{k\geq m+1}B_kf(X(s))ds\\
&=&f(X(t))-f(X(0))-\sum_{k=1}^m\int_0^t\lambda_k(f(H_k(X(s-),U_k(
s-)))-f(X(s-)))ds\\
&&\qquad -\int_0^t\sum_{k\geq m+1}B_kf(X(s))ds-\sum_{k=1}^m\int_0^
t(f(H_k(X(s-),U_k(s-)))-f(X(s-)))d\tilde {N}_k(s)\\
&=&M_f^m(t)+\sum_{k=1}^m\frac 12\int_0^t(f(H_k(X(s-),U_k(s-)))-f(
X(s-)))\Theta_k(s-)dM_{\theta_k}(s)\end{eqnarray*}
is a $\{{\cal F}_t^m\}$-martingale for 
${\cal F}_t^m=\sigma (\hat {X}(s),U_1(s),\ldots ,U_m(s),\Theta_1(
s),\ldots ,\Theta_m(s):s\leq t)$.

Note that
\begin{eqnarray*}
\langle M_f^m\rangle&=&\sum_{k=1}^m\int_0^t\lambda_k(f^2(H_k(X(s)
,U_k(s)))-f^2(X(s))\\
&&\qquad -2f(X(s))\lambda_k(f(H_k(X(s-),U_k(s-)))-f(X(s-))))ds\\
&&\qquad +\int_0^t\sum_{k\geq m+1}(B_kf^2(X(s))-2f(X(s))B_kf(X(s)
))ds,\\
\langle M_{\theta_k}\rangle_t&=&4\lambda_kt,\\
\langle M_f^m,M_{\theta_k}\rangle&=&\int_0^t\Big(\sum_{1\leq l\neq 
k\leq m}\Theta_k(s)\lambda_l(f(H_l(X(s),U_l(s)))-f(X(s))\\
&&\qquad -\lambda_k(\Theta_k(s)(f(H_k(X(s),U(s))+f(X(s)))\\
&&\qquad +\Theta_k(s)\sum_{l\geq m+1}B_lf(X(s))-\Theta_k(s)\sum_{
l\geq m+1}B_lf(X(s))\\
&&\qquad -\Theta_k(s)\sum_{l=1}^m\lambda_l(f(H_l(X(s),U_l(s)))-f(
X(s))\\
&&\qquad +2\lambda_k\Theta_k(s)f(X(s))\Big)\\
&=&-\int_0^t2\lambda_k\Theta_k(s)(f(H_k(X(s),U_k(s)))-f(X(s)))ds.\end{eqnarray*}
Consequently,
\begin{eqnarray*}
\langle\hat {M}^m_f\rangle_t&=&\langle M_f^m\rangle_t+\sum_{k=1}^
m\int_0^t(f(H_k(X(s),U_k(s)))-f(X(s)))\Theta_k(s)d\langle M_f^m,M_{
\theta_k}\rangle_s\\
&&\qquad +\sum_{k=1}^m\int_0^t\frac 14(f(H_k(X(s),U_k(s)))-f(X(s)
))^2d\langle M_{\theta_k}\rangle_s\\
&=&\int_0^t\Big(\sum_{k=1}^m\lambda_k(f^2(H_k(X(s),U_k(s)))-f^2(X
(s)))\\
&&\qquad\qquad -2f(X(s))\sum_{k=1}^m\lambda_k(f(H_k(X(s-),U_k(s-)
))-f(X(s-))))ds\\
&&\qquad\qquad +\sum_{k\geq m+1}(B_kf^2(X(s))-2f(X(s))B_kf(X(s)))\\
&&\qquad\qquad -2\sum_{k=1}^m\lambda_k(f(H_k(X(s),U_k(s)))-f(X(s)
))^2\\
&&\qquad\qquad +\sum_{k=1}^m\lambda_k(f(H_k(X(s),U_k(s)))-f(X(s))
)^2\Big)ds\\
&=&\int_0^t\sum_{k\geq m+1}(B_kf^2(X(s))-2f(X(s))B_kf(X(s))).\end{eqnarray*}

\begin{theorem}\label{sumeq}
Let $\{B_k\}$ be a sequence of  bounded generators of the 
form (\ref{bk}), and 
assume that
Condition \ref{sumcnd2} holds.  Suppose that $\hat {X}$  is a 
solution of the martingale problem for $A$ satisfying 
\[E[\int_0^t\psi (\hat {X}(s))ds]<\infty ,\quad t\geq 0.\]
Then, for each $f\in {\cal D}$, 
\[f(X(t))=f(X(0))+\sum_{k=1}^{\infty}\int_0^t(f(H_k(X(s-),U_k(s-)
))-f(X(s-)))dN_k(s),\]
in the sense that, for each  $T\geq 0$,
\[\lim_{m\rightarrow\infty}\sup_{t\leq T}|f(X(t))-f(X(0))-\sum_{k
=1}^m\int_0^t(f(H_k(X(s-),U_k(s-)))-f(X(s-)))dN_k(s)|=0\]
in probability.
\end{theorem}

\begin{proof}
Since $\langle\hat {M}_f^m\rangle_t\rightarrow 0$, it follows that $\sup_{
t\leq T}|\hat {M}_f^m(t)|\rightarrow 0$, and 
since
\begin{eqnarray*}
\hat {M}_f^m(t)&=&f(X(t))-f(X(0))-\sum_{k=1}^m\int_0^t(f(H_k(X(s-
),U_k(s-)))-f(X(s-)))dN_k(s)\\
&&\qquad -\int_0^t\sum_{k\geq m+1}B_kf(X(s))ds\end{eqnarray*}
and the last term goes to zero, the lemma 
follows.
\end{proof}